\documentclass[12pt]{article}
\usepackage{amsfonts}
\usepackage{}
\usepackage{graphicx}
\usepackage{mathrsfs}
\usepackage{booktabs}
\usepackage{multirow}
\usepackage{amsthm}
\usepackage{amsmath}
\usepackage{amsmath,amsthm,amsfonts,amssymb}
\usepackage{xcolor,graphicx,float}
\usepackage{mathrsfs}
\usepackage{caption}
\usepackage{subcaption}
\usepackage{cases}
\usepackage{latexsym,euscript,dsfont}
\usepackage{txfonts}
\usepackage{bbm}
\usepackage{colortbl}
\usepackage{cite}
\usepackage{enumerate}
\usepackage{xcolor}
\usepackage{tikz}
\usetikzlibrary{arrows,shapes,chains}
\allowdisplaybreaks[4]
\renewcommand\baselinestretch{1.2}

\input amssym.def
\input amssym

\title{{\LARGE\bf
               Open-loop and closed-loop solvabilities for discrete-time mean-field stochastic linear \\quadratic optimal control problems\thanks{This work was supported by the State Key Program of National Natural Science Foundation of China (Grant No. 12231008)}
               }
}
\author{ \large \ Teng Song\thanks{Corresponding author, E-mail: songt@whut.edu.cn.}, \  Bin Liu\thanks{E-mail: binliu@mail.hust.edu.cn.} \\
      $^{\dag}$ Department of Mathematics, Wuhan University of Technology \\
       $^{\ddagger}$ School of Mathematics and Statistics, Hubei Key Laboratory of Engineering\\ Modeling and Scientific Computing, Huazhong University of Science\\
        and Technology, Wuhan 430074, Hubei, P. R. China\\
       }
\date{}
\begin{document}
\maketitle
\par\noindent

\par\noindent
{\bf Abstract}
 \par

This paper discusses the discrete-time mean-field stochastic linear quadratic optimal control problems, whose weighting matrices in the cost functional are not assumed to be definite. The open-loop solvability is characterized by the existence of the solution to a mean-field forward-backward stochastic difference equations with a convexity condition and a stationary condition. The closed-loop solvability is presented by virtue of the existences of the regular solution to the generalized Riccati equations and the solution to the linear recursive equation, which is also shown by the uniform convexity of the cost functional. Moreover, based on a family of uniformly convex cost functionals, the finiteness of the problem is characterized. Also, it turns out that a minimizing sequence, whose convergence is equivalent to the open-loop solvability of the problem.
Finally, some examples are given to illustrate the theory developed.



\newpage
\section{Introduction}\indent
Let $(\Omega,\mathfrak{F},\mathbb{P})$ be a complete probability space on which a sequence $\mathfrak{F}$-measurable
$\mathbb{R}$-valued random variables $\{\omega_{k}:k=l,\ldots,N-1\}$ are defined such that $\mathfrak{F}_{k}\subset\mathfrak{F}$ is the $\sigma$-field generated by $\xi,\omega_{l},\omega_{l+1},\ldots,\omega_{k}$, i.e., $\mathfrak{F}_{k}=\sigma\{\xi,\omega_{l},\omega_{l+1},\ldots,\omega_{k}\}$, $\mathfrak{F}_{0}=\{\emptyset,\Omega\}$.
Consider the following discrete-time stochastic control system
\begin{eqnarray}
  \left\{
  \begin{array}{ll}
x_{k+1}=A_{k}x_{k}+\bar{A}_{k}\mathbb{E}x_{k}+B_{k}u_{k}+\bar{B}_{k}\mathbb{E}u_{k}+b_{k}\\
~~~~~~~~~~~~~~~~~+(C_{k}x_{k}+\bar{C}_{k}\mathbb{E}x_{k}+D_{k}u_{k}+\bar{D}_{k}\mathbb{E}u_{k}+\sigma_{k})\omega_{k},\\
x_{l}=\xi,~~k=l,\ldots,N-1,~~l\in\{0,1,\ldots,N-1\},
\end{array}\right.
\end{eqnarray}
and the corresponding cost functional
\begin{align}
J(l,\xi;u)&=\mathbb{E}\Bigg\{\langle Gx_{N},x_{N}\rangle+2\langle g,x_{N}\rangle+\langle \bar{G}\mathbb{E}x_{N},\mathbb{E}x_{N}\rangle+2\langle \bar{g},\mathbb{E}x_{N}\rangle\nonumber\\
&~~~~~~~~~+\sum_{k=l}^{N-1}\left[\left\langle  \left(
                            \begin{array}{cc}
                              Q_{k} & S_{k}' \\
                              S_{k} & R_{k}\\
                            \end{array}
                          \right)\left(
                                   \begin{array}{c}
                                     x_{k} \\
                                     u_{k} \\
                                   \end{array}
                                 \right),\left(
                                   \begin{array}{c}
                                     x_{k} \\
                                     u_{k} \\
                                   \end{array}
                                 \right)\right\rangle+2\left\langle  \left(
                                                            \begin{array}{c}
                                                              q_{k} \\
                                                              \rho_{k} \\
                                                            \end{array}
                                                          \right),\left(
                                   \begin{array}{c}
                                     x_{k} \\
                                     u_{k} \\
                                   \end{array}
                                 \right)
                                 \right\rangle\right]\nonumber\\
&~~~~~~~~~+\sum_{k=l}^{N-1}\left[\left\langle  \left(
                            \begin{array}{cc}
                              \bar{Q}_{k} & \bar{S}_{k}' \\
                              \bar{S}_{k} & \bar{R}_{k} \\
                            \end{array}
                          \right)\left(
                                   \begin{array}{c}
                                     \mathbb{E}x_{k} \\
                                     \mathbb{E}u_{k} \\
                                   \end{array}
                                 \right),\left(
                                   \begin{array}{c}
                                     \mathbb{E}x_{k} \\
                                     \mathbb{E}u_{k} \\
                                   \end{array}
                                 \right)\right\rangle+2\left\langle  \left(
                                                            \begin{array}{c}
                                                              \bar{q}_{k} \\
                                                              \bar{\rho}_{k} \\
                                                            \end{array}
                                                          \right),\left(
                                   \begin{array}{c}
                                     \mathbb{E}x_{k} \\
                                     \mathbb{E}u_{k} \\
                                   \end{array}
                                 \right)
                                 \right\rangle\right]
\Bigg\},
\end{align}
where $x_{k}\in \mathbb{R}^{n}$, $u_{k}\in \mathbb{R}^{m}$ are the state and control processes, respectively. $\mathbb{E}$ means the expectation operator, and
$A_{k}, \bar{A}_{k}, B_{k}, \bar{B}_{k}, C_{k}, \bar{C}_{k}, D_{k}, \bar{D}_{k}$, $Q_{k}, \bar{Q}_{k}, S_{k}, \bar{S}_{k}, R_{k}$, $\bar{R}_{k}$, $G, \bar{G}$ are deterministic matrix-valued functions of proper dimensions with $Q_{k}'=Q_{k}$, $\bar{Q}_{k}'=\bar{Q}_{k}$, $R_{k}'=R_{k}$, $\bar{S}_{k}'=\bar{S}_{k}$, $G=G', \bar{G}=\bar{G}'$; $b_{k}, \sigma_{k}$ are vector-valued functions, $g, \bar{g}$ are random variables, and $q, \bar{q}, \rho, \bar{\rho}$ are vector-valued random processes. $\{\omega_{k}\}$ means the martingale difference sequence in the sense that
$\mathbb{E}[\omega_{k+1}|\mathfrak{F}_{k}]=0,~\mathbb{E}[\omega_{k+1}^{2}|\mathfrak{F}_{k}]=1.$ $\xi$ is a square integrable random variable, and independent of $\{\omega_{k}\}$.
Indeed, the system (1.1) is a discrete-time stochastic difference equation (SDE) of Mckean-Vlasov type, also called a mean-field SDE, which could be applied to investigate particle systems at the mesoscopic level.
Then we could propose the stochastic optimal control problem as follows.

$\mathbf{Problem~(MF\mbox{-}LQ).}$ For any given initial pair $(l,\xi)$, find a $u^{\ast}\in \mathscr{U}_{ad}$ (which would be defined in Section 2) such that
\begin{eqnarray}
J(l,\xi;u^{\ast})=\inf_{u\in\mathscr{U}_{ad}}J(l,\xi;u)\doteq V(l,\xi).
\end{eqnarray}
It is widely shared that any $u^{\ast}$ satisfying (1.3) is called an optimal control of Problem (MF-LQ) at the initial pair $(l,\xi)$, and the corresponding $x^{\ast}$, $(x^{\ast},u^{\ast})$ are called an optimal state process and an optimal pair, respectively. $V(l,\xi)$ is named the value function of Problem (MF-LQ). In the special case when $b_{k}, \sigma_{k}, q_{k}, \bar{q}_{k}, \rho_{k}, \bar{\rho}_{k}, g, \bar{g}=0$, we denote the corresponding Problem (MF-LQ) by Problem (MF-LQ)$^{0}$. Likewise, the corresponding cost functional and value function are denoted by $J^{0}(l,\xi;u)$ and $V^{0}(l,\xi)$, respectively. We would like to mention that our results also remain for discrete-time stochastic systems with multiplicative noise.

Linear quadratic (LQ) optimal control problems can be traced back to the work of Kalman \cite{K60} in 1960, and is now a classical as well as essential problem in modern control theory. Stochastic LQ optimal control has attracted wide attention and effort thereafter \cite{CLZ98,D15,LLY20,LZ01,NCZZ17,SLY16,T03,SY20,WSY19}. Recently, more and more researchers have been interested in mean-field LQ problems. The mean-field theory is developed to discuss the collective behavior caused by individuals' mutual interactions in various sociological and physical dynamical systems, and it has been applied in many fields of finance, engineering, economics and games, in particular, in optimal control field; see for example, \cite{B23,BSYY16,S21,SY14,DNW22,LLL20,NLZ16b,LSY16}. The continuous-time mean-field LQ stochastic optimal control problem was initially carried out by Yong \cite{Y13} and then extended by Huang, Li and Yong \cite{HLY15}, Yong \cite{Y17}, Sun \cite{S17}, Lin, Jiang and Zhang \cite{LJZ18}, Li, Sun and Xiong \cite{LSX19}, L\"{u} \cite{L19b}, Sun, Xiong and Yong \cite{SXY21}, and references therein. To some extent, the study of mean-field stochastic LQ problem for continuous-time framework are relatively mature.

Discrete-time optimal control problems are more pertinent to economic, engineering, biomedical, operation research, optimizing complex technological problems and so on. It is more natural and sufficient to describe a system by a discrete-time model. For the discrete-time mean-field LQ stochastic optimal control problems, one of the main motivations comes from the mean-field games \cite{BDT19,L19a,SWW21,M20,S23}. Roughly speaking, mean-field games adopt decentralized controls, to put it differently, the controls are chose to accomplish each individual's own purpose by means of local information. By contrast, mean-field games could be viewed as a standard control problem and an equilibrium, while the mean-field LQ optimal control problem is a nonstandard control problem.
The second motivation is from the Markowitz's mean-variance portfolio selection problems in financial applications \cite{BC20}.
Thirdly, one attempts to retain the variances $var(x)$ and $var(u)$ small so that the control and state processes are not too sensitive to the possible variation of the random events. Naturally, one needs to consider the following cost functional
\begin{align}
\bar{J}(l,\xi;u)=\mathbb{E}\Bigg\{\langle Gx_{N},x_{N}\rangle+\tilde{G}var(x_{N})\!+\!\sum_{k=l}^{N-1}\Big[\langle Q_{k}x_{k},x_{k}\rangle
+\tilde{Q}_{k}var(x_{k})+\langle R_{k}u_{k},u_{k}\rangle+\tilde{R}_{k}var(u_{k})\Big]\Bigg\}.
\end{align}
Since $var(x_{k})=\mathbb{E}|x_{k}|^{2}-(\mathbb{E}x_{k})^{2}$, $var(u_{k})=\mathbb{E}|u_{k}|^{2}-(\mathbb{E}u_{k})^{2}$, then the optimal control problem with cost functional (1.4) can be regarded as a mean-field LQ optimal control problem. For more motivations to discuss this type of questions, see \cite{BFY13,G16,SL23} for details. Recently,
by the kernel-range decomposition and pseudo-inverse, Elliott, Li and Ni \cite{ELN13} discussed the mean-field LQ stochastic optimal control problem over finite horizon and derived the necessary and sufficient solvability conditions. Later, the same authors in \cite{NEL15} generalized the results to the infinite horizon case.
Most recently, Zhang, Qi and Fu \cite{ZQF18} concluded that the mean-field system is $L^{2}$-stabilizable if and only if two coupled generalized algebraic Riccati equations have a unique positive definite (semi-definite) solution under the exact observability (exact detectability) assumption.
Unfortunately, in the above works, the following standard condition was taken for granted
\begin{eqnarray}
  \left\{
  \begin{array}{ll}
G\geq 0,~R_{k}>0,~Q_{k}-S_{k}'R_{k}^{-1}S_{k}\geq0,\\
G+\bar{G}\geq 0,~R_{k}+\bar{R}_{k}>0,~(S_{k}+\bar{S}_{k})'(R_{k}+\bar{R}_{k})^{-1}(S_{k}+\bar{S}_{k})\geq0.
\end{array}\right.
\end{eqnarray}
In that case, Problem (MF-LQ) exists a unique optimal control and the corresponding Riccati difference equation is uniquely solvable. Interestingly, Chen, Li and Zhou \cite{CLZ98} clarified that Problem (MF-LQ) may still be well-posed even though $R_{k}$ is not positive semidefinite. Hence, it is more worthwhile and constructive to explore Problem (MF-LQ) without (1.5). Ni, Zhang and Li \cite{NZL15} found that the well-posedness and the solvability of the mean-field LQ problem are both equivalent to the solvability of two coupled difference Riccatic equations and a linear recursive equation. Subsequently, the same authors \cite{NLZ16a} studied the finite and infinite horizon indefinite mean-field LQ control problems. In the past decade, some new and interesting results were obtained, among them, we refer the readers to \cite{BDT20,LLL18,MZZ19,ZDS23,WZL18,GLZ23}. Nevertheless, there are few works that deal with the open-loop and closed-loop solvabilities problems for discrete-time mean-field systems.

On the other hand, it can be asserted that most of the techniques applied in continuous-time case cannot be directly adopted to the discrete-time framework. For instance, It\^{o} formula, which is the basic method in continuous-time stochastic systems, fails in our setting. In addition, the classical results on continuous-time systems cannot be easily extended to discrete-time case. So we must propose a subtle and characteristic techniques to investigate the discrete-time stochastic system.

In this paper, we shall investigate the open-loop and closed-loop solvabilities for discrete-time mean-field stochastic LQ optimal control problems. We have two main goals in this paper: (i) establish a theory for Problem (MF-LQ) parallel to that of continuous-time case; (ii) discuss the differences between the closed-loop and open-loop solvabilities. Firstly, we characterize the open-loop solvabilty by the mean-field forward-backward stochastic difference equations (MF-FBSDEs) with a convexity condition and a stationary condition. Then, the closed-loop solvability of Problem (MF-LQ) is obtained by the existences of a regular solution to the generalized Riccati equations (GRE) and a solution to the linear recursive equation (LRE). Besides, we also build the equivalence between the uniform convexity of the cost functional and the strongly regular solvability of the GRE. The finiteness of Problem (MF-LQ) is portrayed via the convergence of the solutions to a family of GRE. At last but not least, by a minimizing sequence, we present the open-loop solvability of the problem (MF-LQ) again.


The rest of the paper is organized as follows. Section 2 presents some preliminary results. In Section 3, we are devoted to discussing the representation of the cost functional. Section 4 characterizes the open-loop and closed-loop solvabilities of Problem (MF-LQ), as well as the uniform convexity of the cost functional. In Section 5, we consider the finiteness of Problem (MF-LQ) and the convexity of cost functional. Section 6 studies the open-loop solvability of Problem (MF-LQ) based on the convergence of minimizing sequences. Two examples are shown in Section 7 to illustrate the theory results in previous sections. Finally, some perspectives and open problems conclude the paper.

\setcounter{equation}{0}
\section{Preliminaries}
Let $\mathbb{N}=\{l,l+1,\ldots,N-1\}$, $\mathbb{N}_{0}=\{0,1,\ldots,N-1\}$, $\bar{\mathbb{N}}=\{l,l+1,\ldots,N\}$, and $\mathbb{R}^{n\times m}$ be the Euclidean space of all $n\times m$ real matrices, endowed with the inner product $\langle A,B\rangle= Tr(A'B)$, in which $A'$ and $Tr(A)$ represent the transpose and trace of $A$, respectively. $Ker(A)$ and $\mathcal{R}(A)$ stand for the kernel space and the range of $A$, respectively.
We denote the inner space in possibly different Hilbert spaces by $\langle \cdot,\cdot\rangle$ when there is no confusion, let $L^{2}(\mathbb{N};\mathbb{H})$ be the space of all $\mathbb{H}$-valued functions that are square integrable on $\mathbb{N}$ and
\begin{align*}
&L_{\mathfrak{F}}^{2}(k;\mathbb{R}^{p})=\Big\{\xi~|~\xi\in\mathbb{R}^{p},~\xi~\mbox{is}~\mathfrak{F}_{k}\mbox{-measurable},
~\mathbb{E}|\xi|^{2}<\infty\Big\},\\
&L_{\mathfrak{F}}^{2}(\mathbb{N};\mathbb{R}^{p})=\Bigg\{\phi=(\phi_{l},\ldots,\phi_{N-1})~|~\phi_{k}\in\mathbb{R}^{p},~k\in \mathbb{N},~\phi_{k}~\mbox{is}~\mathfrak{F}_{k}\mbox{-measurable},
~\mathbb{E}\sum\limits_{k=l}^{N-1}|\phi_{k}|^{2}<\infty\Bigg\}.
\end{align*}
Also, we introduce the following admissible control set
$$\mathscr{U}_{ad}=\Bigg\{(u_{l},u_{l+1},\ldots,u_{N-1})~|~u_{k}~\mbox{is}~\mathfrak{F}_{k}\mbox{-measurable},~u_{k}\in \mathbb{R}^{m},~\mathbb{E}\sum\limits_{k=l}^{N-1}|u_{k}|^{2}<\infty\Bigg\}.$$
Now, we recall some important notions of LQ control problems.

\vskip3mm\noindent
\textbf{Definition 2.1.} Problem (MF-LQ) is said to be finite at initial pair $(l,\xi)\in\mathbb{N}_{0}\times \mathbb{R}^{n}$ if
\begin{align}
V(l,\xi)>-\infty.
\end{align}
\noindent
(i) Problem (MF-LQ) is said to be finite at $l\in\mathbb{N}_{0}$ if (2.1) holds for all $\xi\in \mathbb{R}^{n}$.

\noindent
(ii) Problem (MF-LQ) is said to be finite if (2.1) holds for all $(l,\xi)\in\mathbb{N}_{0}\times \mathbb{R}^{n}$.

\vskip3mm\noindent
\textbf{Definition 2.2.} $u^{\ast}$ is called to be an open-loop optimal control of Problem (MF-LQ) for the initial pair $(l,\xi)\in\mathbb{N}_{0}\times \mathbb{R}^{n}$ if
\begin{align}
J(l,\xi;u^{\ast})\leq J(l,\xi;u),~~\forall~u\in \mathscr{U}_{ad}.
\end{align}
(i) Problem (MF-LQ) is said to be (uniquely) open-loop solvable at $(l,\xi)$ if an open-loop optimal
\indent control (uniquely) exists at $l\in\mathbb{N}_{0}$.

\noindent
(ii) Problem (MF-LQ) is said to be (uniquely) open-loop solvable at $l$ if for given $l$, (2.2) holds for
\indent all $\xi\in \mathbb{R}^{n}$.

\noindent
(iii) Problem (MF-LQ) is said to be (uniquely) open-loop solvable if it is (uniquely) open-loop
\indent solvable at all $(l,\xi)\in\mathbb{N}_{0}\times \mathbb{R}^{n}$.

\vskip3mm\noindent
\textbf{Definition 2.3.}

\noindent
(i) Let $u=\{\Theta_{k}x_{k}+\bar{\Theta}_{k}\mathbb{E}x_{k}+v_{k},~k\in \mathbb{N}\}$ be a control of (1.1) with $\Theta_{k}, \bar{\Theta}_{k}\in \mathbb{R}^{m\times n}$, $v_{k}\in \mathbb{R}^{m\times m}$ being
\indent deterministic matrices. If $\Theta_{k}, \bar{\Theta}_{k}, v_{k}$ are independent of all the initial pairs $(l,\xi)$, then the triple
\indent $(\Theta, \bar{\Theta}, v)=\{(\Theta_{k}, \bar{\Theta}_{k}, v_{k}),~k\in \mathbb{N}\}$ is named a closed-loop strategy of (1.1).

\noindent
(ii) A closed-loop strategy $(\Theta, \bar{\Theta}, v)=\{(\Theta_{k}, \bar{\Theta}_{k}, v_{k}),~k\in \mathbb{N}\}$ is called a closed-loop optimal strategy \indent of Problem (MF-LQ) on $\mathbb{N}$ if
\begin{align}
J(l,\xi;\Theta_{k}x_{k}^{\ast}+\bar{\Theta}_{k}\mathbb{E}x_{k}^{\ast}+v_{k})\leq J(l,\xi;u),~~\forall~u\in \mathscr{U}_{ad},
\end{align}
where
\begin{eqnarray}
  \left\{
  \begin{array}{ll}
x_{k+1}^{\ast}=(A_{k}+B_{k}\Theta_{k})x_{k}^{\ast}+[\bar{A}_{k}+B_{k}\bar{\Theta}_{k}+\bar{B}_{k}(\Theta_{k}+\bar{\Theta}_{k})]\mathbb{E}x_{k}^{\ast}+B_{k}v_{k}+\bar{B}_{k}\mathbb{E}v_{k}+b_{k}\\
~~~~~~~~~~~~+\Big\{(C_{k}+D_{k}\Theta_{k})x_{k}^{\ast}+[\bar{C}_{k}+D_{k}\bar{\Theta}_{k}+\bar{D}_{k}(\Theta_{k}+\bar{\Theta}_{k})]\mathbb{E}x_{k}^{\ast}
+D_{k}v_{k}+\bar{D}_{k}\mathbb{E}v_{k}+\sigma_{k}\Big\}\omega_{k},\\
x_{l}^{\ast}=\xi.
\end{array}\right.
\end{eqnarray}
Problem (MF-LQ) is called to be (uniquely) closed-loop solvable at $\mathbb{N}$ if a closed-loop optimal strategy (uniquely) exists at $\mathbb{N}$;
Problem (MF-LQ) is called to be (uniquely) closed-loop solvable if it is (uniquely) closed-loop solvable at any $\mathbb{N}$.

\vskip2mm
With the above definitions, we proceed to make the following assumptions.

\noindent
${\rm{(\mathbf{A1})}}$ $J^{0}(l,0;u)\geq0,~\forall~u\in \mathscr{U}_{ad}$.

\noindent
${\rm{(\mathbf{A2})}}$ There exists a constant $\alpha>0$ such that $J^{0}(l,0;u)\geq\alpha\mathbb{E}\sum_{k=l}^{N-1}|u_{k}|^{2}$.

\noindent
${\rm{(\mathbf{A3})}}$ There exists a constant $\alpha>0$ such that $J(0,0;u)\geq\alpha|u_{k}|^{2}$.

\vskip3mm\noindent
\textbf{Lemma 2.1.} \cite{ACZ02} Let $F=F'$, $H=H'$ and $G$ be given deterministic matrices, $\gamma, \rho$ be given deterministic vectors with appropriate size. Consider the following quadratic form
$$f(x,u)=\mathbb{E}[x'Fx+2x'Gu+u'Hu+2\gamma'x+2\rho'u],$$
where $x, u$ are square integrable random variables defined on a probability space.
Then the following statements are equivalent:

\noindent
(i) For any $x$$, \inf_{u}f(x,u)>-\infty$.

\noindent
(ii) $H\geq0$, $Ker(H)\subseteq Ker(G)$, $\rho \in \mathcal{R}(H)$.

\noindent
(iii) $H\geq0$, $G=GHH^{\dagger}$, $\rho'=\rho'HH^{\dagger}$.

\vskip3mm\noindent
\textbf{Lemma 2.2.} \cite{NZL15} Let $H=H'$, $\bar{H}=\bar{H}'$, $G, \bar{G}$ be given deterministic matrices, and $a$ be given deterministic vector with appropriate size. Consider the following quadratic form
$$\bar{f}(x,u)=\mathbb{E}[2(x-\mathbb{E}x)'G(u-\mathbb{E}u)+(u-\mathbb{E}u)'H(u-\mathbb{E}u)+2(\mathbb{E}x)'\bar{G}\mathbb{E}u+(\mathbb{E}u)'\bar{H}\mathbb{E}u
+2a'\mathbb{E}u],$$
where $x, u$ are square integrable random variables defined on a probability space, then the following statements are equivalent:

\noindent
(i) For any $x$$, \inf_{u}\bar{f}(x,u)>-\infty$.

\noindent
(ii) $H\geq0$, $Ker(H)\subseteq Ker(G)$, $\bar{H}\geq0$, $Ker(\bar{H})\subseteq Ker(\bar{G})$, $a \in \mathcal{R}(H)$.

\noindent
(iii) $H\geq0$, $G=GHH^{\dagger}$, $\bar{H}\geq0$, $\bar{G}=\bar{G}\bar{H}\bar{H}^{\dagger}$, $a'=a'HH^{\dagger}$.

\vskip3mm\noindent
\textbf{Lemma 2.3.} For any $u\in \mathscr{U}_{ad}$, let $x^{(u)}$ be the solution of
\begin{eqnarray}
  \left\{
  \begin{array}{ll}
x_{k+1}^{(u)}=A_{k}x_{k}^{(u)}+\bar{A}_{k}\mathbb{E}x_{k}^{(u)}+B_{k}u_{k}+\bar{B}_{k}\mathbb{E}u_{k}+(C_{k}x_{k}^{(u)}+\bar{C}_{k}\mathbb{E}x_{k}^{(u)}+D_{k}u_{k}+\bar{D}_{k}\mathbb{E}u_{k})\omega_{k},\\
x_{l}^{(u)}=0.
\end{array}\right.
\end{eqnarray}
Then for any $\Theta_{k},~\bar{\Theta}_{k}\in L^{2}(\mathbb{N};\mathbb{R}^{m\times n})$, there is a constant $\mu>0$ such that
\begin{eqnarray}
  \left.
  \begin{array}{ll}
\mathbb{E}\sum\limits_{k=l}^{N-1}\left|u_{k}-\Theta_{k}\left(x^{(u)}_{k}-\mathbb{E}x^{(u)}_{k}\right)\right|^{2}\geq \mu\mathbb{E}\sum\limits_{k=l}^{N-1}|u_{k}|^{2},~~\forall~u_{k}\in \mathscr{U}_{ad},\\
~~~~~~~~~~~~~\mathbb{E}\sum\limits_{k=l}^{N-1}\left|u_{k}-\bar{\Theta}_{k}\mathbb{E}x^{(u)}_{k}\right|^{2}\geq \mu\mathbb{E}\sum\limits_{k=l}^{N-1}|u_{k}|^{2},~~\forall~u_{k}\in \mathscr{U}_{ad}.
\end{array}\right.
\end{eqnarray}
\begin{proof} Define a bounded linear operator $\mathcal{A}:\mathscr{U}_{ad}\rightarrow \mathscr{U}_{ad}$ by
$\mathcal{A}u=u-\Theta(x^{(u)}-\mathbb{E}x^{(u)}),$
then $\mathcal{A}$ is bijective and its inverse $\mathcal{A}^{-1}$ is given as
$$\mathcal{A}^{-1}u=u+\Theta (\tilde{x}^{(u)}-\mathbb{E}\tilde{x}^{(u)}),$$
where $\tilde{x}^{(u)}$ is the solution of
\begin{eqnarray}
  \left\{
  \begin{array}{ll}
\tilde{x}^{(u)}_{k+1}=(A_{k}+B_{k}\Theta_{k})\tilde{x}^{(u)}_{k}+B_{k}u_{k}+(\bar{A}_{k}-B_{k}\bar{\Theta}_{k})
\mathbb{E}\tilde{x}_{k}^{(u)}+\bar{B}_{k}\mathbb{E}u_{k}\nonumber\\
~~~~~~~~~~+\Big\{(C_{k}+D_{k}\Theta_{k})\tilde{x}^{(u)}_{k}+D_{k}u_{k}+(\bar{C}_{k}-D_{k}\bar{\Theta}_{k})\mathbb{E}\tilde{x}_{k}^{(u)}
+\bar{D}_{k}\mathbb{E}u_{k}\Big\}\omega_{k},\\
\tilde{x}^{(u)}_{l}=0.
\end{array}\right.
\end{eqnarray}
By virtue of the bounded inverse theorem, $\mathcal{A}^{-1}$ is bounded with norm $\|\mathcal{A}^{-1}\|>0$.
Therefore,
\begin{align*}
\mathbb{E}\sum\limits_{k=l}^{N-1}|u_{k}|^{2}=\mathbb{E}\sum\limits_{k=l}^{N-1}\left\|\mathcal{A}^{-1}\mathcal{A}u_{k}\right\|^{2}\leq \|\mathcal{A}^{-1}\|\mathbb{E}\sum\limits_{k=l}^{N-1}\|\mathcal{A}u_{k}\|^{2}
=\|\mathcal{A}^{-1}\|\mathbb{E}\sum\limits_{k=l}^{N-1}\left|u_{k}-\Theta_{k}\left(x^{(u)}_{k}-\mathbb{E}x^{(u)}_{k}\right)\right|^{2},~~u_{k}\in \mathscr{U}_{ad},
\end{align*}
which implies the first inequality of (2.6). Next, we shall show the second one. For any $v_{k}\in \mathscr{U}_{ad}$, let $z^{(v)}$ be the solution to the following difference equation
\begin{eqnarray}
  \left\{
  \begin{array}{ll}
z_{k+1}^{(v)}=(A_{k}+\bar{A}_{k})z_{k}^{(v)}+(B_{k}+\bar{B}_{k})v_{k},\\
z_{l}^{(v)}=0.
\end{array}\right.
\end{eqnarray}
For any $\bar{\Theta}_{k}\in L^{2}(\mathbb{N};\mathbb{R}^{m\times n})$, we define a bounded linear operator $\mathcal{B}:\mathscr{U}_{ad}\rightarrow \mathscr{U}_{ad}$ by
$\mathcal{B}v=v-\bar{\Theta} z^{(v)}.$
Similarly, $\mathcal{B}$ is invertible and
\begin{align*}
\mathbb{E}\sum\limits_{k=l}^{N-1}|v_{k}-\bar{\Theta}_{k}z_{k}^{(v)}|^{2}\geq \frac{1}{\mathcal{B}^{-1}}\mathbb{E}\sum\limits_{k=l}^{N-1}|v_{k}|^{2},~~v_{k} \in \mathscr{U}_{ad}.
\end{align*}
The proof is completed since $z_{k}^{(v)}$ satisfies (2.5) with $z_{k}^{(v)}=\mathbb{E}x^{(u)}_{k}$.
\end{proof}

\setcounter{equation}{0}
\section{Representation of the cost functional}
In this section, we shall discuss the cost functional from a Hilbert space aspect, by which some basic conditions will be obtained for finiteness, open-loop and closed-loop solvabilities of Problem (MF-LQ).

\vskip3mm\noindent
\textbf{Lemma 3.1.} Let $l\in \mathbb{N}_{0}$ be given, for any $\xi\in \mathbb{R}^{n}$, $\lambda\in \mathbb{R}$, and $u_{k}$, $v_{k}\in \mathscr{U}_{ad}$, then we have
\begin{align}
J(l,\xi;u_{k}+\lambda v_{k})&=J(l,\xi;u_{k})+\lambda^{2}J^{0}(l,0;v_{k})+2\lambda \mathbb{E}\sum\limits_{k=l}^{N-1}\langle B'_{k}\mathbb{E}(y_{k+1}|\mathfrak{F}_{k-1})+\bar{B}'_{k}\mathbb{E}y_{k+1}+\rho_{k}+\bar{\rho}_{k}\nonumber\\
&~~~~~~+D'_{k}\mathbb{E}(y_{k+1}\omega_{k}|\mathfrak{F}_{k-1})+\bar{D}'_{k}\mathbb{E}(y_{k+1}\omega_{k})
+S'_{k}x_{k}+\bar{S}'_{k}\mathbb{E}x_{k}+R_{k}u_{k}+\bar{R}_{k}\mathbb{E}u_{k},v_{k}\rangle,
\end{align}
where $(x_{k},y_{k})$ is the solution of the following MF-FBSDEs:
\begin{eqnarray}
  \left\{
  \begin{array}{ll}
x_{k+1}=A_{k}x_{k}+\bar{A}_{k}\mathbb{E}x_{k}+B_{k}u_{k}+\bar{B}_{k}\mathbb{E}u_{k}+b_{k}
+(C_{k}x_{k}+\bar{C}_{k}\mathbb{E}x_{k}+D_{k}u_{k}+\bar{D}_{k}\mathbb{E}u_{k}+\sigma_{k})\omega_{k},\\
y_{k}=A'_{k}\mathbb{E}(y_{k+1}|\mathfrak{F}_{k-1})+\bar{A}'_{k}\mathbb{E}y_{k+1}
+C'_{k}\mathbb{E}(y_{k+1}\omega_{k}|\mathfrak{F}_{k-1})+\bar{C}'_{k}\mathbb{E}(y_{k+1}\omega_{k})\\
~~~~~~~~~~~~~~~~~+Q_{k}x_{k}+\bar{Q}_{k}\mathbb{E}x_{k}+S_{k}u_{k}+\bar{S}_{k}\mathbb{E}u_{k}+q_{k}+\bar{q}_{k},\\
x_{l}=\xi,~~y_{N}=Gx_{N}+\bar{G}\mathbb{E}x_{N}+g+\bar{g}.
\end{array}\right.
\end{eqnarray}
Therefore, the map $u\mapsto J(l,\xi;u)$ is Gateaux differentiable with the Gateaux derivative given as
\begin{align*}
dJ(l,\xi;u_{k})&=2[B'_{k}\mathbb{E}(y_{k+1}|\mathfrak{F}_{k-1})+\bar{B}'_{k}\mathbb{E}y_{k+1}+\rho_{k}+\bar{\rho}_{k}+D'_{k}\mathbb{E}(y_{k+1}\omega_{k}|\mathfrak{F}_{k-1})\nonumber\\
&~~~~~~~~+\bar{D}'_{k}\mathbb{E}(y_{k+1}\omega_{k})
+S'_{k}x_{k}+\bar{S}'_{k}\mathbb{E}x_{k}+R_{k}u_{k}+\bar{R}_{k}\mathbb{E}u_{k}],
\end{align*}
and (3.1) could also be represented as
\begin{align}
J(l,\xi;u_{k}+\lambda v_{k})&=J(l,\xi;u_{k})+\lambda^{2}J^{0}(l,0;v_{k})+\lambda \mathbb{E}\sum\limits_{k=l}^{N-1}\langle dJ(l,\xi;u_{k}),v_{k}\rangle.
\end{align}
\begin{proof} If no confusion is likely, we shall use the following notations just in this proof. For $L^{2}_{\mathfrak{F}}(\mathbb{N};\mathbb{R}^{p})$ $(p=m, n)$, we introduce an inner product
\begin{align*}
\Big( z^{(1)},z^{(2)}\Big)=\mathbb{E}\sum\limits_{k=l}^{N-1}\left\langle z^{(1)},z^{(2)}\right\rangle,~~~z^{(1)},~z^{(2)}\in L^{2}_{\mathfrak{F}}(\mathbb{N};\mathbb{R}^{p}),
\end{align*}
and use the convention
$(Tx)(\cdot)=T.x.,~(\mathbb{E}v)(\cdot)=\mathbb{E}v$.
For $L^{2}_{\mathfrak{F}}(N-1;\mathbb{R}^{p})$ $(p=m, n)$, an inner product is given as
\begin{align*}
\Big(z^{(1)},z^{(2)}\Big)=\mathbb{E}\left\langle z^{(1)},z^{(2)}\right\rangle,~~~~z^{(1)},~z^{(2)}\in L^{2}_{\mathfrak{F}}(N-1;\mathbb{R}^{p}).
\end{align*}
The cost functional (1.2) is then rewritten as
\begin{align}
J(l,\xi;u)&=(Qx,x)+(\bar{Q}\mathbb{E}x,\mathbb{E}x)+(Ru,u)+(\bar{R}\mathbb{E}u,\mathbb{E}u)
+(Gx_{N},x_{N})+(\bar{G}\mathbb{E}x_{N},\mathbb{E}x_{N})+2(S'x,u)\nonumber\\
&~~~~~~~+2(\bar{S}'\mathbb{E}x,\mathbb{E}u)+2(q,x)+2(\rho,u)+2(\bar{q},\mathbb{E}x)+2(\bar{\rho},\mathbb{E}u)+2(g,x_{N})+2(\bar{g},\mathbb{E}x_{N}).
\end{align}
Denote by $x^{\lambda}$ the solution to (1.1) with control $u+\lambda v$, then
\begin{eqnarray*}
  \left\{
  \begin{array}{ll}
\frac{x_{k+1}^{\lambda}-x_{k+1}}{\lambda}=A_{k}\frac{x_{k}^{\lambda}-x_{k}}{\lambda}+\bar{A}_{k}\frac{\mathbb{E}x_{k}^{\lambda}-\mathbb{E}x_{k}}{\lambda}
+B_{k}v_{k}+\bar{B}_{k}\mathbb{E}v_{k}+b_{k}\\
~~~~~~~~~~~~~~~~~+\Big(C_{k}\frac{x_{k}^{\lambda}-x_{k}}{\lambda}+\bar{C}_{k}\frac{\mathbb{E}x_{k}^{\lambda}-\mathbb{E}x_{k}}{\lambda}
+D_{k}v_{k}+\bar{D}_{k}\mathbb{E}v_{k}+\sigma_{k}\Big)\omega_{k},\\
\frac{x_{l}^{\lambda}-x_{l}}{\lambda}=0,~~k\in \mathbb{N},~~l\in\mathbb{N}_{0}.
\end{array}\right.
\end{eqnarray*}
Set $z_{k}=\frac{x_{k}^{\lambda}-x_{k}}{\lambda}$, which is independent of $u$ and $\lambda$, then we get
\begin{eqnarray}
  \left\{
  \begin{array}{ll}
z_{k+1}=A_{k}z_{k}+\bar{A}_{k}\mathbb{E}z_{k}+B_{k}v_{k}+\bar{B}_{k}\mathbb{E}v_{k}+b_{k}
+(C_{k}z_{k}+\bar{C}_{k}\mathbb{E}z_{k}+D_{k}v_{k}+\bar{D}_{k}\mathbb{E}v_{k}+\sigma_{k})\omega_{k},\\
z_{l}=0,~~k\in \mathbb{N},~~l\in \mathbb{N}_{0}.
\end{array}\right.
\end{eqnarray}
By (3.4) and some calculations, we deduce that
\begin{align}
J(l,\xi;u+\lambda v)-J(l,\xi;u)&=2\lambda(Qx,z)+2\lambda(\bar{Q}\mathbb{E}x,\mathbb{E}z)+2\lambda(Ru,v)+2\lambda( \bar{R}\mathbb{E}u,\mathbb{E}v)+2\lambda(S'x,v)\nonumber\\
&~~~~~+2\lambda(Gx_{N},z_{N})+2\lambda(\bar{G}\mathbb{E}x_{N},\mathbb{E}z_{N})+2\lambda(\bar{S}'\mathbb{E}x,\mathbb{E}v)+2\lambda(S'z,u)\nonumber\\
&~~~~~+2\lambda(\bar{S}'\mathbb{E}z,\mathbb{E}u)+\lambda^{2}(Rv,v)+\lambda^{2}(\bar{R}\mathbb{E}v,\mathbb{E}v)+\lambda^{2}(Qz,z)+\lambda^{2}(\bar{Q}\mathbb{E}z,\mathbb{E}z)\nonumber\\
&~~~~~+2\lambda^{2}(S'z,v)+2\lambda^{2}(\bar{S}'\mathbb{E}z,\mathbb{E}v)+\lambda^{2}(Gz_{N},z_{N})+\lambda^{2}(\bar{G}\mathbb{E}z_{N},\mathbb{E}z_{N})
+2\lambda(q,z)\nonumber\\
&~~~~~+2\lambda(\rho,v)+2\lambda(\bar{q},\mathbb{E}z)+2\lambda(\bar{\rho},\mathbb{E}v)+2\lambda(g,z_{N})+2\lambda(\bar{g},\mathbb{E}z_{N}),
\end{align}
and
\begin{align}
J^{0}(l,0;v)&=(Rv,v)+(\bar{R}\mathbb{E}v,\mathbb{E}v)
+(Qz,z)+(\bar{Q}\mathbb{E}z,\mathbb{E}z)\nonumber\\
&~~~+2(S'z,v)+2(\bar{S}'\mathbb{E}z,\mathbb{E}v)
+(Gz_{N},z_{N})+(\bar{G}\mathbb{E}z_{N},\mathbb{E}z_{N}).
\end{align}
Notice that the first order directional derivative with the direction $v$ is obtained as
\begin{align}
dJ(l,\xi;u;v)&=\lim_{\lambda\downarrow0}\frac{J(l,\xi;u+\lambda v)-J(l,\xi;u)}{\lambda}\nonumber\\
&=2(Qx,z)+2(\bar{Q}\mathbb{E}x,\mathbb{E}z)+2(Ru,v)+2(\bar{R}\mathbb{E}u,\mathbb{E}v)+2(Gx_{N},z_{N})+2(\bar{G}\mathbb{E}x_{N},\mathbb{E}z_{N})\nonumber\\
&~~~~~~~~~~+2(S'x,v)+2(\bar{S}'\mathbb{E}x,\mathbb{E}v)+2(S'z,u)+2(\bar{S}'\mathbb{E}z,\mathbb{E}u)
+2(q,z)\nonumber\\
&~~~~~~~~~~+2(\rho,v)+2(\bar{q},\mathbb{E}z)+2(\bar{\rho},\mathbb{E}v)+2(g,z_{N})+2(\bar{g},\mathbb{E}z_{N}).
\end{align}
From (3.6)-(3.8),  (3.3) follows. On the other hand, using (3.5) yields
\begin{eqnarray*}
  \left\{
  \begin{array}{ll}
\mathbb{E}z_{k+1}=(A_{k}+\bar{A}_{k})\mathbb{E}z_{k}+(B_{k}+\bar{B}_{k})\mathbb{E}v_{k},\\
\mathbb{E}z_{l}=0,
\end{array}\right.
\end{eqnarray*}
and
\begin{eqnarray*}
  \left\{
  \begin{array}{ll}
z_{k+1}-\mathbb{E}z_{k+1}=A_{k}(z_{k}-\mathbb{E}z_{k})+B_{k}(v_{k}-\mathbb{E}v_{k})
+[C_{k}(z_{k}-\mathbb{E}z_{k})+(C_{k}+\bar{C}_{k})\mathbb{E}z_{k}\\
~~~~~~~~~~~~~~~~~~~~~~~~~+D_{k}(v_{k}-\mathbb{E}v_{k})+(D_{k}+\bar{D}_{k})\mathbb{E}v_{k}]\omega_{k},\\
z_{l}-\mathbb{E}z_{l}=0.
\end{array}\right.
\end{eqnarray*}
Besides, by (3.2), we get
\begin{eqnarray*}
  \left\{
  \begin{array}{ll}
y_{k}-\mathbb{E}y_{k}=A'_{k}[\mathbb{E}(y_{k+1}|\mathfrak{F}_{k-1})-\mathbb{E}y_{k+1}]
+C'_{k}[\mathbb{E}(y_{k+1}\omega_{k}|\mathfrak{F}_{k-1})-\mathbb{E}(y_{k+1}\omega_{k})]\\
~~~~~~~~~~~~~~~~~~+Q_{k}(x_{k}-\mathbb{E}x_{k})+S_{k}(u_{k}-\mathbb{E}u_{k})+(q_{k}+\bar{q}_{k})-\mathbb{E}(q_{k}+\bar{q}_{k}),\\
y_{N}-\mathbb{E}y_{N}=G(x_{N}-\mathbb{E}x_{N})+g+\bar{g}-\mathbb{E}(g+\bar{g}),
\end{array}\right.
\end{eqnarray*}
and
\begin{eqnarray*}
  \left\{
  \begin{array}{ll}
\mathbb{E}y_{k}=(A_{k}+\bar{A}_{k})'\mathbb{E}y_{k+1}
+(C_{k}+\bar{C}_{k})'\mathbb{E}(y_{k+1}\omega_{k})+(Q_{k}+\bar{Q}_{k})\mathbb{E}x_{k}+(S_{k}+\bar{S}_{k})\mathbb{E}u_{k}+\mathbb{E}(q_{k}+\bar{q}_{k}),\\
\mathbb{E}y_{N}=(G+\bar{G})\mathbb{E}x_{N}+\mathbb{E}(g+\bar{g}).
\end{array}\right.
\end{eqnarray*}
Thus we see
\begin{align*}
dJ(l,\xi;u;v)&=2\mathbb{E}\sum_{k=l}^{N-1}\Big\{[A'_{k}(\mathbb{E}(y_{k+1}|\mathfrak{F}_{k-1})-\mathbb{E}y_{k+1})
+C'_{k}(\mathbb{E}(y_{k+1}\omega_{k}|\mathfrak{F}_{k-1})-\mathbb{E}(y_{k+1}\omega_{k}))+Q_{k}(x_{k}-\mathbb{E}x_{k})\nonumber\\
&~~~~~~~~+S_{k}(u_{k}-\mathbb{E}u_{k})-(y_{k}-\mathbb{E}y_{k})+(q_{k}+\bar{q}_{k})-\mathbb{E}(q_{k}+\bar{q}_{k})]'(z_{k}-\mathbb{E}z_{k})+\rho_{k}+\bar{\rho}_{k}\nonumber\\
&~~~~~~~~+[B'_{k}(\mathbb{E}(y_{k+1}|\mathfrak{F}_{k-1})-\mathbb{E}y_{k+1})+D'_{k}(\mathbb{E}(y_{k+1}\omega_{k}|\mathfrak{F}_{k-1})-\mathbb{E}(y_{k+1}\omega_{k}))
+S'_{k}(x_{k}-\mathbb{E}x_{k})\\
&~~~~~~~~+R_{k}(u_{k}-\mathbb{E}u_{k})]'
(v_{k}-\mathbb{E}v_{k})\Big\}+2\mathbb{E}\sum_{k=l}^{N-1}\Big\{[(A_{k}+\bar{A}_{k})'\mathbb{E}y_{k+1}-\mathbb{E}y_{k}+\mathbb{E}(q_{k}+\bar{q}_{k})\nonumber\\
&~~~~~~~~+(C_{k}+\bar{C}_{k})'\mathbb{E}(y_{k+1}\omega_{k})+(Q_{k}+\bar{Q}_{k})\mathbb{E}x_{k}+(S_{k}+\bar{S}_{k})\mathbb{E}u_{k}]\mathbb{E}z_{k}
+[(B_{k}+\bar{B}_{k})'\mathbb{E}y_{k+1}\nonumber\\
&~~~~~~~~+\mathbb{E}(\rho_{k}+\bar{\rho}_{k})+(D_{k}+\bar{D}_{k})'\mathbb{E}(y_{k+1}\omega_{k})
+(S_{k}+\bar{S}_{k})'\mathbb{E}x_{k}+(R_{k}+\bar{R}_{k})\mathbb{E}u_{k}]'\mathbb{E}v_{k}\Big\}\nonumber\\
&=2\mathbb{E}\sum_{k=l}^{N-1}\Big\{[B'_{k}(\mathbb{E}(y_{k+1}|\mathfrak{F}_{k-1})-\mathbb{E}y_{k+1})
+D'_{k}(\mathbb{E}(y_{k+1}\omega_{k}|\mathfrak{F}_{k-1})-\mathbb{E}(y_{k+1}\omega_{k}))+\rho_{k}+\bar{\rho}_{k}\nonumber\\
&~~~~~~~~+S'_{k}(x_{k}-\mathbb{E}x_{k})+R_{k}(u_{k}-\mathbb{E}u_{k})]'(v_{k}-\mathbb{E}v_{k})+2[(B_{k}+\bar{B}_{k})'\mathbb{E}y_{k+1}+\mathbb{E}(\rho_{k}+\bar{\rho}_{k})\nonumber\\
&~~~~~~~~+(D_{k}+\bar{D}_{k})'\mathbb{E}(y_{k+1}\omega_{k})
+(S_{k}+\bar{S}_{k})'\mathbb{E}x_{k}+(R_{k}+\bar{R}_{k})\mathbb{E}u_{k}]'\mathbb{E}v_{k}\Big\}\nonumber\\
&=2\mathbb{E}\sum_{k=l}^{N-1}\Big\{[B'_{k}\mathbb{E}(y_{k+1}|\mathfrak{F}_{k-1})+\bar{B}'_{k}\mathbb{E}y_{k+1}+\rho_{k}+\bar{\rho}_{k}+D'_{k}\mathbb{E}(y_{k+1}\omega_{k}|\mathfrak{F}_{k-1})\nonumber\\
&~~~~~~~~+\bar{D}'_{k}\mathbb{E}(y_{k+1}\omega_{k})
+S'_{k}x_{k}+\bar{S}'_{k}\mathbb{E}x_{k}+R_{k}u_{k}+\bar{R}_{k}\mathbb{E}u_{k}]'v_{k}\Big\},
\end{align*}
i.e., the proof is finished.
\end{proof}

\vskip0mm\noindent
\textbf{Remark 3.1.} For any $l\in\mathbb{N}_{0}$, there are bounded linear operators $\mathcal{L}:\mathscr{U}_{ad}\rightarrow L^{2}_{\mathfrak{F}}(\mathbb{N};\mathbb{R}^{n})$ and $\bar{\mathcal{L}}:\mathscr{U}_{ad}\rightarrow L^{2}_{\mathfrak{F}}(N-1;\mathbb{R}^{n})$  such that
\begin{eqnarray*}
J^{0}(l,0;u)=\langle \mathcal{M}u,u\rangle,
\end{eqnarray*}
where
\begin{align*}
&\mathcal{M}=\mathcal{L}^{\ast}(Q+\mathbb{E}^{\ast}\bar{Q}\mathbb{E})\mathcal{L}
+\bar{\mathcal{L}}^{\ast}(G+\mathbb{E}^{\ast}\bar{G}\mathbb{E})\bar{\mathcal{L}}+R+\mathbb{E}^{\ast}\bar{R}\mathbb{E}+(S+\mathbb{E}^{\ast}\bar{S}\mathbb{E})\mathcal{L}+\mathcal{L}^{\ast}(S'+\mathbb{E}^{\ast}\bar{S}'\mathbb{E}).
\end{align*}
Using Lemma 3.1, the cost functional could be expressed as
\begin{align}
J(l,\xi;u)=\langle \mathcal{M}u,u\rangle+J(l,\xi;0)+\mathbb{E}\sum\limits_{k=l}^{N-1}\langle dJ(l,\xi;0),u\rangle.
\end{align}

\vskip0mm\noindent
\textbf{Corollary 3.1.} The following statements are equivalent:

\vskip0mm\noindent
(i) $J^{0}(l,0;u)\geq 0$ for all $u\in \mathscr{U}_{ad}$.

\noindent
(ii) The map $u\mapsto J(l,\xi;u)$ is convex for some $\xi\in \mathbb{R}^{n}$.

\noindent
(iii) The map $u\mapsto J(l,\xi;u)$ is convex for any $\xi\in \mathbb{R}^{n}$.

\noindent
(iv) The map $u\mapsto J^{0}(l,\xi;u)$ is convex for some $\xi\in \mathbb{R}^{n}$.

\noindent
(v) The map $u\mapsto J^{0}(l,\xi;u)$ is convex for any $\xi\in \mathbb{R}^{n}$.

\noindent
(vi) $\mathcal{M}\geq0$.

\vskip2mm\noindent
\textbf{Lemma 3.2.} Let the standard condition (1.5) hold, then the map $u\longmapsto J^{0}(l,0;u)$ is uniformly convex for any $l\in \mathbb{N}_{0}$.
\begin{proof} For any $u\in \mathscr{U}_{ad}$, let $x^{(u)}$ be the solution of (2.5), then we obtain
\begin{align*}
J^{0}(l,0;u)&=\mathbb{E}\left\{\sum_{k=l}^{N-1}\left\langle  \left(
                            \begin{array}{cc}
                              Q_{k} & S_{k}' \\
                              S_{k} & R_{k}\\
                            \end{array}
                          \right)\left(
                                   \begin{array}{c}
                                     x_{k}^{(u)}-\mathbb{E}x_{k}^{(u)} \\
                                     u_{k}-\mathbb{E}u_{k} \\
                                   \end{array}
                                 \right),\left(
                                   \begin{array}{c}
                                     x_{k}^{(u)}-\mathbb{E}x_{k}^{(u)}  \\
                                     u_{k}-\mathbb{E}u_{k} \\
                                   \end{array}
                                 \right)\right\rangle+\left\langle G\left(x_{N}^{(u)}-\mathbb{E}x_{N}^{(u)}\right),x_{N}^{(u)}-\mathbb{E}x_{N}^{(u)}\right\rangle\right\}\nonumber\\
&~~~~~~~~+\mathbb{E}\left\{\sum_{k=l}^{N-1}\left\langle  \left(
                            \begin{array}{cc}
                              Q_{k}+\bar{Q}_{k} & S_{k}'+\bar{S}_{k}' \\
                              S_{k}+\bar{S}_{k} & R_{k}+\bar{R}_{k} \\
                            \end{array}
                          \right)\left(
                                   \begin{array}{c}
                                     \mathbb{E}x_{k}^{(u)} \\
                                     \mathbb{E}u_{k} \\
                                   \end{array}
                                 \right),\left(
                                   \begin{array}{c}
                                     \mathbb{E}x_{k}^{(u)} \\
                                     \mathbb{E}u_{k} \\
                                   \end{array}
                                 \right)\right\rangle+\left\langle (G+\bar{G})\mathbb{E}x_{N}^{(u)},\mathbb{E}x_{N}^{(u)}\right\rangle\right\}\nonumber\\
&\geq \mathbb{E}\sum_{k=l}^{N-1}\bigg\{\left\langle Q_{k}\left(x_{k}^{(u)}-\mathbb{E}x_{k}^{(u)} \right),x_{k}^{(u)}-\mathbb{E}x_{k}^{(u)}\right\rangle+
2\left\langle S_{k}\left(x_{k}^{(u)}-\mathbb{E}x_{k}^{(u)} \right),u_{k}-\mathbb{E}u_{k}\right\rangle\nonumber\\
&~~~~~~~~+\Big\langle R_{k}\left(u_{k}-\mathbb{E}u_{k} \right),u_{k}-\mathbb{E}u_{k} \Big\rangle+\left\langle (Q_{k}+\bar{Q}_{k})\mathbb{E}x_{k}^{(u)},\mathbb{E}x_{k}^{(u)}\right\rangle
+2\Big\langle (S_{k}+\bar{S}_{k})\mathbb{E}u_{k},\mathbb{E}u_{k} \Big\rangle
\bigg\}\nonumber\\
&=\mathbb{E}\sum_{k=l}^{N-1}\Big\{ \Big\langle (R_{k}+\bar{R}_{k})\left[\mathbb{E}u_{k}+(R_{k}+\bar{R}_{k})^{-1}(S_{k}+\bar{S}_{k})\mathbb{E}x_{k}^{(u)}\right],\mathbb{E}u_{k}+(R_{k}+\bar{R}_{k})^{-1}\nonumber\\
&~~~~~~~~\times(S_{k}+\bar{S}_{k})\mathbb{E}x_{k}^{(u)} \Big\rangle+\left\langle (Q_{k}-S'_{k}R^{-1}_{k}S_{k})\left(x_{k}^{(u)}-\mathbb{E}x_{k}^{(u)} \right),x_{k}^{(u)}-\mathbb{E}x_{k}^{(u)}\right\rangle\nonumber\\
&~~~~~~~~+\left\langle R_{k}\left[u_{k}-\mathbb{E}u_{k}+R^{-1}_{k}S_{k}\left(x_{k}^{(u)}-\mathbb{E}x_{k}^{(u)} \right)\right],u_{k}-\mathbb{E}u_{k}+R^{-1}_{k}S_{k}\left(x_{k}^{(u)}-\mathbb{E}x_{k}^{(u)} \right) \right\rangle\nonumber\\
&~~~~~~~~+\left\langle [Q_{k}+\bar{Q}_{k}-(S_{k}+\bar{S}_{k})'(R_{k}+\bar{R}_{k})^{-1}(S_{k}+\bar{S}_{k})]\mathbb{E}x_{k}^{(u)},\mathbb{E}x_{k}^{(u)}\right\rangle\Big\}.
\end{align*}
Hence, combining (1.5) with Lemma 2.3, there exist $\delta>0$ and $\mu>0$ such that
\begin{align*}
J^{0}(l,0;u)
&\geq \delta\mathbb{E}\sum_{k=l}^{N-1}\left\{\left|u_{k}-\mathbb{E}u_{k}+R^{-1}_{k}S_{k}\left(x_{k}^{(u)}-\mathbb{E}x_{k}^{(u)}\right)\right|^{2}
+\left|\mathbb{E}u_{k}+(R_{k}+\bar{R}_{k})^{-1}(S_{k}+\bar{S}_{k})\mathbb{E}x_{k}^{(u)}\right|^{2}\right\}\\
&\geq \delta\mathbb{E}\sum_{k=l}^{N-1}\left\{\left|u_{k}-\mathbb{E}u_{k}+R^{-1}_{k}S_{k}\left(x_{k}^{(u)}-\mathbb{E}x_{k}^{(u)}\right)\right|^{2}+\mu|\mathbb{E}u_{k}|^{2}\right\}\nonumber\\
&\geq \frac{\delta\mu}{1+\mu}\mathbb{E}\sum_{k=l}^{N-1}\left|u_{k}+R^{-1}_{k}S_{k}\left(x_{k}^{(u)}-\mathbb{E}x_{k}^{(u)}\right)\right|^{2}\nonumber\\
&\geq \frac{\delta\mu^{2}}{1+\mu}\mathbb{E}\sum_{k=l}^{N-1}|u_{k}|^{2},~~\forall~u_{k}\in \mathscr{U}_{ad}.
\end{align*}
The proof is completed.
\end{proof}

\setcounter{equation}{0}
\section{Solvabilities of Problem (MF-LQ), uniform convexity of the cost functional, and GRE}
We begin with the open-loop optimal control of Problem (MF-LQ) on the basis of MF-FBSDEs.

\vskip3mm\noindent
\textbf{Theorem 4.1.} Given the initial pair $(l,\xi)$, let $u^{\ast}\in \mathscr{U}_{ad}$ and $(x,y)$ be the solution to the following MF-FBSDEs
\begin{eqnarray}\label{41}
  \left\{
  \begin{array}{ll}
x_{k+1}=A_{k}x_{k}+\bar{A}_{k}\mathbb{E}x_{k}+B_{k}u_{k}^{\ast}+\bar{B}_{k}\mathbb{E}u_{k}^{\ast}+b_{k}\\
~~~~~~~~~~~~+(C_{k}x_{k}+\bar{C}_{k}\mathbb{E}x_{k}+D_{k}u_{k}^{\ast}+\bar{D}_{k}\mathbb{E}u_{k}^{\ast}+\sigma_{k})\omega_{k},\\
y_{k}=A'_{k}\mathbb{E}(y_{k+1}|\mathfrak{F}_{k-1})+\bar{A}'_{k}\mathbb{E}y_{k+1}
+C'_{k}\mathbb{E}(y_{k+1}\omega_{k}|\mathfrak{F}_{k-1})+\bar{C}'_{k}\mathbb{E}(y_{k+1}\omega_{k})\\
~~~~~~~~~~~~+Q_{k}x_{k}+\bar{Q}_{k}\mathbb{E}x_{k}+S_{k}u_{k}^{\ast}+\bar{S}_{k}\mathbb{E}u_{k}^{\ast}+q_{k}+\bar{q}_{k},\\
x_{l}=\xi,~~y_{N}=Gx_{N}+\bar{G}\mathbb{E}x_{N}+g+\bar{g}.
\end{array}\right.
\end{eqnarray}
Then $u^{\ast}$ is an open-loop optimal control of Problem (MF-LQ) if and only if
\begin{align}\label{42}
J^{0}(l,0;u)\geq0,~~\forall~u~\in \mathscr{U}_{ad},
\end{align}
and satisfies the following stationarity condition
\begin{align}\label{43}
&B'_{k}\mathbb{E}(y_{k+1}|\mathfrak{F}_{k-1})+\bar{B}'_{k}\mathbb{E}y_{k+1}+\rho_{k}+\bar{\rho}_{k}+D'_{k}\mathbb{E}(y_{k+1}\omega_{k}|\mathfrak{F}_{k-1})\nonumber\\
&~~~~~~~~~~~~~~~~~~+\bar{D}'_{k}\mathbb{E}(y_{k+1}\omega_{k})
+S'_{k}x_{k}+\bar{S}'_{k}\mathbb{E}x_{k}+R_{k}u_{k}^{\ast}+\bar{R}_{k}\mathbb{E}u_{k}^{\ast}=0.
\end{align}
In this case, the optimal value of Problem (MF-LQ) at $(l,\xi)$ is derived as
\begin{align*}
V(l,\xi)=\mathbb{E}(y_{l}'\xi).
\end{align*}
\begin{proof} According to Lemma 3.1, we can immediately derive that $u_{k}^{\ast}$ is an optimal control of Problem (MF-LQ) for the initial pair $(l,\xi)$ if and only if
\begin{align*}
\lambda^{2}J^{0}(l,0;v_{k})+\lambda dJ(l,\xi;u_{k};v_{k})=J(l,\xi;u_{k}+\lambda v_{k})-J(l,\xi;u_{k})
\end{align*}
holds for all $\lambda\in \mathbb{R}$ and $u\in \mathscr{U}_{ad}$,
which is equivalent to (4.2) and $dJ(l,\xi;u_{k};v_{k})\leq0.$
Since the above inequality holds for all $u\in \mathscr{U}_{ad}$ if and only if $dJ(l,\xi;u_{k};v_{k})=0$, then (4.3) follows. Next, define the Lyapunov function candidate as
\begin{align*}
V_{N}(k,x_{k})=\mathbb{E}[y_{k}'x_{k}],
\end{align*}
then we obtain
\begin{align}\label{44}
&V_{N}(k,x_{k})-V_{N}(k+1,x_{k+1})
=\mathbb{E}[y_{k}'x_{k}]-\mathbb{E}[y_{k+1}'x_{k+1}]\nonumber\\
&~~~~~~=\mathbb{E}\Big\{[A'_{k}\mathbb{E}(y_{k+1}|\mathfrak{F}_{k-1})+\bar{A}'_{k}\mathbb{E}y_{k+1}
+C'_{k}\mathbb{E}(y_{k+1}\omega_{k}|\mathfrak{F}_{k-1})+\bar{C}'_{k}\mathbb{E}(y_{k+1}\omega_{k})+Q_{k}x_{k}+\bar{Q}_{k}\mathbb{E}x_{k}\nonumber\\
&~~~~~~~~~~~~~~+S_{k}u_{k}+\bar{S}_{k}\mathbb{E}u_{k}+q_{k}+\bar{q}_{k}]'x_{k}\Big\}-\mathbb{E}\Big\{y_{k+1}'[A_{k}x_{k}
+\bar{A}_{k}\mathbb{E}x_{k}+B_{k}u_{k}^{\ast}+\bar{B}_{k}\mathbb{E}u_{k}^{\ast}+b_{k}\nonumber\\
&~~~~~~~~~~~~~~+(C_{k}x_{k}+\bar{C}_{k}\mathbb{E}x_{k}+D_{k}u_{k}^{\ast}+\bar{D}_{k}\mathbb{E}u_{k}^{\ast}+\sigma_{k})\omega_{k}]'\Big\}\nonumber\\
&~~~~~~=x_{k}'Q_{k}x_{k}+(\mathbb{E}x_{k})'\bar{Q}_{k}\mathbb{E}x_{k}
+2x_{k}'S_{k}u_{k}+2(\mathbb{E}x_{k})'\bar{S}_{k}\mathbb{E}u_{k}
+u_{k}'R_{k}u_{k}+(\mathbb{E}u_{k})'\bar{R}_{k}\mathbb{E}u_{k}\nonumber\\
&~~~~~~~~~~~~~~+2q_{k}'x_{k}+2\rho_{k}'u_{k}+2\bar{q}_{k}'\mathbb{E}x_{k}+2\bar{\rho}_{k}'\mathbb{E}u_{k}-\mathbb{E}[B'_{k}\mathbb{E}(y_{k+1}|\mathfrak{F}_{k-1})+\bar{B}'_{k}\mathbb{E}y_{k+1}
+\rho_{k}+\bar{\rho}_{k}\nonumber\\
&~~~~~~~~~~~~~~+D'_{k}\mathbb{E}(y_{k+1}\omega_{k}|\mathfrak{F}_{k-1})
+\bar{D}'_{k}\mathbb{E}(y_{k+1}\omega_{k})+S'_{k}x_{k}+\bar{S}'_{k}\mathbb{E}x_{k}+R_{k}u_{k}+\bar{R}_{k}\mathbb{E}u_{k}].
\end{align}
Adding from $k=l$ to $k=N-1$ on the both sides of (4.4), we get
\begin{align*}
\mathbb{E}[y_{l}'x_{l}]&=\mathbb{E}\sum\limits_{k=l}^{N-1}\{x_{k}'Q_{k}x_{k}+(\mathbb{E}x_{k})'\bar{Q}_{k}\mathbb{E}x_{k}
+2x_{k}'S_{k}u_{k}^{\ast}+2(\mathbb{E}x_{k})'\bar{S}_{k}\mathbb{E}u_{k}^{\ast}
+(u_{k}^{\ast})'R_{k}u_{k}^{\ast}+(\mathbb{E}u_{k}^{\ast})'\bar{R}_{k}\mathbb{E}u_{k}^{\ast}\\
&~~~~~~~~~~~~+2q_{k}'x_{k}+2\rho_{k}'u_{k}^{\ast}+2\bar{q}_{k}'\mathbb{E}x_{k}+2\bar{\rho}_{k}'\mathbb{E}u_{k}^{\ast}\}
+\mathbb{E}(x_{N}'Gx_{N})+(\mathbb{E}x_{N})'\bar{G}\mathbb{E}x_{N}\\
&=J(l;\xi;u^{\ast}).
\end{align*}
Thus we can see that the optimal value $V(l,\xi)=\mathbb{E}(y_{l}'\xi)$.
\end{proof}

\vskip0mm\noindent
\textbf{Remark 4.1.} The stationarity condition (4.3) takes a coupling into the MF-FBSDEs (4.1). (4.2) and (4.3) make up of the optimality system for the open-loop optimal control of Problem (MF-LQ).

\vskip3mm\noindent
\textbf{Lemma 4.1.} Assume that ($\mathbf{A2}$) holds, then Problem (MF-LQ) is uniquely open-loop solvable at the initial pair $(l,\xi)$, and the optimal control
\begin{align*}
u^{\ast}_{k}&=-\frac{1}{2}\mathcal{M}^{-1}dJ(l,\xi;0),
\end{align*}
besides, the optimal value
\begin{align}\label{45}
V(l,\xi)=J(l,\xi;0)-\frac{1}{4}\Big|\mathcal{M}^{-\frac{1}{2}}dJ(l,\xi;0)\Big|^{2}.
\end{align}
Moreover, there admits a constant $\beta\in \mathbb{R}$ (not necessarily nonnegative) such that
\begin{align}\label{46}
V^{0}(l,\xi)\geq \beta|\xi|^{2}.
\end{align}
\begin{proof} Under ($\mathbf{A2}$), the operator $\mathcal{M}$ is invertible, and
\begin{align*}
J(l,\xi;u)&=\langle \mathcal{M}u,u\rangle+J(l,\xi;0)+\mathbb{E}\sum\limits_{k=l}^{N-1}\langle dJ(l,\xi;0),u\rangle\\
&=\left|\mathcal{M}^{\frac{1}{2}}u+\frac{1}{2}\mathcal{M}^{-\frac{1}{2}}dJ(l,\xi;0)\right|^{2}
-\frac{1}{4}\left|\mathcal{M}^{-\frac{1}{2}}dJ(l,\xi;0)\right|^{2}+J(l,\xi;0)\\
&\geq J(l,\xi;0)-\frac{1}{4}\left|\mathcal{M}^{-\frac{1}{2}}dJ(l,\xi;0)\right|^{2}.
\end{align*}
Since the above equality holds for all $\xi\in L_{\mathfrak{F}}^{2}(l;\mathbb{R}^{n})$ and $u\in \mathscr{U}_{ad}$ if and only if
$u^{\ast}_{k}=-\frac{1}{2}\mathcal{M}^{-1}dJ(l,\xi;0)$,
and in this case, the optimal value (4.5) follows. On the other hand, by the classical results on convex analysis, the uniformly convexity of map $u\mapsto J^{0}(0,0;u)$ is equivalent to the uniformly convexity of map $u\mapsto J^{0}(l,\xi;u)$. Then we have
\begin{align*}
J(l,\xi;u_{k})&=J(l,\xi;0)+J^{0}(l,0;u_{k})+\mathbb{E}\sum\limits_{k=l}^{N-1}\langle dJ(l,\xi;0),u_{k}\rangle\nonumber\\
&\geq J(l,\xi;0)+J^{0}(l,0;u_{k})-\frac{\delta}{2}\mathbb{E}\sum\limits_{k=l}^{N-1}|u_{k}|^{2}-\frac{1}{2\delta}\mathbb{E}\sum\limits_{k=l}^{N-1}|dJ(l,\xi;0)|^{2}\\
&\geq J(l,\xi;0)+\frac{\delta}{2}\mathbb{E}\sum\limits_{k=l}^{N-1}|u_{k}|^{2}-\frac{1}{2\delta}\mathbb{E}\sum\limits_{k=l}^{N-1}|dJ(l,\xi;0)|^{2}.
\end{align*}
When $b_{k}, \sigma_{k}, q_{k}, \bar{q}_{k}, \rho_{k}, \bar{\rho}_{k}, g, \bar{g}=0$, we further get that
\begin{align*}
V^{0}(l,\xi)\geq J^{0}(l,\xi;0)-\frac{1}{2\delta}\mathbb{E}\sum\limits_{k=l}^{N-1}|dJ^{0}(l,\xi;0)|^{2}.
\end{align*}
Because $J^{0}(l,\xi;0)$ and $dJ^{0}(l,\xi;0)$ are quadratic in $\xi$ and continuous in $l$, (4.6) follows.
\end{proof}

For future references, we introduce the following GRE associated with Problem (MF-LQ)
\begin{eqnarray}\label{47}
  \left\{
  \begin{array}{ll}
P_{k}=Q_{k}+A_{k}'P_{k+1}A_{k}+C_{k}'P_{k+1}C_{k}-(B_{k}'P_{k+1}A_{k}+D_{k}'P_{k+1}C_{k}+S'_{k})'\\
~~~~~~~~~~\times(R_{k}+B_{k}'P_{k+1}B_{k}+D_{k}'P_{k+1}D_{k})^{\dagger}
(B_{k}'P_{k+1}A_{k}+D_{k}'P_{k+1}C_{k}+S'_{k}),\\
\Pi_{k}=Q_{k}+\bar{Q}_{k}+(A_{k}+\bar{A}_{k})'\Pi_{k+1}(A_{k}+\bar{A}_{k})+(C_{k}+\bar{C}_{k})'P_{k+1}(C_{k}+\bar{C}_{k})\\
~~~~~~~~~~-[(B_{k}+\bar{B}_{k})'\Pi_{k+1}(A_{k}+\bar{A}_{k})+(D_{k}+\bar{D}_{k})'P_{k+1}(C_{k}+\bar{C}_{k})+S'_{k}+\bar{S}'_{k}]'\\
~~~~~~~~~~\times[R_{k}+\bar{R}_{k}+(B_{k}+\bar{B}_{k})'\Pi_{k+1}(B_{k}+\bar{B}_{k})+(D_{k}+\bar{D}_{k})'P_{k+1}(D_{k}+\bar{D}_{k})]^{\dagger}\\
~~~~~~~~~~\times[(B_{k}+\bar{B}_{k})'\Pi_{k+1}(A_{k}+\bar{A}_{k})+(D_{k}+\bar{D}_{k})'P_{k+1}(C_{k}+\bar{C}_{k})+S'_{k}+\bar{S}'_{k}],\\
P_{N}=G,~~\Pi_{N}=G+\bar{G}.
\end{array}\right.
\end{eqnarray}

\vskip0mm\noindent
\textbf{Definition 4.1.}

\noindent
(i) A solution $(P_{k},\Pi_{k})$ of the GRE (4.7) is called to be regular if
\begin{align}
  \left\{
  \begin{array}{ll}
\Upsilon_{k}=R_{k}+B'_{k}P_{k+1}B_{k}+D'_{k}P_{k+1}D_{k}\geq0,\\
\bar{\Upsilon}_{k}=R_{k}+\bar{R}_{k}+(B_{k}+\bar{B}_{k})'\Pi_{k+1}(B_{k}+\bar{B}_{k})+(D_{k}+\bar{D}_{k})'P_{k+1}(D_{k}+\bar{D}_{k})\geq0,
\end{array}\right.\\
  \left\{
  \begin{array}{ll}
\mathcal{R}(B_{k}'P_{k+1}A_{k}+D_{k}'P_{k+1}C_{k}+S_{k}')\subseteq \mathcal{R}(R_{k}+B_{k}'P_{k+1}B_{k}+D_{k}'P_{k+1}D_{k}),\\
\mathcal{R}[(B_{k}+\bar{B}_{k})'\Pi_{k+1}(A_{k}+\bar{A}_{k})+(D_{k}+\bar{D}_{k})'P_{k+1}(C_{k}+\bar{C}_{k})+S_{k}'+\bar{S}_{k}']\\
~~~\subseteq \mathcal{R}[R_{k}+\bar{R}_{k}+(B_{k}+\bar{B}_{k})'\Pi_{k+1}(B_{k}+\bar{B}_{k})+(D_{k}+\bar{D}_{k})'P_{k+1}(D_{k}+\bar{D}_{k})],
\end{array}\right.
\end{align}
and
\begin{align}\label{410}
  \left\{
  \begin{array}{ll}
(R_{k}+B_{k}'P_{k+1}B_{k}+D_{k}'P_{k+1}D_{k})^{\dagger}(B_{k}'P_{k+1}A_{k}+D_{k}'P_{k+1}C_{k}+S_{k}')\in L^{2}(\mathbb{N};\mathbb{R}^{m\times n}),\\
\Big[R_{k}+\bar{R}_{k}+(B_{k}+\bar{B}_{k})'\Pi_{k+1}(B_{k}+\bar{B}_{k})+(D_{k}+\bar{D}_{k})'P_{k+1}(D_{k}+\bar{D}_{k})\Big]^{\dagger}\Big[S_{k}'+\bar{S}_{k}'\\
~~~~~~~~+(B_{k}+\bar{B}_{k})'\Pi_{k+1}(A_{k}+\bar{A}_{k})+(D_{k}+\bar{D}_{k})'P_{k+1}(C_{k}+\bar{C}_{k})\Big]\in L^{2}(\mathbb{N};\mathbb{R}^{m\times n}).
\end{array}\right.
\end{align}
(ii) A solution $(P_{k},\Pi_{k})$ of the GRE (4.7) is called to be strongly regular if for some $\alpha>0$,
\begin{align}\label{411}
\Upsilon_{k}\geq\alpha I,~~~~~~\bar{\Upsilon}_{k}\geq\alpha I.
\end{align}
The GRE (4.7) is called to be regularly solvable if it exists a regular solution, and it is called to be strongly regularly solvable if it exists a strongly regular solution.

\vskip2mm
Let $\Theta_{k}^{\ast}, \bar{\Theta}_{k}^{\ast}\in \mathbb{R}^{m\times n}$, we consider the following state equation
\begin{eqnarray*}
  \left\{
  \begin{array}{ll}
x_{k+1}=(A_{k}+B_{k}\Theta_{k}^{\ast})x_{k}+[\bar{A}_{k}+B_{k}\bar{\Theta}_{k}^{\ast}+\bar{B}_{k}(\Theta_{k}^{\ast}+\bar{\Theta}_{k}^{\ast})]\mathbb{E}x_{k}
+B_{k}u_{k}+\bar{B}_{k}\mathbb{E}u_{k}+b_{k}\\
~~~~~~~~~~~~~~+\Big\{(C_{k}+D_{k}\Theta_{k}^{\ast})x_{k}+[\bar{C}_{k}+D_{k}\bar{\Theta}_{k}^{\ast}+\bar{D}_{k}(\Theta_{k}^{\ast}+\bar{\Theta}_{k}^{\ast})]\mathbb{E}x_{k}
+D_{k}u_{k}+\bar{D}_{k}\mathbb{E}u_{k}+\sigma_{k}\Big\}\omega_{k},\\
x_{l}=\xi,
\end{array}\right.
\end{eqnarray*}
and cost functional
\begin{align}\label{412}
\tilde{J}(l,\xi;u)&\doteq\mathbb{E}\Big\{\langle Gx_{N},x_{N}\rangle+2\langle g,x_{N}\rangle+\langle \bar{G}\mathbb{E}x_{N},\mathbb{E}x_{N}\rangle+2\langle \bar{g},\mathbb{E}x_{N}\rangle\Big\}\nonumber\\
&~~~~~~~+\mathbb{E}\sum_{k=l}^{N-1}\left\langle  \left(
                            \begin{array}{cc}
                              Q_{k} & S_{k}' \\
                              S_{k} & R_{k}\\
                            \end{array}
                          \right)\left(
                                   \begin{array}{c}
                                    x_{k} \\
                                     \Theta^{\ast}_{k}x_{k}+\bar{\Theta}^{\ast}_{k}\mathbb{E}x_{k}+u_{k} \\
                                   \end{array}
                                 \right),\left(
                                   \begin{array}{c}
                                     x_{k} \\
                                    \Theta^{\ast}_{k}x_{k}+\bar{\Theta}^{\ast}_{k}\mathbb{E}x_{k}+u_{k} \\
                                  \end{array}
                                 \right)\right\rangle\nonumber\\
&~~~~~~~+\mathbb{E}\sum_{k=l}^{N-1}\left\langle  \left(
                           \begin{array}{cc}
                             \bar{Q}_{k} & \bar{S}_{k}' \\
                              \bar{S}_{k} & \bar{R}_{k} \\
                           \end{array}
                          \right)\left(
                                   \begin{array}{c}
                                     \mathbb{E}x_{k} \\
                                    (\Theta^{\ast}_{k}+\bar{\Theta}^{\ast}_{k})\mathbb{E}x_{k}+\mathbb{E}u_{k} \\
                                   \end{array}
                                 \right),\left(
                                   \begin{array}{c}
                                     \mathbb{E}x_{k} \\
                                     (\Theta^{\ast}_{k}+\bar{\Theta}^{\ast}_{k})\mathbb{E}x_{k}+\mathbb{E}u_{k} \\
                                  \end{array}
                                 \right)\right\rangle\nonumber\\
&~~~~~~~+2\mathbb{E}\sum_{k=l}^{N-1}\left[\left\langle  \left(
                                                            \begin{array}{c}
                                                              q_{k} \\
                                                              \rho_{k} \\
                                                            \end{array}
                                                          \right),\left(
                                   \begin{array}{c}
                                     x_{k} \\
                                     \Theta^{\ast}_{k}x_{k}+\bar{\Theta}^{\ast}_{k}\mathbb{E}x_{k}+u_{k} \\
                                   \end{array}
                                 \right)
                                \right\rangle+\left\langle  \left(
                                                            \begin{array}{c}
                                                             \bar{q}_{k} \\
                                                              \bar{\rho}_{k} \\
                                                           \end{array}
                                                         \right),\left(
                                  \begin{array}{c}
                                    \mathbb{E}x_{k} \\
                                     (\Theta^{\ast}_{k}+\bar{\Theta}^{\ast}_{k})\mathbb{E}x_{k}+\mathbb{E}u_{k} \\
                                  \end{array}
                                 \right)
                                 \right\rangle\right]\nonumber\\
&=\mathbb{E}\Big\{\langle Gx_{N},x_{N}\rangle+2\langle g,x_{N}\rangle+\langle \bar{G}\mathbb{E}x_{N},\mathbb{E}x_{N}\rangle+2\langle \bar{g},\mathbb{E}x_{N}\rangle\Big\}\nonumber\\
&~~~~~~~+\mathbb{E}\sum_{k=l}^{N-1}\left[\left\langle  \left(
                            \begin{array}{cc}
                              \tilde{Q}_{k} & \tilde{S}_{k}' \\
                              \tilde{S}_{k} & R_{k}\\
                            \end{array}
                          \right)\left(
                                  \begin{array}{c}
                                     x_{k} \\
                                    u_{k} \\
                                  \end{array}
                                \right),\left(
                                   \begin{array}{c}
                                    x_{k} \\
                                    u_{k} \\
                                 \end{array}
                                \right)\right\rangle+2\left\langle  \left(
                                                           \begin{array}{c}
                                                             \tilde{q}_{k} \\
                                                              \rho_{k} \\
                                                            \end{array}
                                                          \right),\left(
                                  \begin{array}{c}
                                     x_{k} \\
                                    u_{k} \\
                                   \end{array}
                                 \right)
                                 \right\rangle\right]\nonumber\\
&~~~~~~~+\mathbb{E}\sum_{k=l}^{N-1}\left[\left\langle  \left(
                            \begin{array}{cc}
                              \hat{Q}_{k} & \hat{S}_{k}' \\
                              \hat{S}_{k} & \bar{R}_{k} \\
                            \end{array}
                          \right)\left(
                                   \begin{array}{c}
                                     \mathbb{E}x_{k} \\
                                    \mathbb{E}u_{k} \\
                                   \end{array}
                                 \right),\left(
                                   \begin{array}{c}
                                     \mathbb{E}x_{k} \\
                                    \mathbb{E}u_{k} \\
                                  \end{array}
                                 \right)\right\rangle+2\left\langle  \left(
                                                            \begin{array}{c}
                                                              \hat{q}_{k} \\
                                                              \bar{\rho}_{k} \\
                                                           \end{array}
                                                          \right),\left(
                                   \begin{array}{c}
                                    \mathbb{E}x_{k} \\
                                    \mathbb{E}u_{k} \\
                                  \end{array}
                                 \right)
                                \right\rangle\right],
\end{align}
where
\begin{eqnarray*}
  \left\{
  \begin{array}{ll}
\tilde{Q}_{k}=Q_{k}+(\Theta^{\ast}_{k})'S_{k}+S_{k}'\Theta^{\ast}_{k}+(\Theta^{\ast}_{k})'R_{k}\Theta_{k}^{\ast},~\tilde{S}_{k}=S_{k}+R_{k}\Theta_{k}^{\ast},
~\tilde{q}_{k}=q_{k}+(\Theta^{\ast}_{k})'\rho_{k},\\
\hat{Q}_{k}=\bar{Q}_{k}+(\Theta^{\ast}_{k}+\bar{\Theta}^{\ast}_{k})'\bar{S}_{k}+\bar{S}_{k}'(\Theta^{\ast}_{k}+\bar{\Theta}^{\ast}_{k})+(\Theta^{\ast}_{k}+\bar{\Theta}^{\ast}_{k})'
\bar{R}_{k}(\Theta^{\ast}_{k}+\bar{\Theta}^{\ast}_{k})\\
~~~~~~~~~~~~+(\bar{\Theta}^{\ast}_{k})'R_{k}\bar{\Theta}^{\ast}_{k}+(\bar{\Theta}^{\ast}_{k})'S_{k}+S_{k}'\bar{\Theta}^{\ast}_{k}
+(\bar{\Theta}^{\ast}_{k})'R_{k}\Theta^{\ast}_{k}+(\Theta^{\ast}_{k})'R_{k}\bar{\Theta}^{\ast}_{k},\\
\hat{S}_{k}=\bar{S}_{k}+\bar{R}_{k}(\Theta^{\ast}_{k}+\bar{\Theta}^{\ast}_{k})+R_{k}\bar{\Theta}^{\ast}_{k},
~\hat{q}_{k}=\bar{q}_{k}+(\Theta^{\ast}_{k}+\bar{\Theta}^{\ast}_{k})'\bar{\rho}_{k}+(\bar{\Theta}^{\ast}_{k})'\mathbb{E}\rho_{k}.
\end{array}\right.
\end{eqnarray*}
\vskip1mm\noindent
\textbf{Theorem 4.2.}

\noindent
(i) Given the initial pair $(l,\xi)$, $u^{\ast}\in \mathscr{U}_{ad}$, and $(x^{\ast},y^{\ast})$ is the solution to the following MF-FBSDEs:
\begin{eqnarray}\label{1}
  \left\{
  \begin{array}{ll}
x_{k+1}^{\ast}=(A_{k}+B_{k}\Theta_{k}^{\ast})x_{k}^{\ast}+[\bar{A}_{k}+B_{k}\bar{\Theta}_{k}^{\ast}+\bar{B}_{k}(\Theta_{k}^{\ast}+\bar{\Theta}_{k}^{\ast})]\mathbb{E}x_{k}^{\ast}
+B_{k}u_{k}^{\ast}+\bar{B}_{k}\mathbb{E}u_{k}^{\ast}+b_{k}\\
~~~~~~~~~~~+\Big\{(C_{k}+D_{k}\Theta_{k}^{\ast})x_{k}^{\ast}+[\bar{C}_{k}+D_{k}\bar{\Theta}_{k}^{\ast}+\bar{D}_{k}(\Theta_{k}^{\ast}+\bar{\Theta}_{k}^{\ast})]\mathbb{E}x_{k}^{\ast}
+D_{k}u_{k}^{\ast}+\bar{D}_{k}\mathbb{E}u_{k}^{\ast}+\sigma_{k}\Big\}\omega_{k},\\
y_{k}^{\ast}=(A_{k}+B_{k}\Theta_{k}^{\ast})'\mathbb{E}(y_{k+1}^{\ast}|\mathfrak{F}_{k-1})+[\bar{A}_{k}+B_{k}\bar{\Theta}_{k}^{\ast}+\bar{B}_{k}
(\Theta_{k}^{\ast}+\bar{\Theta}_{k}^{\ast})]'\mathbb{E}y_{k+1}^{\ast}\\
~~~~~~~~~~~+(C_{k}+D_{k}\Theta_{k}^{\ast})'\mathbb{E}(y_{k+1}^{\ast}\omega_{k}|\mathfrak{F}_{k-1})+[\bar{C}_{k}+D_{k}\bar{\Theta}_{k}^{\ast}
+\bar{D}_{k}(\Theta_{k}^{\ast}+\bar{\Theta}_{k}^{\ast})]'
\mathbb{E}(y_{k+1}^{\ast}\omega_{k})\\
~~~~~~~~~~~+\tilde{q}_{k}+\hat{q}_{k}+\tilde{Q}_{k}x_{k}^{\ast}+\hat{Q}_{k}\mathbb{E}x_{k}^{\ast}
+\tilde{S}_{k}u_{k}^{\ast}+\hat{S}_{k}\mathbb{E}u_{k}^{\ast},\\
x_{l}^{\ast}=\xi,~~y_{N}^{\ast}=Gx_{N}^{\ast}+\bar{G}\mathbb{E}x_{N}^{\ast}+g+\bar{g}.
\end{array}\right.
\end{eqnarray}
Then $u^{\ast}$ is a closed-loop optimal control of Problem (MF-LQ) if and only if
\begin{align}\label{2}
&\mathbb{E}\{\langle Gx_{N},x_{N}\rangle+\langle \bar{G}\mathbb{E}x_{N},\mathbb{E}x_{N}\rangle\}+\mathbb{E}\sum_{k=l}^{N-1}[\langle \tilde{Q}_{k}x_{k},x_{k}\rangle+2\langle \tilde{S}_{k}x_{k},u_{k}\rangle +\langle R_{k}u_{k},u_{k}\rangle\nonumber\\
&~~~~~~~~~~+\langle \hat{Q}_{k}\mathbb{E}x_{k},\mathbb{E}x_{k}\rangle+2\langle \hat{S}_{k}\mathbb{E}x_{k},\mathbb{E}u_{k}\rangle +\langle \bar{R}_{k}\mathbb{E}u_{k},\mathbb{E}u_{k}\rangle]\geq0,~~\forall~u_{k}~\in \mathscr{U}_{ad},
\end{align}
where $x_{k}$ satisfies
\begin{eqnarray}\label{3}
  \left\{
  \begin{array}{ll}
x_{k+1}=(A_{k}+B_{k}\Theta_{k}^{\ast})x_{k}+[\bar{A}_{k}+B_{k}\bar{\Theta}_{k}^{\ast}+\bar{B}_{k}(\Theta_{k}^{\ast}+\bar{\Theta}_{k}^{\ast})]\mathbb{E}x_{k}
+B_{k}u_{k}+\bar{B}_{k}\mathbb{E}u_{k}\\
~~~~~~~~~~~~+\Big\{(C_{k}+D_{k}\Theta_{k}^{\ast})x_{k}+[\bar{C}_{k}+D_{k}\bar{\Theta}_{k}^{\ast}+\bar{D}_{k}(\Theta_{k}^{\ast}+\bar{\Theta}_{k}^{\ast})]\mathbb{E}x_{k}
+D_{k}u_{k}+\bar{D}_{k}\mathbb{E}u_{k}\Big\}\omega_{k},\\
x_{l}=0,
\end{array}\right.
\end{eqnarray}
and the following stationarity condition holds
\begin{align}\label{4}
&B'_{k}\mathbb{E}(y_{k+1}^{\ast}|\mathfrak{F}_{k-1})+\bar{B}'_{k}\mathbb{E}y_{k+1}^{\ast}+\rho_{k}+\bar{\rho}_{k}+D'_{k}\mathbb{E}(y_{k+1}^{\ast}\omega_{k}|\mathfrak{F}_{k-1})\nonumber\\
&~~~~~~~~~~~~~~~~~~+\bar{D}'_{k}\mathbb{E}(y_{k+1}^{\ast}\omega_{k})
+\tilde{S}'_{k}x_{k}^{\ast}+\hat{S}'_{k}\mathbb{E}x_{k}^{\ast}+R_{k}u_{k}^{\ast}+\bar{R}_{k}\mathbb{E}u_{k}^{\ast}=0.
\end{align}
(ii) If $(\Theta^{\ast},\bar{\Theta}^{\ast},u^{\ast})$ is an optimal closed-loop strategy of Problem (MF-LQ), then $(\Theta^{\ast},\bar{\Theta}^{\ast},0)$ is an optimal closed-loop strategy of Problem (MF-LQ)$^{0}$.
\begin{proof} (i) By Theorem 4.1 and the definition of closed-loop solvability, the proof of (i) is clear.

(ii) For any $\xi \in \mathbb{R}^{n}$, \eqref{1}-\eqref{4} hold.
As $(\Theta^{\ast},\bar{\Theta}^{\ast},u^{\ast})$ is independent of $\xi$, by means of subtracting solutions corresponding $\xi$ and 0, the latter from the former, we obtain that for any $\xi \in \mathbb{R}^{n}$, the following MF-FBSDEs
\begin{eqnarray}\label{5}
  \left\{
  \begin{array}{ll}
x_{k+1}=(A_{k}+B_{k}\Theta_{k}^{\ast})x_{k}+[\bar{A}_{k}+B_{k}\bar{\Theta}_{k}^{\ast}+\bar{B}_{k}(\Theta_{k}^{\ast}+\bar{\Theta}_{k}^{\ast})]\mathbb{E}x_{k}\\
~~~~~~~~~~~~~~+\{(C_{k}+D_{k}\Theta_{k}^{\ast})x_{k}+[\bar{C}_{k}+D_{k}\bar{\Theta}_{k}^{\ast}+\bar{D}_{k}(\Theta_{k}^{\ast}+\bar{\Theta}_{k}^{\ast})]\mathbb{E}x_{k}\}\omega_{k},\\
y_{k}=(A_{k}+B_{k}\Theta_{k}^{\ast})'\mathbb{E}(y_{k+1}|\mathfrak{F}_{k-1})+[\bar{A}_{k}+B_{k}\bar{\Theta}_{k}^{\ast}+\bar{B}_{k}(\Theta_{k}^{\ast}+\bar{\Theta}_{k}^{\ast})]'\mathbb{E}y_{k+1}
+\tilde{Q}_{k}x_{k}+\hat{Q}_{k}\mathbb{E}x_{k}\\
~~~~~~~~~~~~~~+(C_{k}+D_{k}\Theta_{k}^{\ast})'\mathbb{E}(y_{k+1}\omega_{k}|\mathfrak{F}_{k-1})+[\bar{C}_{k}+D_{k}\bar{\Theta}_{k}^{\ast}+\bar{D}_{k}
(\Theta_{k}^{\ast}+\bar{\Theta}_{k}^{\ast})]'
\mathbb{E}(y_{k+1}\omega_{k}),\\
x_{l}=\xi,~~y_{N}=Gx_{N}+\bar{G}\mathbb{E}x_{N},
\end{array}\right.
\end{eqnarray}
also exists a solution $(x,y)$ satisfying
\begin{align*}
B'_{k}\mathbb{E}(y_{k+1}|\mathfrak{F}_{k-1})+\bar{B}'_{k}\mathbb{E}y_{k+1}+D'_{k}\mathbb{E}(y_{k+1}\omega_{k}|\mathfrak{F}_{k-1})
+\bar{D}'_{k}\mathbb{E}(y_{k+1}\omega_{k})
+\tilde{S}'_{k}x_{k}+\hat{S}'_{k}\mathbb{E}x_{k}=0.
\end{align*}
Using (i), it is obvious that $(\Theta^{\ast},\bar{\Theta}^{\ast},0)$ is an optimal closed-loop strategy of Problem (MF-LQ)$^{0}$.
\end{proof}

\vskip0mm\noindent
\textbf{Theorem 4.3.} Assume that Problem (MF-LQ) is closed-loop solvable, then the GRE (4.7) is regularly solvable.
\begin{proof} Assume that $(\Theta^{\ast},\bar{\Theta}^{\ast},u^{\ast})$ is an optimal closed-loop strategy of Problem (MF-LQ), then by Theorem 4.2, $(\Theta^{\ast},\bar{\Theta}^{\ast},0)$ is an optimal closed-loop strategy of Problem (MF-LQ)$^{0}$. For Problem (MF-LQ)$^{0}$, the closed-loop system becomes
\begin{eqnarray}
  \left\{
  \begin{array}{ll}
x_{k+1}=(A_{k}+B_{k}\Theta_{k})x_{k}+[\bar{A}_{k}+B_{k}\bar{\Theta}_{k}+\bar{B}_{k}(\Theta_{k}+\bar{\Theta}_{k})]\mathbb{E}x_{k}\\
~~~~~~~~~~~~+\{(C_{k}+D_{k}\Theta_{k})x_{k}+[\bar{C}_{k}+D_{k}\bar{\Theta}_{k}+\bar{D}_{k}(\Theta_{k}+\bar{\Theta}_{k})]\mathbb{E}x_{k}\}\omega_{k},\\
x_{l}=\xi.
\end{array}\right.
\end{eqnarray}
Set $X_{k}=\mathbb{E}(x_{k}x_{k}')$, $Y_{k}=\mathbb{E}x_{k}(\mathbb{E}x_{k})'$, we have
\begin{eqnarray}
  \left\{
  \begin{array}{ll}
X_{k+1}=[\bar{A}_{k}+\bar{B}_{k}\Theta_{k}+(B_{k}+\bar{B}_{k})\bar{\Theta}_{k}]Y_{k}(A_{k}+B_{k}\Theta_{k})'+(A_{k}+B_{k}\Theta_{k})X_{k}(A_{k}+B_{k}\Theta_{k})'\\
~~~~~~~~~~~~~+(A_{k}+B_{k}\Theta_{k})Y_{k}[\bar{A}_{k}+\bar{B}_{k}\Theta_{k}+(B_{k}+\bar{B}_{k})\bar{\Theta}_{k})]'+[\bar{A}_{k}+\bar{B}_{k}\Theta_{k}+(B_{k}+\bar{B}_{k})\bar{\Theta}_{k}]Y_{k}\\
~~~~~~~~~~~~~\times[\bar{A}_{k}+\bar{B}_{k}\Theta_{k}+(B_{k}+\bar{B}_{k})\bar{\Theta}_{k}]'+[\bar{C}_{k}+\bar{D}_{k}\Theta_{k}+(D_{k}+\bar{D}_{k})\bar{\Theta}_{k}]Y_{k}(C_{k}+D_{k}\Theta_{k})'\\
~~~~~~~~~~~~~+(C_{k}+D_{k}\Theta_{k})X_{k}(C_{k}+D_{k}\Theta_{k})'+(C_{k}+D_{k}\Theta_{k})Y_{k}[\bar{C}_{k}+\bar{D}_{k}\Theta_{k}+(D_{k}+\bar{D}_{k})\bar{\Theta}_{k})]'\\
~~~~~~~~~~~~~+[\bar{C}_{k}+\bar{D}_{k}\Theta_{k}+(D_{k}+\bar{D}_{k})\bar{\Theta}_{k}]Y_{k}[\bar{C}_{k}+\bar{D}_{k}\Theta_{k}+(D_{k}+\bar{D}_{k})\bar{\Theta}_{k}]'\\
~~~~~~~\doteq\mathcal{X}_{k}(\Theta_{k},\bar{\Theta}_{k},X_{k},Y_{k}),\\
Y_{k+1}=[(A_{k}+\bar{A}_{k})+(B_{k}+\bar{B}_{k})(\Theta_{k}+\bar{\Theta}_{k})]Y_{k}[(A_{k}+\bar{A}_{k})+(B_{k}+\bar{B}_{k})(\Theta_{k}+\bar{\Theta}_{k})]'\\
~~~~~~~\doteq\mathcal{Y}_{k}(\Theta_{k},\bar{\Theta}_{k},X_{k},Y_{k}),\\
X_{l}=\mathbb{E}(\xi\xi'),~~Y_{l}=\mathbb{E}\xi(\mathbb{E}\xi)'.
\end{array}\right.
\end{eqnarray}
The cost functional $J^{0}(l,\xi;\Theta_{k} x_{k}+\bar{\Theta}_{k}\mathbb{E}x_{k})$ can be expressed as
\begin{align}
J^{0}(l,\xi;\Theta_{k} x_{k}+\bar{\Theta}_{k}\mathbb{E}x_{k})=\sum\limits_{k=l}^{N-1}\Big\{Tr[(Q_{k}+\Theta'_{k}R_{k}\Theta_{k}+\Theta'_{k}S_{k}+S'_{k}\Theta_{k})X_{k}+\Phi_{k}Y_{k}]\Big\}
+Tr[GX_{N}+\bar{G}Y_{N}],
\end{align}
where
\begin{align*}
\Phi_{k}&=\bar{Q}_{k}+(\Theta_{k}+\bar{\Theta}_{k})'\bar{R}_{k}(\Theta_{k}+\bar{\Theta}_{k})+\Theta'_{k}R_{k}\bar{\Theta}_{k}
+\bar{\Theta}'_{k}R_{k}\Theta_{k}+\bar{\Theta}'_{k}R_{k}\bar{\Theta}_{k}\\
&~~~~~~~~~~~~~~+(\Theta_{k}+\bar{\Theta}_{k})'\bar{S}_{k}+\bar{S}_{k}'(\Theta_{k}+\bar{\Theta}_{k})+\bar{\Theta}_{k}'S_{k}+S_{k}'\bar{\Theta}_{k}.
\end{align*} Now, we introduce the Lagrangian function
$\mathfrak{L}=\sum_{k=l}^{N-1}\mathfrak{L}_{k}+Tr[GX_{N}+\bar{G}Y_{N}],$
where
\begin{align}\label{7}
\mathfrak{L}_{k}&=Tr[(Q_{k}+S'_{k}R_{k}S_{k}+\Theta_{k}'S_{k}+S'_{k}\Theta_{k})X_{k}+\Phi_{k}Y_{k}]
+Tr[P_{k+1}(\mathcal{X}_{k}(\Theta_{k},\bar{\Theta}_{k},X_{k},Y_{k})-X_{k+1})]\nonumber\\
&~~~~~~~~~~+Tr[\Pi_{k+1}(\mathcal{Y}_{k}(\Theta_{k},\bar{\Theta}_{k},X_{k},Y_{k})-Y_{k+1})].
\end{align}
Clearly, using the matrix minimum principle, the lagrangian multipliers $P, \bar{P}$ and the optimal feedback gains $(\Theta^{\ast},\bar{\Theta}^{\ast})$ satisfy the following first-order necessary conditions
\begin{eqnarray}\label{8}
  \left\{
  \begin{array}{ll}
P_{k}=\frac{\partial \mathfrak{L}_{k}}{\partial X_{k}}\big|_{\Theta_{k}=\Theta_{k}^{\ast}},~~\bar{P}_{k}=\frac{\partial \mathfrak{L}_{k}}{\partial Y_{k}}\big|_{\Theta_{k}=\Theta_{k}^{\ast},\bar{\Theta}_{k}=\bar{\Theta}_{k}^{\ast}},~~
\frac{\partial \mathfrak{L}_{k}}{\partial \Theta_{k}}=0,~~\frac{\partial \mathfrak{L}_{k}}{\partial \bar{\Theta}_{k}}=0,~~k\in \mathbb{N},\\
P_{N}=G,~~\bar{P}_{N}=\bar{G}.
\end{array}\right.
\end{eqnarray}
By some calculations,
\begin{eqnarray}\label{a}
  \left\{
  \begin{array}{ll}
P_{k}=(A_{k}+B_{k}\Theta_{k}^{\ast})'P_{k+1}(A_{k}+B_{k}\Theta_{k}^{\ast})+(C_{k}+D_{k}\Theta_{k}^{\ast})'P_{k+1}(C_{k}+D_{k}\Theta_{k}^{\ast})\\
~~~~~~~~~~~+Q_{k}+(\Theta_{k}^{\ast})'R_{k}\Theta_{k}^{\ast}+(\Theta_{k}^{\ast})'S_{k}+S'_{k}\Theta_{k}^{\ast},\\
\bar{P}_{k}=\bar{Q}_{k}+(\Theta_{k}^{\ast}+\bar{\Theta}_{k}^{\ast})'\bar{R}_{k}(\Theta_{k}^{\ast}+\bar{\Theta}_{k}^{\ast})+(\Theta_{k}^{\ast})'R_{k}\bar{\Theta}_{k}^{\ast}
+(\bar{\Theta}_{k})^{\ast}R_{k}\Theta_{k}^{\ast}+(\bar{\Theta}_{k}^{\ast})'R_{k}\bar{\Theta}_{k}^{\ast}+(\Theta_{k}^{\ast}+\bar{\Theta}_{k}^{\ast})'\bar{S}_{k}\\
~~~~~~~~~~~+\bar{S}_{k}'(\Theta_{k}^{\ast}+\bar{\Theta}_{k}^{\ast})+(\bar{\Theta}_{k}^{\ast})'S_{k}+S_{k}'\bar{\Theta}_{k}^{\ast}+(A_{k}+B_{k}\Theta_{k}^{\ast})'P_{k+1}
[\bar{A}_{k}+B_{k}\bar{\Theta}_{k}^{\ast}+\bar{B}_{k}(\Theta_{k}^{\ast}+\bar{\Theta}_{k}^{\ast})]\\
~~~~~~~~~~~+[\bar{A}_{k}+B_{k}\bar{\Theta}_{k}^{\ast}+\bar{B}_{k}(\Theta_{k}^{\ast}+\bar{\Theta}_{k}^{\ast})]'P_{k+1}(A_{k}+B_{k}\Theta_{k}^{\ast})
+[\bar{A}_{k}+B_{k}\bar{\Theta}_{k}^{\ast}+\bar{B}_{k}(\Theta_{k}^{\ast}+\bar{\Theta}_{k}^{\ast})]'\\
~~~~~~~~~~~\times P_{k+1}[\bar{A}_{k}+B_{k}\bar{\Theta}_{k}^{\ast}+\bar{B}_{k}(\Theta_{k}^{\ast}+\bar{\Theta}_{k}^{\ast})]
+(C_{k}+D_{k}\Theta_{k}^{\ast})'P_{k+1}[\bar{C}_{k}+D_{k}\bar{\Theta}_{k}^{\ast}+\bar{D}_{k}(\Theta_{k}^{\ast}+\bar{\Theta}_{k}^{\ast})]\\
~~~~~~~~~~~+[\bar{C}_{k}+D_{k}\bar{\Theta}_{k}^{\ast}+\bar{D}_{k}(\Theta_{k}^{\ast}+\bar{\Theta}_{k}^{\ast})]'P_{k+1}(C_{k}+D_{k}\Theta_{k}^{\ast})
+[\bar{C}_{k}+D_{k}\bar{\Theta}_{k}^{\ast}+\bar{C}_{k}(\Theta_{k}^{\ast}+\bar{\Theta}_{k}^{\ast})]'\\
~~~~~~~~~~~\times P_{k+1}[\bar{C}_{k}+D_{k}\bar{\Theta}_{k}^{\ast}+\bar{D}_{k}(\Theta_{k}^{\ast}+\bar{\Theta}_{k}^{\ast})]+[A_{k}+\bar{A}_{k}+(B_{k}+\bar{B}_{k})(\Theta_{k}^{\ast}+\bar{\Theta}_{k}^{\ast})]'\\
~~~~~~~~~~~\times\bar{P}_{k+1}[A_{k}+\bar{A}_{k}+(B_{k}+\bar{B}_{k})(\Theta_{k}^{\ast}+\bar{\Theta}_{k}^{\ast})].
\end{array}\right.
\end{eqnarray}
Set $\Pi_{k}=P_{k}+\bar{P}_{k},~\Lambda_{k}=\Theta_{k}^{\ast}+\bar{\Theta}_{k}^{\ast},$
according to \eqref{7}-\eqref{a}, we have
\begin{align}\label{10}
\Pi_{k}&=\Lambda'_{k}(R_{k}+\bar{R}_{k})\Lambda_{k}+Q_{k}+\bar{Q}_{k}+\Lambda'_{k}(S_{k}+\bar{S}_{k})+(S_{k}+\bar{S}_{k})'\Lambda_{k}+[C_{k}+\bar{C}_{k}+(D_{k}+\bar{D}_{k})\Lambda_{k}]'P_{k+1}\nonumber\\
&~~~~~\times[C_{k}+\bar{C}_{k}+(D_{k}+\bar{D}_{k})\Lambda_{k}]+[A_{k}+\bar{A}_{k}+(B_{k}+\bar{B}_{k})\Lambda_{k}]'\Pi_{k+1}[A_{k}+\bar{A}_{k}+(B_{k}+\bar{B}_{k})\Lambda_{k}],
\end{align}
and
\begin{align}\label{11}
0&=2\Big\{R_{k}\Theta_{k}+B'_{k}P_{k+1}(A_{k}+B_{k}\Theta_{k})+D'_{k}P_{k+1}(C_{k}+D_{k}\Theta_{k})+S_{k}'\Big\}X_{k}+2\Big\{[R_{k}+\bar{R}_{k}\nonumber\\
&~~~~~~~+(B_{k}+\bar{B}_{k})'\Pi_{k+1}(B_{k}+\bar{B}_{k})+(D_{k}+\bar{D}_{k})'
P_{k+1}(D_{k}+\bar{D}_{k})]\bar{\Theta}_{k}-B'_{k}P_{k+1}A_{k}+[\bar{R}_{k}\nonumber\\
&~~~~~~~+(B_{k}+\bar{B}_{k})'
\Pi_{k+1}(B_{k}+\bar{B}_{k})-B'_{k}P_{k+1}B_{k}+\bar{D}'_{k}P_{k+1}D_{k}+D'_{k}P_{k+1}\bar{D}_{k}+\bar{D}'_{k}P_{k+1}\bar{D}_{k}]\Theta_{k}\nonumber\\
&~~~~~~~+(B_{k}+\bar{B}_{k})'\Pi_{k+1}(A_{k}+\bar{A}_{k})+(D_{k}+\bar{D}_{k})'
P_{k+1}(C_{k}+\bar{C}_{k})-D'_{k}P_{k+1}C_{k}+\bar{S}_{k}'\Big\}Y_{k}\nonumber\\
&=2\Big\{(R_{k}+B_{k}'P_{k+1}B_{k}+D_{k}'P_{k+1}D_{k})\Theta_{k}+B_{k}'P_{k+1}A_{k}+D_{k}'P_{k+1}C_{k}+S_{k}'\Big\}X_{k}\nonumber\\
&~~~~~~~+2\Big\{[R_{k}+\bar{R}_{k}+(B_{k}+\bar{B}_{k})'\Pi_{k+1}(B_{k}+\bar{B}_{k})+(D_{k}+\bar{D}_{k})'P_{k+1}(D_{k}+\bar{D}_{k})]\bar{\Theta}_{k}\nonumber\\
&~~~~~~~+[\bar{R}_{k}+(B_{k}+\bar{B}_{k})'\Pi_{k+1}(B_{k}+\bar{B}_{k})-B_{k}'P_{k+1}B_{k}+\bar{D}_{k}'P_{k+1}D_{k}+D_{k}'P_{k+1}\bar{D}_{k}
+\bar{D}_{k}'P_{k+1}\bar{D}_{k}]\Theta_{k}\nonumber\\
&~~~~~~~+(B_{k}+\bar{B}_{k})'\Pi_{k+1}(A_{k}+\bar{A}_{k})+(D_{k}+\bar{D}_{k})'P_{k+1}(C_{k}+\bar{C}_{k})-B_{k}'P_{k+1}A_{k}-D_{k}'P_{k+1}C_{k}\Big\}Y_{k}\nonumber\\
&=2\Big\{(R_{k}+B_{k}'P_{k+1}B_{k}+D_{k}'P_{k+1}D_{k})\Theta_{k}+B_{k}'P_{k+1}A_{k}+D_{k}'P_{k+1}C_{k}+S_{k}'\Big\}(X_{k}-Y_{k})\nonumber\\
&~~~~~~~+2\Big\{[R_{k}+\bar{R}_{k}+(B_{k}+\bar{B}_{k})'\Pi_{k+1}(B_{k}+\bar{B}_{k})+(D_{k}+\bar{D}_{k})'P_{k+1}(D_{k}+\bar{D}_{k})]\Lambda_{k}\nonumber\\
&~~~~~~~+(B_{k}+\bar{B}_{k})'\Pi_{k+1}(A_{k}+\bar{A}_{k})+(D_{k}+\bar{D}_{k})'P_{k+1}(C_{k}+\bar{C}_{k})+S_{k}'+\bar{S}_{k}'\Big\}Y_{k},
\end{align}
besides,
\begin{align}\label{12}
0&=2\Big\{[R_{k}+\bar{R}_{k}+(B_{k}+\bar{B}_{k})'
\Pi_{k+1}(B_{k}+\bar{B}_{k})+(D_{k}+\bar{D}_{k})'
P_{k+1}(D_{k}+\bar{D}_{k})]\Lambda_{k}\nonumber\\
&~~~~~+(B_{k}+\bar{B}_{k})'
\Pi_{k+1}(A_{k}+\bar{A}_{k})+(D_{k}+\bar{D}_{k})'P_{k+1}(C_{k}+\bar{C}_{k})+S_{k}'+\bar{S}_{k}'\Big\}Y_{k}.
\end{align}
Notice that \eqref{11}-\eqref{12} hold for all initial values $Y_{l}=\mathbb{E}(\xi\xi')$ and $Y_{l}=\mathbb{E}\xi(\mathbb{E}\xi)'$, respectively, then it could be reduced to
\begin{align*}
&\Upsilon_{k}\Theta_{k}^{\ast}+B_{k}'P_{k+1}A_{k}+D_{k}'P_{k+1}C_{k}+S_{k}'=0,\\
&\bar{\Upsilon}_{k}\Lambda_{k}+(B_{k}+\bar{B}_{k})'
\Pi_{k+1}(A_{k}+\bar{A}_{k})+(D_{k}+\bar{D}_{k})'P_{k+1}(C_{k}+\bar{C}_{k})+S_{k}'+\bar{S}_{k}'=0.
\end{align*}
These further indicate that
\begin{align*}
&\mathcal{R}(B_{k}'P_{k+1}A_{k}+D_{k}'P_{k+1}C_{k}+S_{k}')\subseteq \mathcal{R}(\Upsilon_{k}),\\
&\mathcal{R}[(B_{k}+\bar{B}_{k})'
\Pi_{k+1}(A_{k}+\bar{A}_{k})+(D_{k}+\bar{D}_{k})'P_{k+1}(C_{k}+\bar{C}_{k})+S_{k}'+\bar{S}_{k}']\subseteq \mathcal{R}(\bar{\Upsilon}_{k}).
\end{align*}
Moreover, $\Upsilon^{\dagger}_{k}\Upsilon_{k}$ and $\bar{\Upsilon}^{\dagger}_{k}\bar{\Upsilon}_{k}$ are orthogonal projections, it yields that
\begin{align*}
&\Upsilon^{\dagger}_{k}(B_{k}'P_{k+1}A_{k}+D_{k}'P_{k+1}C_{k}+S_{k}')\in L^{2}(\mathbb{N};\mathbb{R}^{m\times n}),\\
&\bar{\Upsilon}^{\dagger}_{k}[(B_{k}+\bar{B}_{k})'\Pi_{k+1}(A_{k}+\bar{A}_{k})+(D_{k}+\bar{D}_{k})'P_{k+1}(C_{k}+\bar{C}_{k})+S_{k}'+\bar{S}_{k}']\in L^{2}(\mathbb{N};\mathbb{R}^{m\times n}),
\end{align*}
and for some $\phi, \varphi\in L^{2}(\mathbb{N};\mathbb{R}^{m\times n})$,
\begin{align}
&\Theta_{k}^{\ast}=-\Upsilon_{k}^{\dagger}(B_{k}'P_{k+1}A_{k}+D_{k}'P_{k+1}C_{k}+S_{k}')+(I-\Upsilon_{k}^{\dagger}\Upsilon_{k})\phi,\\
&\Lambda_{k}=-\bar{\Upsilon}_{k}^{\dagger}[(B_{k}+\bar{B}_{k})'\Pi_{k+1}(A_{k}+\bar{A}_{k})+(D_{k}+\bar{D}_{k})'P_{k+1}(C_{k}+\bar{C}_{k})+S_{k}'+\bar{S}_{k}']+
(I-\Upsilon_{k}^{\dagger}\Upsilon_{k})\varphi.
\end{align}
By substituting (4.27)-(4.28) to \eqref{a}-\eqref{10}, we can immediately get that $(P,\Pi)$ satisfies GRE (4.7).
Next, it remains to show $\Upsilon_{k}\geq0$, $\bar{\Upsilon}_{k}\geq0$. To this end, we consider the state equation
\begin{eqnarray}\label{14}
  \left\{
  \begin{array}{ll}
x_{k+1}=A_{k}x_{k}+\bar{A}_{k}\mathbb{E}x_{k}+B_{k}u_{k}+\bar{B}_{k}\mathbb{E}u_{k}
+(C_{k}x_{k}+\bar{C}_{k}\mathbb{E}x_{k}+D_{k}u_{k}+\bar{D}_{k}\mathbb{E}u_{k})\omega_{k},\\
x_{l}=0.
\end{array}\right.
\end{eqnarray}
Define the Lyapunov function candidate as
\begin{align}\label{15}
\bar{V}_{N}(k,x_{k})=\mathbb{E}[(x_{k}-\mathbb{E}x_{k})'P_{k}(x_{k}-\mathbb{E}x_{k})]+(\mathbb{E}x_{k})'\Pi_{k}\mathbb{E}x_{k},
\end{align}
it follows that
\begin{align*}
&\bar{V}_{N}(k,x_{k})-\bar{V}_{N}(k+1,x_{k+1})\\
&~~=\mathbb{E}\Big\{(x_{k}-\mathbb{E}x_{k})'P_{k}(x_{k}-\mathbb{E}x_{k})+(\mathbb{E}x_{k})'\Pi_{k}\mathbb{E}x_{k}-x_{k}'(A'_{k}P_{k+1}A_{k}+C'_{k}P_{k+1}C_{k})x_{k}
-u_{k}'(B_{k}'P_{k+1}B_{k}\\
&~~~~~~~+D_{k}'P_{k+1}D_{k})u_{k}-x_{k}'(A_{k}'P_{k+1}B_{k}+C_{k}'P_{k+1}D_{k})u_{k}-u_{k}'(B_{k}'P_{k+1}A_{k}+D_{k}'P_{k+1}C_{k})x_{k}\\
&~~~~~~~-(\mathbb{E}x_{k})'[(A_{k}+\bar{A}_{k})'\Pi_{k+1}(A_{k}+\bar{A}_{k})+(C_{k}+\bar{C}_{k})'P_{k+1}(C_{k}+\bar{C}_{k})-A'_{k}P_{k+1}A_{k}-C'_{k}P_{k+1}C_{k}]\mathbb{E}x_{k}\\
&~~~~~~~-(\mathbb{E}x_{k})'[C_{k}'P_{k+1}\bar{D}_{k}
+\bar{C}_{k}'P_{k+1}D_{k}+\bar{C}_{k}'P_{k+1}\bar{D}_{k}+(A_{k}+\bar{A}_{k})'\Pi_{k+1}(B_{k}+\bar{B}_{k})-A_{k}'P_{k+1}B_{k}]\mathbb{E}u_{k}\\
&~~~~~~~-(\mathbb{E}u_{k})'[D_{k}'P_{k+1}\bar{C}_{k}+\bar{D}_{k}'P_{k+1}C_{k}+\bar{D}_{k}'P_{k+1}\bar{C}_{k}
+(B_{k}+\bar{B}_{k})'\Pi_{k+1}(A_{k}+\bar{A}_{k})-B_{k}'P_{k+1}A_{k}]\mathbb{E}x_{k}\\
&~~~~~~~-(\mathbb{E}u_{k})'[D_{k}'P_{k+1}\bar{D}_{k}+\bar{D}_{k}'P_{k+1}D_{k}
+\bar{D}_{k}'P_{k+1}\bar{D}_{k}
+(B_{k}+\bar{B}_{k})'\Pi_{k+1}(B_{k}+\bar{B}_{k})-B_{k}'P_{k+1}B_{k}]\mathbb{E}u_{k}\Big\}\\
&~~=\mathbb{E}\bigg\{x_{k}'[Q_{k}-(B_{k}'P_{k+1}A_{k}+D_{k}'P_{k+1}C_{k}+S'_{k})'\Upsilon_{k}^{\dagger}(B_{k}'P_{k+1}A_{k}+D_{k}'P_{k+1}C_{k}+S'_{k})]x_{k}\\
&~~~~~~~+(\mathbb{E}x_{k})'\Big\{\bar{Q}_{k}+(B_{k}'P_{k+1}A_{k}+D_{k}'P_{k+1}C_{k}+S'_{k})'\Upsilon_{k}^{\dagger}(B_{k}'P_{k+1}A_{k}+D_{k}'P_{k+1}C_{k}+S'_{k})\\
&~~~~~~~-[(B_{k}+\bar{B}_{k})'\Pi_{k+1}(A_{k}+\bar{A}_{k})+(D_{k}+\bar{D}_{k})'P_{k+1}(C_{k}+\bar{C}_{k})+S'_{k}+\bar{S}'_{k}]'\bar{\Upsilon}_{k}^{\dagger}\\
&~~~~~~~\times[(B_{k}+\bar{B}_{k})'\Pi_{k+1}(A_{k}+\bar{A}_{k})+(D_{k}+\bar{D}_{k})'P_{k+1}(C_{k}+\bar{C}_{k})+S'_{k}+\bar{S}'_{k}]\Big\}\mathbb{E}x_{k}\\
&~~~~~~~-x_{k}' (B_{k}'P_{k+1}A_{k}+D_{k}'P_{k+1}C_{k}+S'_{k})'u_{k}-u_{k}'(B_{k}'P_{k+1}A_{k}+D_{k}'P_{k+1}C_{k}+S'_{k})x_{k}\\
&~~~~~~~-(\mathbb{E}x_{k})'\Big\{[(B_{k}+\bar{B}_{k})'\Pi_{k+1}(A_{k}+\bar{A}_{k})+(D_{k}+\bar{D}_{k})'P_{k+1}
(C_{k}+\bar{C}_{k})+S'_{k}+\bar{S}'_{k}]'\\
&~~~~~~~-(B_{k}'P_{k+1}A_{k}+D_{k}'P_{k+1}C_{k}+S'_{k})'\Big\}\mathbb{E}u_{k}-(\mathbb{E}u_{k})'
\Big\{[(B_{k}+\bar{B}_{k})'\Pi_{k+1}(A_{k}+\bar{A}_{k})\\
&~~~~~~~+(D_{k}+\bar{D}_{k})'P_{k+1}(C_{k}+\bar{C}_{k})+S'_{k}+\bar{S}'_{k}]-
[B_{k}'P_{k+1}A_{k}+D_{k}'P_{k+1}C_{k}+S'_{k}]\Big\}\mathbb{E}x_{k}\\
&~~~~~~~-u_{k}'\Upsilon_{k}u_{k}-(\mathbb{E}u_{k})'(\bar{\Upsilon}_{k}-\Upsilon_{k})\mathbb{E}u_{k}+u_{k}'R_{k}u_{k}+(\mathbb{E}u_{k})'\bar{R}_{k}\mathbb{E}u_{k}\bigg\}\\
&~~=\mathbb{E}\Big\{x_{k}'Q_{k}x_{k}+(\mathbb{E}x_{k})'Q_{k}\mathbb{E}x_{k}+u_{k}'R_{k}u_{k}
+(\mathbb{E}u_{k})'R_{k}\mathbb{E}u_{k}+2x_{k}'S_{k}u_{k}+2(\mathbb{E}x_{k})'\bar{S}_{k}\mathbb{E}u_{k}\Big\}\\
&~~~~~~~-\mathbb{E}\Big\{[u_{k}-\mathbb{E}u_{k}-\Theta_{k}(x_{k}-\mathbb{E}x_{k})]'\Upsilon_{k}[u_{k}-\mathbb{E}u_{k}-\Theta_{k}(x_{k}-\mathbb{E}x_{k})]\Big\}\\
&~~~~~~~-\Big\{\mathbb{E}u_{k}-(\Theta_{k}+\bar{\Theta}_{k})\mathbb{E}x_{k}]'\bar{\Upsilon}_{k}[\mathbb{E}u_{k}-(\Theta_{k}+\bar{\Theta}_{k})\mathbb{E}x_{k}\Big\}.
\end{align*}
Adding from $k=l$ to $k=N$ on both sides of the above equality, we see
\begin{align*}
J^{0}(l,0;u_{k})&=\mathbb{E}\sum\limits_{k=l}^{N-1}\Big\{[u_{k}-\mathbb{E}u_{k}-\Theta_{k}(x_{k}-\mathbb{E}x_{k})]'\Upsilon_{k}[u_{k}-\mathbb{E}u_{k}-\Theta_{k}(x_{k}-\mathbb{E}x_{k})]\\
&~~~~~~~~~~~~+[\mathbb{E}u_{k}-(\Theta_{k}+\bar{\Theta}_{k})\mathbb{E}x_{k}]'\bar{\Upsilon}_{k}
[\mathbb{E}u_{k}-(\Theta_{k}+\bar{\Theta}_{k})\mathbb{E}x_{k}]\Big\}.
\end{align*}
Observe that $(\Theta^{\ast},\bar{\Theta}^{\ast},0)$ is an optimal closed-loop strategy of Problem (MF-LQ)$^{0}$, one has
\begin{align*}
&\mathbb{E}\sum\limits_{k=l}^{N-1}\Big\{[u_{k}-\mathbb{E}u_{k}-\Theta_{k}(x_{k}-\mathbb{E}x_{k})]'\Upsilon_{k}[u_{k}-\mathbb{E}u_{k}-\Theta_{k}(x_{k}-\mathbb{E}x_{k})]\Big\}\\
&~~~~~~~~+\mathbb{E}\sum\limits_{k=l}^{N-1}\Big\{[\mathbb{E}u_{k}-(\Theta_{k}+\bar{\Theta}_{k})\mathbb{E}x_{k}]'\bar{\Upsilon}_{k}
[\mathbb{E}u_{k}-(\Theta_{k}+\bar{\Theta}_{k})\mathbb{E}x_{k}]\Big\}\\
&~~~~~~~~~~~~~~=J^{0}(l,0;u_{k})\geq J^{0}(l,0;\Theta^{\ast}_{k}x_{k}+\bar{\Theta}^{\ast}_{k}\mathbb{E}x_{k})=0,~~~~~\forall~u_{k}\in \mathscr{U}_{ad}.
\end{align*}
Since for any $u_{k}\in \mathscr{U}_{ad}$,
$$u_{k}=\Theta_{k}^{\ast} x_{k}+v_{k}\omega_{k},~~~~v_{k}\in L^{2}_{\mathfrak{F}}(\mathbb{N};\mathbb{R}^{m}),$$
the solution $x_{k}$ of \eqref{14} satisfies $\mathbb{E}x_{k}=0$, and then $\mathbb{E}u_{k}=0$. Furthermore, we could deduce that there exists a constant $c>0$ such that
\begin{align*}
0&\leq\mathbb{E}\sum\limits_{k=l}^{N-1}\langle \Upsilon_{k}[u_{k}-\Theta_{k}x_{k}],u_{k}-\Theta_{k}x_{k}\rangle
=\mathbb{E}\sum\limits_{k=l}^{N-1}\langle \Upsilon_{k}v_{k}\omega_{k},v_{k}\omega_{k}\rangle
\leq c\mathbb{E}\sum\limits_{k=l}^{N-1}\langle \Upsilon_{k}v_{k},v_{k}\rangle,~~\forall~v_{k}\in L^{2}_{\mathfrak{F}}(\mathbb{N};\mathbb{R}^{m}),
\end{align*}
which indicates $\Upsilon_{k}\geq0$. Similarly, for any $u_{k}\in \mathscr{U}_{ad}$ of the representation
$$u_{k}=\Theta_{k}^{\ast}(x_{k}-\mathbb{E}x_{k})+(\Theta_{k}^{\ast}+\bar{\Theta}_{k}^{\ast})\mathbb{E}x_{k}+v_{k},~~v_{k}\in L^{2}_{\mathfrak{F}}(\mathbb{N};\mathbb{R}^{m}),$$
the solution $x_{k}$ of \eqref{14} satisfies
$u_{k}-\mathbb{E}u_{k}=\Theta^{\ast}_{k}(x_{k}-\mathbb{E}x_{k}),~\mathbb{E}u_{k}-(\Theta_{k}^{\ast}+\bar{\Theta}_{k}^{\ast})\mathbb{E}x_{k}=v_{k}.$
Hence, we have
\begin{align*}
0\leq\mathbb{E}\sum\limits_{k=l}^{N-1}\langle \bar{\Upsilon}_{k}[\mathbb{E}u_{k}-(\Theta_{k}^{\ast}+\bar{\Theta}_{k}^{\ast})\mathbb{E}x_{k}],\mathbb{E}u_{k}-(\Theta_{k}^{\ast}+\bar{\Theta}_{k}^{\ast})\mathbb{E}x_{k}\rangle
=\mathbb{E}\sum\limits_{k=l}^{N-1}\langle \bar{\Upsilon}_{k}v_{k},v_{k}\rangle,~~\forall~v_{k}\in L^{2}_{\mathfrak{F}}(\mathbb{N};\mathbb{R}^{m}),
\end{align*}
which shows $\bar{\Upsilon}_{k}\geq0$. The proof is completed.
\end{proof}

\vskip0mm\noindent
\textbf{Theorem 4.4.} Problem (MF-LQ) is closed-loop solvable if and only if GRE (4.7) exists a pair regular solution $(P,\Pi)$,
and the following LRE exists a solution
\begin{eqnarray}\label{16}
  \left\{
  \begin{array}{ll}
\eta_{k}=(C_{k}+\bar{C}_{k})'P_{k+1}\sigma_{k}+(A_{k}+\bar{A}_{k})'(\Pi_{k+1}b_{k}+\eta_{k+1})-[(B_{k}+\bar{B}_{k})'\Pi_{k+1}(A_{k}+\bar{A}_{k})\\
~~~~~~~~~~+(D_{k}+\bar{D}_{k})'P_{k+1}(C_{k}+\bar{C}_{k})+S'_{k}+\bar{S}'_{k}]'\bar{\Upsilon}_{k}^{\dagger}\zeta_{k}+q_{k}+\bar{q}_{k},\\
\eta_{N}=g+\bar{g},
\end{array}\right.
\end{eqnarray}
besides,
\begin{align}\label{17}
\zeta_{k}\in \mathcal{R}(\bar{\Upsilon}_{k}),~~~~~~
-\bar{\Upsilon}_{k}^{\dagger}\zeta_{k}\in L^{2}_{\mathfrak{F}}(\mathbb{N};\mathbb{R}^{m}),
\end{align}
where
\begin{align*}
\zeta_{k}=(D_{k}+\bar{D}_{k})'P_{k+1}\sigma_{k}+(B_{k}+\bar{B}_{k})'(\Pi_{k+1}b_{k}+\eta_{k+1})+\rho_{k}+\bar{\rho}_{k}.
\end{align*}
In this case, the optimal closed-loop strategy $\{(\Theta_{k}^{\ast},\bar{\Theta}_{k}^{\ast},u_{k}^{\ast}),~k\in \mathbb{N}\}$ admits the following form
\begin{eqnarray}\label{18}
  \left\{
  \begin{array}{ll}
u_{k}^{\ast}=\Theta_{k}(x_{k}-\mathbb{E}x_{k})+\bar{\Theta}_{k}\mathbb{E}x_{k}-\Upsilon_{k}^{\dagger}\zeta_{k}, \\
\Theta_{k}=-\Upsilon_{k}^{\dagger}(B_{k}'P_{k+1}A_{k}+D_{k}'P_{k+1}C_{k}+S_{k}'),\\
\bar{\Theta}_{k}=-\bar{\Upsilon}_{k}^{\dagger}[(B_{k}+\bar{B}_{k})'\Pi_{k+1}(A_{k}+\bar{A}_{k})+(D_{k}+\bar{D})'P_{k+1}(C_{k}+\bar{C})+S_{k}'+\bar{S}_{k}'],
\end{array}\right.
\end{eqnarray}
and the value function
\begin{align}\label{19}
V(l,\xi)&=\sum\limits_{k=l}^{N-1}[2\eta_{k+1}'b_{k}+b_{k}'\Pi_{k+1}b_{k}+\sigma_{k}'P_{k+1}\sigma_{k}-\zeta_{k}'\bar{\Upsilon}_{k}^{\dagger}\zeta_{k}]\nonumber\\
&~~~~~~~~~+\mathbb{E}[(x_{l}-\mathbb{E}x_{l})'P_{l}(x_{l}-\mathbb{E}x_{l})]+(\mathbb{E}x_{l})'\Pi_{l} \mathbb{E}x_{l}+2\eta_{l}'\mathbb{E}x_{l}.
\end{align}
\begin{proof} Necessity. Assume that $(\Theta^{\ast},\bar{\Theta}^{\ast},u^{\ast})$ is an optimal closed-loop strategy of Problem (MF-LQ). Using Theorem 4.3, the GRE (4.7) is regularly solvable. We prove this by induction. For $l=N-1$, by some calculations, one has
\begin{align*}
V(N-1,\xi)&=\inf_{u_{N-1}}\mathbb{E}[x_{N-1}'Q_{N-1}x_{N-1}+(\mathbb{E}x_{N-1})'\bar{Q}_{N-1}\mathbb{E}x_{N-1}+u_{N-1}'R_{N-1}u_{N-1}+x_{N}'Gx_{N}\\
&~~~~~~~~+(\mathbb{E}u_{N-1})'\bar{R}_{N-1}\mathbb{E}u_{N-1}
+(\mathbb{E}x_{N})'\bar{G}\mathbb{E}x_{N}+2q_{N-1}'x_{N-1}+2\bar{q}_{N-1}'\mathbb{E}x_{N-1}\\
&~~~~~~~~+2\rho_{N-1}' u_{N-1}+2\bar{\rho}_{N-1}'\mathbb{E}u_{N-1}+2g'x_{N}+2\bar{g}'\mathbb{E}x_{N}]\\
&=\inf_{u_{N-1}}\mathbb{E}\Big\{(x_{N-1}-\mathbb{E}x_{N-1})'(Q_{N-1}+A_{N-1}'P_{N}A_{N-1}+C_{N-1}'P_{N}C_{N-1})(x_{N-1}-\mathbb{E}x_{N-1})\\
&~~~~~~~~+2(x_{N-1}-\mathbb{E}x_{N-1})'(A_{N-1}P_{N}B_{N-1}+C_{N-1}'P_{N}D_{N-1})(u_{N-1}-\mathbb{E}u_{N-1})+\sigma_{N-1} P_{N}\sigma_{N-1}\\
&~~~~~~~~+(u_{N-1}-\mathbb{E}u_{N-1})'
(R_{N-1}+B_{N-1}'P_{N}B_{N-1}+D_{N-1}'P_{N}D_{N-1})(u_{N-1}-\mathbb{E}u_{N-1})+2\eta_{N}' b_{N-1}\\
&~~~~~~~~+(\mathbb{E}x_{N-1})'[Q_{N-1}+\bar{Q}_{N-1}+(C_{N-1}+\bar{C}_{N-1})'P_{N}(C_{N-1}+\bar{C}_{N-1})+(A_{N-1}+\bar{A}_{N-1})'\Pi_{N} \\
&~~~~~~~~\times(A_{N-1}+\bar{A}_{N-1})]\mathbb{E}x_{N-1}
+2(\mathbb{E}x_{N-1})'[(C_{N-1}+\bar{C}_{N-1})'P_{N}(D_{N-1}+\bar{D}_{N-1})+(A_{N-1}\\
&~~~~~~~~+\bar{A}_{N-1})\Pi_{N} (B_{N-1}+\bar{B}_{N-1})]\mathbb{E}u_{N-1}
+(\mathbb{E}u_{N-1})'[R_{N-1}+\bar{R}_{N-1}+(D_{N-1}+\bar{D}_{N-1})'\\
&~~~~~~~~\times P_{N}(D_{N-1}+\bar{D}_{N-1})+(B_{N-1}+\bar{B}_{N-1})\Pi_{N}(B_{N-1}+\bar{B}_{N-1})]\mathbb{E}u_{N-1}+b_{N-1}\Pi_{N} b_{N-1}\\
&~~~~~~~~+2(\mathbb{E}u_{N-1})'[(D_{N-1}+\bar{D}_{N-1})P_{N}\sigma_{N-1}+(B_{N-1}+\bar{B}_{N-1})(\Pi_{N} b_{N-1}+\eta_{N})+\rho_{N-1}+\bar{\rho}_{N-1}]\\
&~~~~~~~~+2(\mathbb{E}x_{N-1})'[(C_{N-1}+\bar{C}_{N-1})P_{N}\sigma_{N-1}+(A_{N-1}+\bar{A}_{N-1})(\Pi_{N} b_{N-1}+\eta_{N})+q_{N-1}+\bar{q}_{N-1}]\Big\}\\
&=\inf_{u_{N-1}}\mathbb{E}\Big\{2(u_{N-1}-\mathbb{E}u_{N-1})'[(B_{N-1}'P_{N}A_{N-1}+D_{N-1}'P_{N}C_{N-1}+S'_{N-1})](x_{N-1}-\mathbb{E}x_{N-1})\\
&~~~~~~~~+(u_{N-1}-\mathbb{E}u_{N-1})'\Upsilon_{N-1}(u_{N-1}-\mathbb{E}u_{N-1})
+2(\mathbb{E}u_{N-1})'[(B_{N-1}+\bar{B}_{N-1})'\Pi_{N}\\
&~~~~~~~~\times(A_{N-1}+\bar{A}_{N-1})+(D_{N-1}+\bar{D}_{N-1})'P_{N}(C_{N-1}+\bar{C}_{N-1})+S'_{N-1}+\bar{S}'_{N-1}]\mathbb{E}x_{N-1}\\
&~~~~~~~~+(\mathbb{E}u_{N-1})'\bar{\Upsilon}_{N-1}\mathbb{E}u_{N-1}+2(\mathbb{E}u_{N-1})'\eta_{N-1}+
(x_{N-1}-\mathbb{E}x_{N-1})'(A_{N-1}'P_{N}A_{N-1}\\
&~~~~~~~~+C_{N-1}'P_{N}C_{N-1}+Q_{N-1})(x_{N-1}-\mathbb{E}x_{N-1})
+(\mathbb{E}x_{N-1})'[(C_{N-1}+\bar{C}_{N-1})'P_{N}\\
&~~~~~~~~\times(C_{N-1}+\bar{C}_{N-1})+(A_{N-1}+\bar{A}_{N-1})\Pi_{N} (A_{N-1}+\bar{A}_{N-1})+Q_{N-1}+\bar{Q}_{N-1}]\mathbb{E}x_{N-1}\\
&~~~~~~~~+2(\mathbb{E}x_{N-1})'[(A_{N-1}+\bar{A}_{N-1})'(\Pi_{N} b_{N-1}+\eta_{N})+(C_{N-1}+\bar{C}_{N-1})'P_{N}\sigma_{N-1}]\\
&~~~~~~~~+b_{N-1}'\Pi_{N} b_{N-1}+\sigma_{N-1}'\Pi_{N} \sigma_{N-1}+2\eta_{N} b_{N-1}+q_{N-1}+\bar{q}_{N-1}\Big\}
>-\infty,
\end{align*}
in which $x_{N-1}=\xi$. Based on Lemma 2.2, we see $\bar{\Upsilon}_{N-1}\bar{\Upsilon}^{\dagger}_{N-1}\eta_{N-1}=\eta_{N-1}$,
which shows that LRE \eqref{16} is solvable at $k=N-1$. Besides, the optimal control
\begin{eqnarray*}
u_{N-1}^{\ast}=\Theta_{N-1}(x_{N-1}-\mathbb{E}x_{N-1})+\bar{\Theta}_{N-1}\mathbb{E}x_{N-1}+v_{N-1},
\end{eqnarray*}
and
\begin{align*}
V(N-1,\xi)&=[2\eta_{N}'b_{N-1}+b_{N-1}'\Pi_{N}b_{N-1}+\sigma_{N-1}'P_{N}\sigma_{N-1}+\zeta_{N-1}'v_{N-1}]+2\eta_{N-1}'\mathbb{E}x_{N-1}\\
&~~~~~+\mathbb{E}[(x_{N-1}-\mathbb{E}x_{N-1})'P_{N}(x_{N-1}-\mathbb{E}x_{N-1})]+(\mathbb{E}x_{N-1})'\Pi_{N} \mathbb{E}x_{N-1}.
\end{align*}
For $l=N-2$, noting that $V(N-1,x_{N-1})=J(N-1,x_{N-1};u_{N-1}^{\ast})$ and the possible time-inconsistency of optimal control, one has
\begin{align*}
V(N-2,\xi)&=\inf_{u_{N-2},u_{N-1}}\mathbb{E}\sum\limits_{k=N-2}^{N-1}[x_{k}'Q_{k}x_{k}+(\mathbb{E}x_{k})'\bar{Q}_{k}\mathbb{E}x_{k}+u_{k}'R_{k}u_{k}
+(\mathbb{E}u_{k})'\bar{R}_{k}\mathbb{E}u_{k}+2q_{k}'x_{k}+2\rho_{k}'u_{k}\\
&~~~~~~~~~+2\bar{q}_{k}'\mathbb{E}x_{k}+2\bar{\rho}_{k}'\mathbb{E}u_{k}+x_{N}'Gx_{N}+(\mathbb{E}x_{N})'\bar{G}\mathbb{E}x_{N}+2g'x_{N}+2\bar{g}'\mathbb{E}x_{N}]\\
&\leq \inf_{u_{N-2}}\mathbb{E}[x_{N-2}'Q_{N-2}x_{N-2}+(\mathbb{E}x_{N-2})'\bar{Q}_{N-2}\mathbb{E}x_{N-2}+u_{N-2}'R_{N-2}u_{N-2}
+(\mathbb{E}u_{N-2})'\bar{R}_{N-2}\mathbb{E}u_{N-2}\\
&~~~~~~~~~+2q_{N-2}'x_{N-2}+2\rho_{N-2}'u_{N-2}+2\bar{q}_{N-2}'\mathbb{E}x_{N-2}+2\bar{\rho}_{N-2}'\mathbb{E}u_{N-2}+2g'x_{N}+2\bar{g}'\mathbb{E}x_{N}\\
&~~~~~~~~~+b_{N-1}'\Pi_{N}b_{N-1}+\sigma_{N-1}'P_{N}\sigma_{N-1}
+(x_{N-1}-\mathbb{E}x_{N-1})'P_{N}(x_{N-1}-\mathbb{E}x_{N-1})\\
&~~~~~~~~~+2\eta_{N}'b_{N-1}+\zeta_{N-1}'v_{N-1}+(\mathbb{E}x_{N-1})'\Pi_{N} \mathbb{E}x_{N-1}+2\eta_{N-1}'\mathbb{E}x_{N-1}]\\
&=\inf_{u_{N-2}}\mathbb{E}\Big\{2(u_{N-2}-\mathbb{E}u_{N-2})'[(B_{N-2}'P_{N-1}A_{N-2}+D_{N-2}'P_{N-1}C_{N-2}+S'_{N-2})](x_{N-2}-\mathbb{E}x_{N-2})\\
&~~~~~~~~~+(u_{N-2}-\mathbb{E}u_{N-2})'\Upsilon_{N-2}(u_{N-2}-\mathbb{E}u_{N-2})
+2\mathbb{E}u_{N-2}'[(B_{N-2}+\bar{B}_{N-2})'\Pi_{N-1}(A_{N-2}+\bar{A}_{N-2})\\
&~~~~~~~~~+(D_{N-2}+\bar{D}_{N-2})'P_{N-1}(C_{N-2}+\bar{C}_{N-2})+S'_{N-2}+\bar{S}'_{N-2}]\mathbb{E}x_{N-2}
+(\mathbb{E}u_{N-2})'\bar{\Upsilon}_{N-2}\mathbb{E}u_{N-2}\\
&~~~~~~~~~+2(\mathbb{E}u_{N-2})'\eta_{N-2}+
(x_{N-2}-\mathbb{E}x_{N-2})'(Q_{N-2}+A_{N-2}'P_{N-1}A_{N-2}+C_{N-2}'P_{N-1}C_{N-2})\\
&~~~~~~~~~\times(x_{N-2}-\mathbb{E}x_{N-2})
+(\mathbb{E}x_{N-2})'[Q_{N-2}+\bar{Q}_{N-2}+(C_{N-2}+\bar{C}_{N-2})'P_{N-1}(C_{N-2}+\bar{C}_{N-2})\\
&~~~~~~~~~+(A_{N-2}+\bar{A}_{N-2})\Pi_{N-1} (A_{N-2}+\bar{A}_{N-2})]\mathbb{E}x_{N-2}
+2(\mathbb{E}x_{N-2})'[(C_{N-2}+\bar{C}_{N-2})'P_{N-1}\sigma_{N-2}\\
&~~~~~~~~~+(A_{N-2}+\bar{A}_{N-2})'(\Pi_{N-1} b_{N-2}+\eta_{N-1})+q_{N-2}+\bar{q}_{N-2}]+b_{N-1}'\Pi_{N} b_{N-1}\\
&~~~~~~~~~+\sigma_{N-1}'\Pi_{N} \sigma_{N-1}+2\eta_{N}' b_{N-1}+2\eta_{N-1}'b_{N-2}
-\eta_{N-1}'\bar{\Upsilon}_{N-1}\eta_{N-1}\Big\}\\
&\doteq \bar{V}(N-2,\xi),
\end{align*}
where $x_{N-2}=\xi$. Using Lemma 2.2 again, we have
$\bar{\Upsilon}_{N-2}\bar{\Upsilon}^{\dagger}_{N-2}\eta_{N-2}=\eta_{N-2}$,
which implies LRE \eqref{16} is solvable for $k=N-2$. Besides,
\begin{align*}
V(N-2,\xi)\leq\bar{V}(N-2,\xi) &=\sum\limits_{k=N-2}^{N-1}[2\eta_{k+1}'b_{k}+b_{N-1}'\Pi_{k+1}b_{k}+\sigma_{k}'P_{k+1}\sigma_{k}+\zeta_{k}'v_{k}]+2\eta_{N-2}'\mathbb{E}x_{N-2}\\
&~~~~~~~~~+\mathbb{E}[(x_{N-2}-\mathbb{E}x_{N-2})'P_{N-2}(x_{N-2}-\mathbb{E}x_{N-2})]+(\mathbb{E}x_{N-2})'\Pi_{N-2}\mathbb{E}x_{N-2}.
\end{align*}
Next, assume that \eqref{16}-\eqref{19} hold at $k=t+1, t+2, \ldots, N-1, N$, and
\begin{align*}
V(t+1,\xi)\leq\bar{V}(t+1,\xi) &=\sum\limits_{k=t+1}^{N-1}[2\eta_{k+1}'b_{k}+b_{N-1}'\Pi_{k+1}b_{k}+\sigma_{k}'P_{k+1}\sigma_{k}+\zeta_{k}'v_{k}]+2\eta_{t+1}'\mathbb{E}x_{t+1}\\
&~~~~~~~~~+\mathbb{E}[(x_{t+1}-\mathbb{E}x_{t+1})'P_{t+1}(x_{t+1}-\mathbb{E}x_{t+1})]+(\mathbb{E}x_{t+1})'\Pi_{t+1} \mathbb{E}x_{t+1}.
\end{align*}
For $k=t$,
\begin{align*}
V(t,\xi)&=\inf_{u_{t}}\mathbb{E}\sum\limits_{k=t}^{N-1}[x_{k}'Q_{k}x_{k}+(\mathbb{E}x_{k})'\bar{Q}_{k}\mathbb{E}x_{k}+u_{k}'R_{k}u_{k}
+(\mathbb{E}u_{k})'\bar{R}_{k}\mathbb{E}u_{k}+2q_{k}'x_{k}+2\rho_{k}'u_{k}\\
&~~~~~~~~~+2\bar{q}_{k}'x_{k}+2\bar{\rho}_{k}'u_{k}+x_{N}'Gx_{N}+(\mathbb{E}x_{N})'\bar{G}\mathbb{E}x_{N}+2g'x_{N}+2\bar{g}'\mathbb{E}x_{N}]\\
&\leq \inf_{u_{t}}\mathbb{E}[x_{t}'Q_{t}x_{t}+(\mathbb{E}x_{t})'\bar{Q}_{t}\mathbb{E}x_{t}+u_{t}'R_{t}u_{t}
+(\mathbb{E}u_{t})'\bar{R}_{t}\mathbb{E}u_{t}+2q_{t}'x_{t}\\
&~~~~~~~~~+2\rho_{t}'u_{t}+2\bar{q}_{t}'\mathbb{E}x_{t}+2\bar{\rho}_{t}'\mathbb{E}u_{t}+V(t+1,x_{t+1})]\\
&\leq \inf_{u_{t}}\mathbb{E}\Bigg\{x_{t}'Q_{t}x_{t}+(\mathbb{E}x_{t})'\bar{Q}_{t}\mathbb{E}x_{t}+u_{t}'R_{t}u_{t}
+(\mathbb{E}u_{t})'\bar{R}_{t}\mathbb{E}u_{t}+2q_{t}'x_{t}+2\rho_{t}'u_{t}+2\bar{q}_{t}'\mathbb{E}x_{t}+2\bar{\rho}_{t}'\mathbb{E}u_{t}\\
&~~~~~~~~~+\mathbb{E}[(x_{t+1}-\mathbb{E}x_{t+1})'P_{t+1}(x_{t+1}-\mathbb{E}x_{t+1})]+(\mathbb{E}x_{t+1})'\Pi_{t+1}\mathbb{E}x_{t+1}
+2\eta_{t+1}'\mathbb{E}x_{t+1}\\
&~~~~~~~~~+\sum\limits_{k=t}^{N-1}(2\eta_{k+1}'b_{k}+b_{k}'\Pi_{k+1}b_{k}+\sigma_{k}'P_{k+1}\sigma_{k}-\eta_{k}'v_{k})\Bigg\}.
\end{align*}
Similarly, we have
$\bar{\Upsilon}_{t}\bar{\Upsilon}^{\dagger}_{t}\eta_{t}=\eta_{t}$,
thus the LRE \eqref{16} is solvable at $k=t$ and
\begin{align*}
V(t,\xi)&\leq\bar{V}(t,\xi) =\sum\limits_{k=t}^{N-1}[2\eta_{k+1}'b_{k}+b_{N-1}'\Pi_{k+1}b_{k}+\sigma_{k}'P_{k+1}\sigma_{k}+\zeta_{k}'v_{k}]\\
&~~~~~~~~~+\mathbb{E}[(x_{t}-\mathbb{E}x_{t})'P_{t}(x_{t}-\mathbb{E}x_{t})]+(\mathbb{E}x_{t})'\Pi_{t} \mathbb{E}x_{t}+2\eta_{t}'\mathbb{E}x_{t}.
\end{align*}
Using the induction method, the necessity is completed.

\noindent
Sufficiency. Assume that $(\Theta^{\ast},\bar{\Theta}^{\ast},u^{\ast})$ is defined by \eqref{18}, we obtain
\begin{eqnarray}\label{20}
  \left\{
  \begin{array}{ll}
B_{k}'P_{k+1}A_{k}+D_{k}'P_{k+1}C_{k}+S_{k}'=-\Upsilon_{k}\Theta^{\ast}_{k},\\
(B_{k}+\bar{B}_{k})'\Pi_{k+1}(A_{k}+\bar{A}_{k})+(D_{k}+\bar{D}_{k})'P_{k+1}(C_{k}+\bar{C}_{k})+S_{k}'+\bar{S}_{k}'=-\bar{\Upsilon}_{k}(\Theta^{\ast}_{k}+\bar{\Theta}^{\ast}_{k}).
\end{array}\right.
\end{eqnarray}
Using (4.7) and \eqref{16}, we have
\begin{align}\label{21}
&\mathbb{E}\{(x_{N}-\mathbb{E}x_{N})'P_{N}(x_{N}-\mathbb{E}x_{N})-(x_{l}-\mathbb{E}x_{l})'P_{l}(x_{l}-\mathbb{E}x_{l})
+(\mathbb{E}x_{N})'\Pi_{N}\mathbb{E}x_{N}-(\mathbb{E}x_{l})'\Pi_{l}\mathbb{E}x_{l}\}\nonumber\\
&~~~~=\mathbb{E}\sum\limits_{k=l}^{N-1}\{(x_{k+1}-\mathbb{E}x_{k+1})'P_{k+1}(x_{k+1}-\mathbb{E}x_{k+1})+(\mathbb{E}x_{k+1})'\Pi_{k+1}\mathbb{E}x_{k+1}\nonumber\\
&~~~~~~~~~~~~~-(x_{k}-\mathbb{E}x_{k})'P_{k}(x_{k}-\mathbb{E}x_{k})-(\mathbb{E}x_{k})'\Pi_{k}\mathbb{E}x_{k}\}\nonumber\\
&~~~~=\mathbb{E}\sum\limits_{k=l}^{N-1}\{(x_{k}-\mathbb{E}x_{k})'(A_{k}'P_{k+1}A_{k}+C_{k}'P_{k+1}C_{k}-P_{k})(x_{k}-\mathbb{E}x_{k})+2(x_{k}-\mathbb{E}x_{k})'(A_{k}'P_{k+1}B_{k}\nonumber\\
&~~~~~~~~~~~~~+C_{k}'P_{k+1}D_{k})(u_{k}-\mathbb{E}u_{k})+(u_{k}-\mathbb{E}u_{k})'(B_{k}'P_{k+1}B_{k}+D_{k}'P_{k+1}D_{k})(u_{k}-\mathbb{E}u_{k})\nonumber\\
&~~~~~~~~~~~~~+(\mathbb{E}x_{k})'[(C_{k}+\bar{C}_{k})'P_{k+1}(C_{k}+\bar{C}_{k})+(A_{k}+\bar{A}_{k})\Pi_{k+1}(A_{k}+\bar{A}_{k})-\Pi_{k+1}]\mathbb{E}x_{k}\nonumber\\
&~~~~~~~~~~~~~+2(\mathbb{E}x_{k})'[(C_{k}+\bar{C}_{k})'P_{k+1}(D_{k}+\bar{D}_{k})+(A_{k}+\bar{A}_{k})\Pi_{k+1}(B_{k}+\bar{B}_{k})]\mathbb{E}u_{k}\nonumber\\
&~~~~~~~~~~~~~+(\mathbb{E}u_{k})'[(D_{k}+\bar{D}_{k})'P_{k+1}(D_{k}+\bar{D}_{k})+(B_{k}+\bar{B}_{k})\Pi_{k+1}(B_{k}+\bar{B}_{k})-\Pi_{k+1}]\mathbb{E}u_{k}\nonumber\\
&~~~~~~~~~~~~~+2(\mathbb{E}x_{k})'(C_{k}+\bar{C}_{k})P_{k+1}\sigma_{k}
+2(\mathbb{E}u_{k})'(D_{k}+\bar{D}_{k})P_{k+1}\sigma_{k}+\sigma_{k}\Pi_{k+1}\sigma_{k}\nonumber\\
&~~~~~~~~~~~~~+2(\mathbb{E}x_{k})'(A_{k}+\bar{A}_{k})\Pi_{k+1} b_{k}
+2(\mathbb{E}u_{k})'(B_{k}+\bar{B}_{k})\Pi_{k+1} b_{k}
+b_{k}'P_{k+1}b_{k}\},
\end{align}
and
\begin{align}\label{22}
\mathbb{E}(\eta_{N}'x_{N})-\mathbb{E}(\eta_{l}'x_{l})
&=\sum\limits_{k=l}^{N-1}[\mathbb{E}(\eta_{k+1}'x_{k+1})-\mathbb{E}(\eta_{k}'x_{k})]\nonumber\\
&=\sum\limits_{k=l}^{N-1}[\eta_{k+1}'(A_{k}+\bar{A}_{k})\mathbb{E}x_{k}-\eta_{k}'\mathbb{E}x_{k}]
+\sum\limits_{k=l}^{N-1}[\eta_{k+1}'(B_{k}+\bar{B}_{k})\mathbb{E}u_{k}+\eta_{k+1}'b_{k}].
\end{align}
By adding \eqref{21}-\eqref{22} into the cost functional (1.2), it follows that
\begin{align*}
&J(l,\xi;u)+\mathbb{E}\Big\{(x_{N}-\mathbb{E}x_{N})'P_{N}(x_{N}-\mathbb{E}x_{N})-(x_{l}-\mathbb{E}x_{l})'P_{l}(x_{l}-\mathbb{E}x_{l})\\
&~~~~+(\mathbb{E}x_{N})'\Pi_{N}\mathbb{E}x_{N}-(\mathbb{E}x_{l})'\Pi_{l}\mathbb{E}x_{l}\Big\}+2\eta_{N}'\mathbb{E}x_{N}-2\eta_{l}'\mathbb{E}x_{l}\\
&~~~~~~~~=\mathbb{E}\sum\limits_{k=l}^{N-1}\Big\{(x_{k}-\mathbb{E}x_{k})'(Q_{k}+A_{k}'P_{k+1}A_{k}+C_{k}'P_{k+1}C_{k}-P_{k})(x_{k}-\mathbb{E}x_{k})\nonumber\\
&~~~~~~~~~~~~~~~~~~~+2(x_{k}-\mathbb{E}x_{k})'(A_{k}'P_{k+1}B_{k}+C_{k}'P_{k+1}D_{k})(u_{k}-\mathbb{E}u_{k})
+(u_{k}-\mathbb{E}u_{k})'R_{k}(u_{k}-\mathbb{E}u_{k})\nonumber\\
&~~~~~~~~~~~~~~~~~~~+(\mathbb{E}x_{k})'[Q_{k}+\bar{Q}_{k}+(C_{k}+\bar{C}_{k})'P_{k+1}(C_{k}+\bar{C}_{k})+(A_{k}+\bar{A}_{k})'\Pi_{k+1}(A_{k}+\bar{A}_{k})-\Pi_{k+1}]\mathbb{E}x_{k}\nonumber\\
&~~~~~~~~~~~~~~~~~~~+2(\mathbb{E}x_{k})'[(C_{k}+\bar{C}_{k})'P_{k+1}(D_{k}+\bar{D}_{k})+(A_{k}+\bar{A}_{k})'\Pi_{k+1}(B_{k}+\bar{B}_{k})]\mathbb{E}u_{k}\nonumber\\
&~~~~~~~~~~~~~~~~~~~+(\mathbb{E}u_{k})'[R_{k}+\bar{R}_{k}+(D_{k}+\bar{D}_{k})'P_{k+1}(D_{k}+\bar{D}_{k})+(B_{k}+\bar{B}_{k})'\Pi_{k+1}(B_{k}+\bar{B}_{k})]\mathbb{E}u_{k}\nonumber\\
&~~~~~~~~~~~~~~~~~~~+2(\mathbb{E}x_{k})'(C_{k}+\bar{C}_{k})P_{k+1}\sigma_{k}
+2(\mathbb{E}u_{k})'(D_{k}+\bar{D}_{k})P_{k+1}\sigma_{k}+b_{k}'\Pi_{k+1}b_{k}+2\rho_{k}'u_{k}+2\bar{\rho}_{k}'\mathbb{E}u_{k}\nonumber\\
&~~~~~~~~~~~~~~~~~~~+2(\mathbb{E}x_{k})'(A_{k}+\bar{A}_{k})\Pi_{k+1} b_{k}
+2(\mathbb{E}u_{k})'(B_{k}+\bar{B}_{k})\Pi_{k+1} b_{k}
+\sigma_{k}'P_{k+1}\sigma_{k}+2q_{k}'x_{k}+2\bar{q}_{k}'\mathbb{E}x_{k}\\
&~~~~~~~~~~~~~~~~~~~+2\eta_{k+1}'(A_{k}+\bar{A}_{k})\mathbb{E}x_{k}-2\eta_{k}'\mathbb{E}x_{k}
+2\eta_{k+1}'(B_{k}+\bar{B}_{k})\mathbb{E}u_{k}
+2\eta_{k+1}'b_{k}\Big\}\\
&~~~~~~~~~~~~~~~~~~~+(x_{N}-\mathbb{E}x_{N})'G(x_{N}-\mathbb{E}x_{N})
+(\mathbb{E}x_{N})'(G+\bar{G})\mathbb{E}x_{N}+2\mathbb{E}(g'x_{N}+\bar{g}'\mathbb{E}x_{N}).
\end{align*}
Furthermore, using the complete squares method, it is not difficult to see
\begin{align*}
J(l,\xi;u)&=\mathbb{E}\sum\limits_{k=l}^{N-1}\Big\{[u_{k}-\mathbb{E}u_{k}+\Upsilon_{k}^{\dagger}(B_{k}'P_{k+1}A_{k}+D_{k}'P_{k+1}C_{k}+S_{k}')(x_{k}-\mathbb{E}x_{k})]'
\Upsilon_{k}\{u_{k}-\mathbb{E}u_{k}\\
&~~~~~~~~~~~~+\Upsilon_{k}^{\dagger}(B_{k}'P_{k+1}A_{k}+D_{k}'P_{k+1}C_{k}+S_{k}')(x_{k}-\mathbb{E}x_{k})\}+\{\mathbb{E}u_{k}+\bar{\Upsilon}_{k}^{\dagger}[S_{k}'+\bar{S}_{k}'\\
&~~~~~~~~~~~~+(B_{k}+\bar{B}_{k})'\Pi_{k+1}(A_{k}+\bar{A}_{k})
+(D_{k}+\bar{D}_{k})'P_{k+1}(C_{k}+\bar{C}_{k})]\mathbb{E}x_{k}+\bar{\Upsilon}_{k}^{\dagger}\eta_{k}\}'\bar{\Upsilon}_{k}\\
&~~~~~~~~~~~~\times\{\mathbb{E}u_{k}+\bar{\Upsilon}_{k}^{\dagger}
[(B_{k}+\bar{B}_{k})'\Pi_{k+1}(A_{k}+\bar{A}_{k})+(D_{k}+\bar{D}_{k})'P_{k+1}(C_{k}+\bar{C}_{k})+S_{k}'\\
&~~~~~~~~~~~~+\bar{S}_{k}']\mathbb{E}x_{k}+\bar{\Upsilon}_{k}^{\dagger}\eta_{k}\}+2\eta_{k+1}'b_{k}-\zeta_{k}'\bar{\Upsilon}_{k}^{\dagger}\zeta_{k}+b_{k}\Pi_{k+1}b_{k}+\sigma_{k}'P_{k+1}\sigma_{k}\Big\}\\
&~~~~~~~~~~~~+\mathbb{E}[(x_{l}-\mathbb{E}x_{l})'P_{l}(x_{l}-\mathbb{E}x_{l})]
+(\mathbb{E}x_{l})'\Pi_{l}\mathbb{E}x_{l}+2\eta_{l}'\mathbb{E}x_{l},
\end{align*}
since $\Upsilon_{k}, \bar{\Upsilon}_{k}\geq0$, one gets
\begin{align*}
J(l,\xi;u)&\geq\mathbb{E}\sum\limits_{k=l}^{N-1}(2\eta_{k+1}'b_{k}-\zeta_{k}'\bar{\Upsilon}_{k}^{\dagger}\zeta_{k}+b_{k}'\Pi_{k+1}b_{k}+\sigma_{k}'P_{k+1}\sigma_{k})\\
&~~~~~~~~~~~~+\mathbb{E}[(x_{l}-\mathbb{E}x_{l})'P_{l}(x_{l}-\mathbb{E}x_{l})]
+(\mathbb{E}x_{l})'\Pi_{l}\mathbb{E}x_{l}+2\eta_{l}'\mathbb{E}x_{l}\\
&=J(l,\xi;\Theta^{\ast}x^{\ast}+\bar{\Theta}^{\ast}\mathbb{E}x^{\ast}+u^{\ast}),~~\forall~(l,\xi)\in \mathbb{N}_{0}\times \mathbb{R}^{n}.
\end{align*}
Thus, $(\Theta^{\ast},\bar{\Theta}^{\ast},u^{\ast})$ ia an optimal closed-loop strategy of Problem (MF-LQ) and \eqref{19} holds.
\end{proof}
\noindent
\textbf{Corollary 4.1.} Problem (MF-LQ)$^{0}$ is closed-loop solvable if and only if the GRE (4.7) exist a pair regular solution $(P,\Pi)$.

\vskip2mm\noindent
\textbf{Lemma 4.2.} Assume that ($\mathbf{A2}$) holds, then for any $\Theta_{k}, \bar{\Theta}_{k}\in \mathbb{R}^{m\times n}$, the solution $(P_{k},\Pi_{k})$ to the Lyapunov equation
\begin{eqnarray}\label{23}
  \left\{
  \begin{array}{ll}
P_{k}=(A_{k}+B_{k}\Theta_{k})'P_{k+1}(A_{k}+B_{k}\Theta_{k})+(C_{k}+D_{k}\Theta_{k})'P_{k+1}(C_{k}+D_{k}\Theta_{k})\\
~~~~~~~~~~~~~+\Theta_{k}'R_{k}\Theta_{k}+S_{k}'\Theta_{k}+\Theta_{k}'S_{k}+Q_{k},\\
\Pi_{k}=Q_{k}+\bar{Q}_{k}+[A_{k}+\bar{A}_{k}+(B_{k}+\bar{B}_{k})(\Theta_{k}+\bar{\Theta}_{k})]'\Pi_{k+1}[A_{k}+\bar{A}_{k}+(B_{k}+\bar{B}_{k})(\Theta_{k}+\bar{\Theta}_{k})]\\
~~~~~~~~~~~~~+[C_{k}+\bar{C}_{k}+(D_{k}+\bar{D}_{k})(\Theta_{k}+\bar{\Theta}_{k})]'P_{k+1}[C_{k}+\bar{C}_{k}+(D_{k}+\bar{D}_{k})(\Theta_{k}+\bar{\Theta}_{k})],\\
P_{N}=G,~~\Pi_{N}=G+\bar{G},
\end{array}\right.
\end{eqnarray}
satisfies
\begin{eqnarray}\label{24}
  \left\{
  \begin{array}{ll}
R_{k}+B'_{k}P_{k}B_{k}+D'_{k}P_{k}D_{k}\geq  \alpha I,~~~~P_{k}\geq \beta I,\\
R_{k}+\bar{R}_{k}+(B_{k}+\bar{B}_{k})'\Pi_{k+1}(B_{k}+\bar{B}_{k})+(D_{k}+\bar{D}_{k})'P_{k+1}(D_{k}+\bar{D}_{k})\geq \alpha I,~~~~\Pi_{k}\geq \beta I,
\end{array}\right.
\end{eqnarray}
where $\beta\in \mathbb{R}$ is the constant in \eqref{46}.
\begin{proof} For any $u\in \mathscr{U}_{ad}$, let $x^{(0)}$ be the solution of
\begin{eqnarray*}
  \left\{
  \begin{array}{ll}
x_{k+1}^{(0)}=(A_{k}+B_{k}\Theta_{k})x_{k}^{(0)}+B_{k}u_{k}+[\bar{A}_{k}+B_{k}\bar{\Theta}_{k}+\bar{B}_{k}(\Theta_{k}+\bar{\Theta}_{k})]\mathbb{E}x_{k}^{(0)}
+\bar{B}_{k}\mathbb{E}u_{k}\\
~~~~~~~~~~~~+\Big\{(C_{k}+D_{k}\Theta_{k})x_{k}^{(0)}+D_{k}u_{k}+[\bar{C}_{k}+D_{k}\bar{\Theta}_{k}+\bar{D}_{k}(\Theta_{k}+\bar{\Theta}_{k})]\mathbb{E}x_{k}^{(0)}
+\bar{D}_{k}\mathbb{E}u_{k}\Big\}\omega_{k},\\
x_{l}^{(0)}=0.
\end{array}\right.
\end{eqnarray*}
Using ($\mathbf{A2}$), it follows that
\begin{align}\label{ssss}
&\beta\mathbb{E}\sum\limits_{k=l}^{N-1}\left|\Theta_{k}(x_{k}^{(0)}-\mathbb{E}x_{k}^{(0)})+\bar{\Theta}_{k}\mathbb{E}x_{k}^{(0)}+u_{k}\right|^{2}\leq J^{0}\Big(l,0;\Theta (x_{k}^{0}-\mathbb{E}x_{k}^{(0)})+\bar{\Theta}_{k}\mathbb{E}x_{k}^{(0)}+u_{k}\Big)\nonumber\\
&~~~~=\mathbb{E}\sum\limits_{k=l}^{N-1}\bigg\{\Big\langle \Upsilon_{k}[\Theta_{k}(x^{(0)}_{k}-\mathbb{E}x^{(0)}_{k})+u_{k}-\mathbb{E}u_{k}],\Theta_{k}(x^{(0)}_{k}-\mathbb{E}x^{(0)}_{k})+u_{k}-\mathbb{E}u_{k}\Big\rangle+\Big\langle \bar{\Theta}_{k}\bar{\Upsilon}_{k}\mathbb{E}x^{(0)}_{k},\mathbb{E}x^{(0)}_{k}\Big\rangle\nonumber\\
&~~~~~~~~~~~~~~~~+2\Big\langle [B_{k}'P_{k+1}A_{k}+D_{k}'P_{k+1}C_{k}+S_{k}'][\Theta_{k}(x^{(0)}_{k}-\mathbb{E}x^{(0)}_{k})+u_{k}-\mathbb{E}u_{k}],x^{(0)}_{k}-\mathbb{E}x^{(0)}_{k}\Big\rangle\nonumber\\
&~~~~~~~~~~~~~~~~+\Big\langle \Theta_{k}\Upsilon_{k}(x^{(0)}_{k}-\mathbb{E}x^{(0)}_{k}),x^{(0)}_{k}-\mathbb{E}x^{(0)}_{k}\Big\rangle+
\Big\langle \bar{\Upsilon}_{k}(\Theta_{k}\mathbb{E}x^{(0)}_{k}+\mathbb{E}u_{k}),\Theta_{k}\mathbb{E}x^{(0)}_{k}+\mathbb{E}u_{k}\Big\rangle\nonumber\\
&~~~~~~~~~~~~~~~~+2\Big\langle [(B_{k}+\bar{B}_{k})'\Pi_{k+1}(A_{k}+\bar{A}_{k})+(D_{k}+\bar{D}_{k})'P_{k+1}(C_{k}+\bar{C}_{k})+S_{k}'+\bar{S}_{k}']\nonumber\\
&~~~~~~~~~~~~~~~~\times(\Theta_{k}\mathbb{E}x^{(0)}_{k}+\mathbb{E}u_{k}),\mathbb{E}x^{(0)}_{k}\Big\rangle\bigg\}\nonumber\\
&~~~~=\mathbb{E}\sum\limits_{k=l}^{N-1}\Big\{2\Big\langle [B_{k}'P_{k+1}A_{k}+D_{k}'P_{k+1}C_{k}+S_{k}'+\Upsilon_{k}\Theta_{k}](x^{(0)}_{k}-\mathbb{E}x^{(0)}_{k}),u_{k}-\mathbb{E}u_{k}\Big\rangle+\Big\langle \Upsilon_{k}u_{k},u_{k}\Big\rangle\Big\}\nonumber\\
&~~~~~~~~~~~~~~~~
+\mathbb{E}\sum\limits_{k=l}^{N-1}\Big\{2\Big\langle (B_{k}+\bar{B}_{k})'\Pi_{k+1}(A_{k}+\bar{A}_{k})+(D_{k}+\bar{D}_{k})'P_{k+1}(C_{k}+\bar{C}_{k})+S_{k}'+\bar{S}_{k}'\nonumber\\
&~~~~~~~~~~~~~~~~+\bar{\Upsilon}_{k}\bar{\Theta}_{k}\mathbb{E}x_{k}^{(0)},\mathbb{E}u_{k}\Big\rangle+\Big\langle\bar{\Upsilon}_{k}\mathbb{E}u_{k},\mathbb{E}u_{k}\Big\rangle\Big\}.
\end{align}
Thus, we have
\begin{align*}
&2\mathbb{E}\sum\limits_{k=l}^{N-1}\bigg\{\Big\langle [B_{k}'P_{k+1}A_{k}+D_{k}'P_{k+1}C_{k}+S_{k}'+(\Upsilon_{k}-\alpha I)\Theta_{k}](x^{(0)}_{k}-\mathbb{E}x^{(0)}_{k}),u_{k}-\mathbb{E}u_{k}\Big\rangle\\
&~~~~~~~~~~+\Big\langle (\Upsilon_{k}-\alpha I)u_{k},u_{k}\Big\rangle+\Big\langle[(B_{k}+\bar{B}_{k})'\Pi_{k+1}(A_{k}+\bar{A}_{k})+(D_{k}+\bar{D}_{k})'P_{k+1}(C_{k}+\bar{C}_{k})\\
&~~~~~~~~~~+S_{k}'+\bar{S}_{k}'+(\bar{\Upsilon}_{k}-\alpha I)\bar{\Theta}_{k}]\mathbb{E}x_{k}^{(0)},\mathbb{E}u_{k}\Big\rangle
+\Big\langle(\bar{\Upsilon}_{k}-\alpha I)\mathbb{E}u_{k},\mathbb{E}u_{k}\Big\rangle\bigg\}\\
&~~~~~~=\beta\mathbb{E}\sum\limits_{k=l}^{N-1}\left|\Theta_{k} (x_{k}^{(0)}-\mathbb{E}x_{k}^{(0)})+\bar{\Theta}_{k}\mathbb{E}x_{k}^{(0)}\right|^{2}\geq0.
\end{align*}
Since the above equality holds for any $u\in \mathscr{U}_{ad}$, we could deduce
$\Upsilon_{k}\geq \alpha I,~\bar{\Upsilon}_{k}\geq \alpha I.$
Next we prove $P_{k}\geq \beta I$, $\Pi_{k}\geq \beta I$.
Let $x_{k}$ be the solution of
\begin{eqnarray*}
  \left\{
  \begin{array}{ll}
x_{k+1}=(A_{k}+B_{k}\Theta_{k})x_{k}+B_{k}u_{k}+[\bar{A}_{k}+B_{k}\bar{\Theta}_{k}+\bar{B}_{k}(\Theta_{k}+\bar{\Theta}_{k})]\mathbb{E}x_{k}
+\bar{B}_{k}\mathbb{E}u_{k}\\
~~~~~~~~~~~~+\Big\{(C_{k}+D_{k}\Theta_{k})x_{k}+D_{k}u_{k}+[\bar{C}_{k}+D_{k}\bar{\Theta}_{k}+\bar{D}_{k}(\Theta_{k}+\bar{\Theta}_{k})]\mathbb{E}x_{k}
+\bar{D}_{k}\mathbb{E}u_{k}\Big\}\omega_{k},\\
x_{l}=\xi,
\end{array}\right.
\end{eqnarray*}
according to Lemma 4.1, it follows that
\begin{align*}
\beta|\xi|^{2}&\leq V^{0}(l,\xi)\leq J^{0}(l,\xi;\Theta_{k} (x_{k}^{0}-\mathbb{E}x_{k}^{0})+\bar{\Theta}_{k}\mathbb{E}x_{k}^{0}+u_{k})\\
&=\mathbb{E}\langle P_{l}(x_{l}-\mathbb{E}x_{l}),x_{l}-\mathbb{E}x_{l}\rangle+\mathbb{E}\langle \Pi_{l}\mathbb{E}x_{l},\mathbb{E}x_{l}\rangle
+\mathbb{E}\sum\limits_{k=l}^{N-1}\Big\{2\langle [B_{k}'P_{k+1}A_{k}+D_{k}'P_{k+1}C_{k}\\
&~~~~~~+S_{k}'+\Upsilon_{k}\Theta_{k}](x^{0}_{k}-\mathbb{E}x^{0}_{k}),u_{k}-\mathbb{E}u_{k}\rangle+\langle \Upsilon_{k}u_{k},u_{k}\rangle
+2\langle (B_{k}+\bar{B}_{k})'\Pi_{k+1}(A_{k}+\bar{A}_{k})\\
&~~~~~~+(D_{k}+\bar{D}_{k})'P_{k+1}(C_{k}+\bar{C}_{k})+S_{k}'+\bar{S}_{k}'+\bar{\Upsilon}_{k}\bar{\Theta}_{k}\mathbb{E}x_{k}^{0},\mathbb{E}u_{k}\rangle
+\langle\bar{\Upsilon}_{k}\mathbb{E}u_{k},\mathbb{E}u_{k}\rangle\Big\}.
\end{align*}
Specially, set $u_{k}=0$, one has
\begin{align*}
\mathbb{E}\langle P_{l}(x_{l}-\mathbb{E}x_{l}),x_{l}-\mathbb{E}x_{l}\rangle+(\mathbb{E}x_{l})'\Pi_{l}\mathbb{E}x_{l}\geq\beta|\xi|^{2},~\forall~(l,\xi)\in\mathbb{N}_{0}\times \mathbb{R}^{n},
\end{align*}
further, let $\mathbb{E}x_{l}=0$, we get $\mathbb{E}\langle P_{l}x_{l},x_{l}\rangle\geq\beta|\xi|^{2}$. Besides, let $x_{l}=\mathbb{E}x_{l}\neq0$,
$(\mathbb{E}x_{l})'\Pi_{l}\mathbb{E}x_{l}\geq\beta|\xi|^{2}$.
\end{proof}

\vskip0mm\noindent
\textbf{Theorem 4.5.} The following statements are equivalent:

\noindent
(i) The map $u\mapsto J^{0}(0,0;u)$ is uniformly convex, in other words, there is a $\alpha>0$ such that
\indent $J^{0}(0,0;u)\geq\alpha\mathbb{E}\sum_{k=0}^{N-1}|u_{k}|^{2}$.

\noindent
(ii) The GRE (4.7) exists a pair strongly regular solution $(P,\Pi)$.
\begin{proof} (i) $\Rightarrow$ (ii) Let $(P^{(0)},\Pi^{(0)})$ be the solution of
\begin{eqnarray}\label{25}
  \left\{
  \begin{array}{ll}
P_{k}^{(0)}=Q_{k}+A_{k}'P_{k+1}^{(0)}A_{k}+C_{k}'P_{k+1}^{(0)}C_{k},~~~P_{N}^{(0)}=G,\\
\Pi_{k}^{(0)}=Q_{k}+\bar{Q}_{k}+(A_{k}+\bar{A}_{k})'\Pi_{k+1}^{(0)}(A_{k}+\bar{A}_{k})+(C_{k}+\bar{C}_{k})'P_{k+1}^{(0)}(C_{k}+\bar{C}_{k}),~~~\Pi_{N}^{(0)}=G+\bar{G},
\end{array}\right.
\end{eqnarray}
according to Lemma 4.2 ($\Theta_{k}, \bar{\Theta}_{k}=0$), we have
\begin{eqnarray*}
  \left\{
  \begin{array}{ll}
R_{k}+B'_{k}P_{k+1}^{(0)}B_{k}+D'_{k}P_{k+1}^{(0)}D_{k}\geq \alpha I,~~~~P_{k}^{(0)}\geq \beta I,\\
R_{k}+\bar{R}_{k}+(B_{k}+\bar{B}_{k})'\Pi_{k+1}^{(0)}(B_{k}+\bar{B}_{k})+(D_{k}+\bar{D}_{k})'P_{k+1}^{(0)}(D_{k}+\bar{D}_{k})\geq \alpha I,~~~~\Pi_{k}^{(0)}\geq \beta I.
\end{array}\right.
\end{eqnarray*}
For $i=0,1,2,\ldots$, by the induction method, let
\begin{eqnarray}\label{b}
  \left\{
  \begin{array}{ll}
\Theta^{(i)}=-(R_{k}+B'_{k}P_{k+1}^{(i)}B_{k}+D'_{k}P_{k+1}^{(i)}D_{k})^{-1}(B_{k}'P_{k+1}^{(i)}A_{k}+D_{k}'P_{k+1}^{(i)}C_{k}+S'_{k}),\\
\bar{\Theta}^{(i)}=[R_{k}+\bar{R}_{k}+(B_{k}+\bar{B}_{k})'\Pi_{k+1}^{(i)}(B_{k}+\bar{B}_{k})+(D_{k}+\bar{D}_{k})'P_{k+1}^{(i)}(D_{k}+\bar{D}_{k})]^{-1}\\
~~~~~~~~~~~~~~~\times[(B_{k}+\bar{B}_{k})'\Pi_{k+1}^{(i)}(A_{k}+\bar{A}_{k})+(D_{k}+\bar{D}_{k})'P_{k+1}^{(i)}(C_{k}+\bar{C}_{k})+S'_{k}+\bar{S}'_{k}],\\
A^{(i)}_{k}=A_{k}+B_{k}\Theta^{(i)}_{k},~\bar{A}^{(i)}_{k}=A_{k}+\bar{A}_{k}+(B_{k}+\bar{B}_{k})(\Theta_{k}^{(i)}+\bar{\Theta}_{k}^{(i)}),\\
C^{(i)}_{k}=C_{k}+D_{k}\Theta^{(i)}_{k},~\bar{C}^{(i)}_{k}=A_{k}+\bar{A}_{k}+(B_{k}+\bar{B}_{k})(\Theta_{k}^{(i)}+\bar{\Theta}_{k}^{(i)}),
\end{array}\right.
\end{eqnarray}
and
\begin{eqnarray*}
  \left\{
  \begin{array}{ll}
P^{(i+1)}_{k}=(A^{(i)}_{k})'P^{(i+1)}_{k+1}A^{(i)}_{k}+(C^{(i)}_{k})'P^{(i+1)}_{k+1}C^{(i)}_{k}
+(\Theta^{(i)}_{k})'R_{k}\Theta^{(i)}_{k}+S_{k}'\Theta^{(i)}_{k}+(\Theta^{(i)}_{k})'S_{k}+Q_{k},\\
\Pi_{k}^{(i+1)}=(A_{k}^{(i)}+\bar{A}_{k}^{(i)})'
\Pi_{k+1}^{(i+1)}(A_{k}^{(i)}+\bar{A}_{k}^{(i)})+(C_{k}^{(i)}+\bar{C}_{k}^{(i)})'P_{k+1}^{(i+1)}
(C_{k}^{(i)}+\bar{C}_{k}^{(i)})+(S_{k}+\bar{S}_{k})'(\Theta^{(i)}_{k}+\bar{\Theta}^{(i)}_{k})\\
~~~~~~~~~~~~~~~~~~+(\Theta^{(i)}_{k}+\bar{\Theta}^{(i)}_{k})'(R_{k}+\bar{R}_{k})(\Theta^{(i)}_{k}+\bar{\Theta}^{(i)}_{k})
+(\Theta^{(i)}_{k}+\bar{\Theta}^{(i)}_{k})'(S_{k}+\bar{S}_{k})+Q_{k}+\bar{Q}_{k},\\
P_{N}^{(i+1)}=G,~~\Pi_{N}^{(i+1)}=G+\bar{G}.
\end{array}\right.
\end{eqnarray*}
By Lemma 4.2, we know
\begin{eqnarray}\label{26}
  \left\{
  \begin{array}{ll}
R_{k}+B'_{k}P_{k+1}^{(i+1)}B_{k}+D'_{k}P_{k+1}^{(i+1)}D_{k}\geq \alpha I,~~~~P_{k}^{(i+1)}\geq \beta I,\\
R_{k}+\bar{R}_{k}+(B_{k}+\bar{B}_{k})'\Pi_{k+1}^{(i+1)}(B_{k}+\bar{B}_{k})+(D_{k}+\bar{D}_{k})'P_{k+1}^{(i+1)}(D_{k}+\bar{D}_{k})\geq \alpha I,~~~~\Pi_{k}^{(i+1)}\geq \beta I.
\end{array}\right.
\end{eqnarray}
Next, we shall show that $\{P^{(i)}\}_{i=1}^{\infty}$ and $\{\Pi^{(i)}\}_{i=1}^{\infty}$ are uniformly convergent. To this end, for $i\geq1$, set
\begin{eqnarray*}
  \left\{
  \begin{array}{ll}
\Delta^{(i)}=P^{(i)}-P^{(i+1)},~~~\nabla^{(i)}=\Theta^{(i-1)}-\Theta^{(i)},\\
\bar{\Delta}^{(i)}=\Pi^{(i)}-\Pi^{(i+1)},~~~\bar{\nabla}^{(i)}=\bar{\Theta}^{(i-1)}-\bar{\Theta}^{(i)},
\end{array}\right.
\end{eqnarray*}
then we obtain
\begin{align}\label{27}
\Delta^{(i)}&=(A^{(i-1)})'P^{(i)}A^{(i-1)}+(C^{(i-1)})'P^{(i)}C^{(i-1)}
+(\Theta^{(i-1)})'R\Theta^{(i-1)}+S'\Theta^{(i-1)}+(\Theta^{(i-1)})'S\nonumber\\
&~~~~~~~~-(A^{(i)})'P^{(i+1)}A^{(i)}-(C^{(i)})'P^{(i+1)}C^{(i)}
-(\Theta^{(i)})'R_{k}\Theta^{(i)}-S'\Theta^{(i)}-(\Theta^{(i)})'S\nonumber\\
&=(A^{(i)})'\Delta^{(i)}A^{(i)}+(A^{(i-1)})'P^{(i)}A^{(i-1)}-(A^{(i)})'P^{(i)}A^{(i)}
+(\Theta^{(i-1)})'R_{k}\Theta^{(i-1)}-(\Theta^{(i)})'R\Theta^{(i)}\nonumber\\
&~~~~~~~~+(C^{(i)})'\Delta^{(i)}C^{(i)}+(C^{(i-1)})'P^{(i)}C^{(i-1)}-(C^{(i)})'P^{(i)}C^{(i)}
+S'\nabla^{(i)}+(\nabla^{(i)})'S.
\end{align}
By virtue of \eqref{b}, it yields that
\begin{eqnarray}\label{28}
  \left\{
  \begin{array}{ll}
A^{(i-1)}-A^{(i)}=B\nabla^{(i)},~~~C^{(i-1)}-C^{(i)}=D\nabla^{(i)},\\
(C^{(i-1)})'P^{(i)}C^{(i-1)}-(C^{(i)})'P^{(i)}C^{(i)}=(\nabla^{(i)})'D'P^{(i)}D\nabla^{(i)}
+(C^{(i)})'P^{(i)}D\nabla^{(i)}+(\nabla^{(i)})'D'P^{(i)}C^{(i)},\\
(\Theta^{(i-1)})'R\Theta^{(i-1)}-(\Theta^{(i)})'R\Theta^{(i)}=(\nabla^{(i)})'R\nabla^{(i)}
+(\nabla^{(i)})'R\Theta^{(i)}+(\Theta^{(i)})'R\nabla^{(i)}.
\end{array}\right.
\end{eqnarray}
Due to
\begin{align*}
B'P^{(i)}A+D'P^{(i)}C+R\Theta^{(i)}+S=B'P^{(i)}A+D'P^{(i)}C+S+(R+B'P^{(i)}B+D'P^{(i)}D)\Theta^{(i)}=0,
\end{align*}
then using \eqref{27}-\eqref{28}, we get
\begin{align}\label{29}
&\Delta^{(i)}-(A^{(i)})'\Delta^{(i)}A^{(i)}-(C^{(i)})'\Delta^{(i)}C^{(i)}\nonumber\\
&~~~~~~=(\nabla^{(i)})'B'P^{(i)}B\nabla^{(i)}
+(A^{(i)})'P^{(i)}B\nabla^{(i)}+(\nabla^{(i)})'B'P^{(i)}A^{(i)}+(\nabla^{(i)})'D'P^{(i)}D\nabla^{(i)}
+(C^{(i)})'P^{(i)}D\nabla^{(i)}\nonumber\\
&~~~~~~~~~~~~~~+(\nabla^{(i)})'D'P^{(i)}C^{(i)}+(\nabla^{(i)})'R\nabla^{(i)}
+(\nabla^{(i)})'R\Theta^{(i)}+(\Theta^{(i)})'R\nabla^{(i)}+S'\nabla^{(i)}+(\nabla^{(i)})'S\nonumber\\
&~~~~~~=(\nabla^{(i)})'(R+B'P^{(i)}B+D'P^{(i)}D)\nabla^{(i)}+[A'P^{(i)}B+C'P^{(i)}D+(\Theta^{(i)})'R+S']\Delta^{(i)}\nonumber\\
&~~~~~~~~~~~~~~+(\Delta^{(i)})'(B'P^{(i)}A+D'P^{(i)}C+R\Theta^{(i)}+S)\nonumber\\
&~~~~~~=(\nabla^{(i)})'(R+B'P^{(i)}B+D'P^{(i)}D)\nabla^{(i)}\geq0,
\end{align}
which, combining with $\Delta^{(i)}_{N}=0$, gives $\Delta^{(i)}_{k}\geq0$. Besides, by \eqref{26} we obtain
\begin{align*}
P_{k}^{(1)}\geq P_{k}^{(i)}\geq P_{k}^{(i+1)}\geq \beta I,~~\forall~k\in\{0,1,\ldots,N\},~~\forall~i\geq1.
\end{align*}
Consequently, the sequence $\{P^{(i)}\}_{i=1}^{\infty}$ is uniformly bounded.
Likewise, we could deduce that
\begin{align*}
&\bar{\Delta}^{(i)}_{k}-(A^{(i)}_{k}+\bar{A}^{(i)}_{k})'\bar{\nabla}^{(i)}_{k+1}(A^{(i)}_{k}+\bar{A}^{(i)}_{k})
-(C^{(i)}_{k}+\bar{C}^{(i)}_{k})'\bar{\Delta}^{(i)}_{k+1}(C^{(i)}_{k}+\bar{C}^{(i)}_{k})\nonumber\\
&~~~~~~=(\bar{\nabla}^{(i)}_{k})'[R_{k}+\bar{R}_{k}+(B_{k}+\bar{B}_{k})'\Pi_{k+1}^{(i)}(B_{k}+\bar{B}_{k})
+(D_{k}+\bar{D}_{k})'P_{k+1}^{(i)}(D_{k}+\bar{D}_{k})]\bar{\nabla}^{(i)}_{k}\geq0,
\end{align*}
then the sequence $\{\Pi^{(i)}\}_{i=1}^{\infty}$ is also uniformly bounded. Namely, there is a constant $d>0$ such that
\begin{eqnarray}\label{30}
  \left\{
  \begin{array}{ll}
|P^{(i)}_{k}|,~|R^{(i)}_{k}|~,|\Pi^{(i)}_{k}|,~|R^{(i)}_{k}+\bar{R}^{(i)}_{k}|\leq d,\\
|\Theta^{(i)}_{k}|\leq d(|B_{k}|+|C_{k}|+|S_{k}|),~~|\bar{\Theta}^{(i)}_{k}|\leq d(|\bar{B}_{k}|+|\bar{C}_{k}|+|\bar{S}_{k}|),\\
|A^{(i)}_{k}|\leq |A_{k}|+d |B_{k}|(|B_{k}|+|C_{k}|+|S_{k}|),~~|C^{(i)}_{k}|\leq |C_{k}|+d (|B_{k}|+|C_{k}|+|S_{k}|),\\
|A^{(i)}_{k}+\bar{A}^{(i)}_{k}|\leq |A_{k}+\bar{A}_{k}|+d |B_{k}+\bar{B}_{k}|(|B_{k}+\bar{B}_{k}|+|C_{k}+\bar{C}_{k}|+|S_{k}+\bar{S}_{k}|),\\
|C^{(i)}_{k}+\bar{C}^{(i)}_{k}|\leq |C_{k}+\bar{C}_{k}|+d (|B_{k}+\bar{B}_{k}|+|C_{k}+\bar{C}_{k}|+|S_{k}+\bar{S}_{k}|),
\end{array}\right.
\end{eqnarray}
where
\begin{align*}
&R^{(i)}_{k}=R_{k}+B_{k}'P^{(i)}_{k+1}B_{k}+D_{k}'P^{(i)}_{k+1}D_{k},\\
&R^{(i)}_{k}+\bar{R}^{(i)}_{k}=R_{k}+\bar{R}_{k}+(B_{k}+\bar{B}_{k})'\Pi_{k+1}^{(i)}(B_{k}+\bar{B}_{k})
+(D_{k}+\bar{D}_{k})'P_{k+1}^{(i)}(D_{k}+\bar{D}_{k}).
\end{align*}
Since
\begin{align*}
\nabla^{(i)}_{k}&=(R^{(i)}_{k})^{-1}(B_{k}'\Delta^{(i-1)}_{k+1}B_{k}+D_{k}'\Delta^{(i-1)}_{k+1}D_{k})(R^{(i-1)}_{k})^{-1}
(B_{k}'P_{k+1}^{(i)}A_{k}+D_{k}'P_{k+1}^{(i)}C_{k}+S_{k}')\\
&~~~~~~~~-(R^{(i-1)}_{k})^{-1}(B'\Delta^{(i-1)}_{k}A_{k}+D'\Delta^{(i-1)}_{k}C_{k}),
\end{align*}
then by \eqref{30}, we get
\begin{align*}
&|(\nabla^{(i)}_{k})'R^{(i)}_{k}\nabla^{(i)}_{k}|\leq (|\Theta^{(i)}_{k}|+|\Theta^{(i-1)}_{k}|)|R^{(i)}_{k}|(|\Theta^{(i-1)}_{k}-\Theta^{(i)}_{k}|)\leq d (|B_{k}|+|C_{k}|+|S_{k}|)^{2}|\Delta_{k}^{(i-1)}|,\\
&|(\bar{\nabla}^{(i)}_{k})'(R^{(i)}_{k}+\bar{R}^{(i)}_{k})\bar{\nabla}^{(i)}_{k}|\leq(|\bar{\Theta}^{(i)}_{k}|+|\bar{\Theta}^{(i-1)}_{k}|)|(R^{(i)}_{k}+\bar{R}^{(i)}_{k})|
(|\bar{\Theta}^{(i-1)}_{k}-\bar{\Theta}^{(i)}_{k}|)\leq d (|\bar{B}_{k}|+|\bar{C}_{k}|+|\bar{S}_{k}|)^{2}|\bar{\Delta}_{k}^{(i-1)}|.
\end{align*}
Combining \eqref{29}-\eqref{30} with $\Delta^{(i)}_{N}=0$, we see
\begin{align*}
|\Delta^{(i)}_{k}|&=(\nabla^{(i)})'R^{i}\nabla^{(i)}+(A^{(i)})'\bar{\nabla}^{(i)}A^{(i)}
+(C^{(i)})'\bar{\Delta}^{(i)}C^{(i)}\leq \phi_{k}(|\nabla^{(i)}_{k}|+|\nabla^{(i-1)}_{k}|),\\
|\bar{\Delta}^{(i)}_{k}|&=(\bar{\nabla}^{(i)})'(R^{(i)}+\bar{R}^{(i)})\bar{\nabla}^{(i)}+(A^{(i)}+\bar{A}^{(i)})'\bar{\nabla}^{(i)}(A^{(i)}+\bar{A}^{(i)})
+(C^{(i)}+\bar{C}^{(i)})'\bar{\Delta}^{(i)}(C^{(i)}+\bar{C}^{(i)})\\
&~~~~\leq \bar{\phi}_{k}(|\bar{\nabla}^{(i)}_{k}|+|\bar{\nabla}^{(i-1)}_{k}|),
\end{align*}
where $\phi_{k}, \bar{\phi}_{k}$ are nonnegative function independent of $\nabla^{(i)}_{k}, \bar{\nabla}^{(i)}_{k}$. These indicate $\{P^{(i)}\}_{i=1}^{\infty}$ and $\{\Pi^{(i)}\}_{i=1}^{\infty}$ are uniformly convergent, and denote it as $P$ and $\Pi$, respectively. Furthermore,
\begin{align*}
&R_{k}+B_{k}'P_{k+1}B_{k}+D_{k}'P_{k+1}D_{k}=\lim_{i\rightarrow\infty}(R_{k}+B_{k}'P_{k+1}^{(i)}B_{k}+D_{k}'P_{k+1}^{(i)}D_{k})\geq \alpha I,\\
&R_{k}+\bar{R}_{k}+(B_{k}+\bar{B}_{k})'\Pi_{k+1}(B_{k}+\bar{B}_{k})
+(D_{k}+\bar{D}_{k})'P_{k+1}(D_{k}+\bar{D}_{k})\\
&~~~~~~~~=\lim_{i\rightarrow\infty}[R_{k}+\bar{R}_{k}+(B_{k}+\bar{B}_{k})'\Pi_{k+1}^{(i)}(B_{k}+\bar{B}_{k})
+(D_{k}+\bar{D}_{k})'P_{k+1}^{(i)}(D_{k}+\bar{D}_{k})]\geq \alpha I,
\end{align*}
and
\begin{eqnarray*}
  \left\{
  \begin{array}{ll}
\Theta^{(i)}\rightarrow -(R_{k}+B'_{k}P_{k+1}B_{k}+D'_{k}P_{k+1}D_{k})^{-1}(B_{k}'P_{k+1}A_{k}+D_{k}'P_{k+1}C_{k}+S'_{k}),\\
\bar{\Theta}^{(i)}\rightarrow [R_{k}+\bar{R}_{k}+(B_{k}+\bar{B}_{k})'\Pi_{k+1}(B_{k}+\bar{B}_{k})+(D_{k}+\bar{D}_{k})'P_{k+1}(D_{k}+\bar{D}_{k})]^{-1}\\
~~~~~~~~~~~~~~~\times[(B_{k}+\bar{B}_{k})'\Pi_{k+1}(A_{k}+\bar{A}_{k})+(D_{k}+\bar{D}_{k})'P_{k+1}(C_{k}+\bar{C}_{k})+S'_{k}+\bar{S}'_{k}],\\
A^{(i)}_{k}\rightarrow A_{k}+B_{k}\Theta_{k},~\bar{A}^{(i)}_{k}\rightarrow A_{k}+\bar{A}_{k}+(B_{k}+\bar{B}_{k})(\Theta_{k}+\bar{\Theta}_{k}),\\
C^{(i)}_{k}\rightarrow C_{k}+D_{k}\Theta_{k},~\bar{C}^{(i)}_{k}\rightarrow A_{k}+\bar{A}_{k}+(B_{k}+\bar{B}_{k})(\Theta_{k}+\bar{\Theta}_{k}).
\end{array}\right.
\end{eqnarray*}
Thus, $(P,\Pi)$ satisfies the following equation
\begin{eqnarray*}
  \left\{
  \begin{array}{ll}
P_{k}=A_{k}'P_{k+1}A_{k}+C_{k}'P_{k+1}C_{k}+\Theta_{k}'R_{k}\Theta_{k}+S_{k}'\Theta_{k}+\Theta_{k}'S_{k}+Q_{k},\\
\Pi_{k}=(A_{k}+\bar{A}_{k})'\Pi_{k+1}(A_{k}+\bar{A}_{k})+(C_{k}+\bar{C}_{k})'P_{k+1}
(C_{k}+\bar{C}_{k})+(\Theta_{k}+\bar{\Theta}_{k})'(S_{k}+\bar{S}_{k})\\
~~~~~~~~~~~~~+(\Theta_{k}+\bar{\Theta}_{k})'(R_{k}+\bar{R}_{k})(\Theta_{k}+\bar{\Theta}_{k})
+(S_{k}+\bar{S}_{k})'(\Theta_{k}+\bar{\Theta}_{k})+Q_{k}+\bar{Q}_{k},\\
P_{N}=G,~~\Pi_{N}=G+\bar{G}.
\end{array}\right.
\end{eqnarray*}
(ii) $\Rightarrow$ (i) Suppose that $(P,\Pi)$ is the strongly regular solution to GRE (4.7), there is a $\alpha>0$ such that
\begin{eqnarray}\label{31}
  \left\{
  \begin{array}{ll}
R_{k}+B'_{k}P_{k+1}B_{k}+D'_{k}P_{k+1}D_{k}\geq \alpha I,~~~~P_{k}\geq \beta I,\\
R_{k}+\bar{R}_{k}+(B_{k}+\bar{B}_{k})'\Pi_{k+1}(B_{k}+\bar{B}_{k})+(D_{k}+\bar{D}_{k})'P_{k+1}(D_{k}+\bar{D}_{k})\geq \alpha I,~~~~\Pi_{k}\geq \beta I.
\end{array}\right.
\end{eqnarray}
Let
\begin{eqnarray*}
  \left\{
  \begin{array}{ll}
\Theta_{k}= -(R_{k}+B'_{k}P_{k}B_{k}+D'_{k}P_{k}D_{k})^{-1}(B'_{k}P_{k}A_{k}+D_{k}'P_{k}C_{k}+S_{k}),\\
\bar{\Theta}_{k}=-[R_{k}+\bar{R}_{k}+(B_{k}+\bar{B}_{k})'\Pi_{k+1}(B_{k}+\bar{B}_{k})+(D_{k}+\bar{D}_{k})'P_{k+1}(D_{k}+\bar{D}_{k})]^{-1}\\
~~~~~~~~~~~~~\times[(B_{k}+\bar{B}_{k})'\Pi_{k+1}(A_{k}+\bar{A}_{k})+(D_{k}+\bar{D}_{k})'P_{k+1}(C_{k}+\bar{C}_{k})+S'_{k}+\bar{S}'_{k}],
\end{array}\right.
\end{eqnarray*}
and for any $u\in \mathscr{U}_{ad}$, $x_{k}^{(u)}$ be the solution of
\begin{eqnarray*}
  \left\{
  \begin{array}{ll}
x_{k+1}^{(u)}=(A_{k}+B_{k}\Theta_{k})x_{k}^{(u)}+B_{k}u_{k}+[\bar{A}_{k}+B_{k}\bar{\Theta}_{k}+\bar{B}_{k}(\Theta_{k}+\bar{\Theta}_{k})]\mathbb{E}x_{k}^{(u)}
+\bar{B}_{k}\mathbb{E}u_{k}\\
~~~~~~~~~~~~~~~~+\Big\{(C_{k}+D_{k}\Theta_{k})x_{k}^{(u)}+D_{k}u_{k}+[\bar{C}_{k}+D_{k}\bar{\Theta}_{k}+\bar{D}_{k}(\Theta_{k}+\bar{\Theta}_{k})]\mathbb{E}x_{k}^{(u)}
+\bar{D}_{k}\mathbb{E}u_{k}\Big\}\omega_{k},\\
x_{0}=\xi.
\end{array}\right.
\end{eqnarray*}
Using \eqref{ssss} again, we obtain
\begin{align*}
J^{0}(0,0;u)&=\mathbb{E}\Bigg\{\langle Gx_{N}^{(u)},x_{N}^{(u)}\rangle+\langle \bar{G}\mathbb{E}x_{N}^{(u)},\mathbb{E}x_{N}^{(u)}\rangle\nonumber\\
&~~~~~~~~~+\sum_{k=0}^{N-1}\left[\left\langle  \left(
                            \begin{array}{cc}
                              Q_{k} & S_{k}' \\
                              S_{k} & R_{k}\\
                            \end{array}
                          \right)\left(
                                   \begin{array}{c}
                                     x_{k}^{(u)} \\
                                     u_{k} \\
                                   \end{array}
                                 \right),\left(
                                   \begin{array}{c}
                                     x_{k}^{(u)} \\
                                     u_{k} \\
                                   \end{array}
                                 \right)\right\rangle
+\left\langle  \left(
                            \begin{array}{cc}
                              \bar{Q}_{k} & \bar{S}_{k}' \\
                              \bar{S}_{k} & \bar{R}_{k} \\
                            \end{array}
                          \right)\left(
                                   \begin{array}{c}
                                     \mathbb{E}x_{k}^{(u)} \\
                                     \mathbb{E}u_{k} \\
                                   \end{array}
                                 \right),\left(
                                   \begin{array}{c}
                                     \mathbb{E}x_{k}^{(u)} \\
                                     \mathbb{E}u_{k} \\
                                   \end{array}
                                 \right)\right\rangle\right]
\Bigg\}\\
&=\mathbb{E}\sum\limits_{k=0}^{N-1}\Big\{\langle \Theta_{k}'(R_{k}+B_{k}'P_{k+1}B_{k}+D_{k}'P_{k+1}D_{k})\Theta_{k}(x_{k}^{(u)}-\mathbb{E}x_{k}^{(u)}),x_{k}^{(u)}-\mathbb{E}x_{k}^{(u)}\rangle+\langle (B_{k}'P_{k+1}B_{k}\\
&~~~~~~~~~+D_{k}'P_{k+1}D_{k}+R_{k})u_{k},u_{k}\rangle-2\langle (R_{k}+B_{k}'P_{k+1}B_{k}+D_{k}'P_{k+1}D_{k})\Theta_{k}(x_{k}^{(u)}-\mathbb{E}x_{k}^{(u)}),u_{k}\rangle\\
&~~~~~~~~~+\langle \bar{\Theta}_{k}'[R_{k}+\bar{R}_{k}+(B_{k}+\bar{B}_{k})'\Pi_{k+1}(B_{k}+\bar{B}_{k})+(D_{k}+\bar{D}_{k})'P_{k+1}(D_{k}+\bar{D}_{k})]
\bar{\Theta}_{k}\mathbb{E}x_{k}^{(u)},\mathbb{E}x_{k}^{(u)}\rangle\\
&~~~~~~~~~-2\langle [R_{k}+\bar{R}_{k}+(B_{k}+\bar{B}_{k})'\Pi_{k+1}(B_{k}+\bar{B}_{k})+(D_{k}+\bar{D}_{k})'P_{k+1}(D_{k}+\bar{D}_{k})]\bar{\Theta}_{k}\mathbb{E}x_{k}^{(u)},\mathbb{E}u_{k}\rangle\\
&~~~~~~~~~+\langle[R_{k}+\bar{R}_{k}+(B_{k}+\bar{B}_{k})'\Pi_{k+1}(B_{k}+\bar{B}_{k})
+(D_{k}+\bar{D}_{k})'P_{k+1}(D_{k}+\bar{D}_{k})]\mathbb{E}u_{k},\mathbb{E}u_{k}\rangle\Big\}\\
&=\mathbb{E}\sum\limits_{k=0}^{N-1}\Big\{\langle[R_{k}+\bar{R}_{k}+(B_{k}+\bar{B}_{k})'\Pi_{k+1}(B_{k}+\bar{B}_{k})
+(D_{k}+\bar{D}_{k})'P_{k+1}(D_{k}+\bar{D}_{k})]\mathbb{E}u_{k},\mathbb{E}u_{k}\rangle\\
&~~~~~~~~~+\langle (R_{k}+B_{k}'P_{k+1}B_{k}+D_{k}'P_{k+1}D_{k})(u_{k}-\Theta_{k}x_{k}^{(u)}),u_{k}-\Theta_{k}x_{k}^{(u)}\rangle\Big\}.
\end{align*}
Combining \eqref{31} with Lemma 2.3, for some $\gamma>0$, we derive
\begin{align*}
J^{0}(0,0;u)\geq \alpha \gamma\mathbb{E}\sum\limits_{k=0}^{N-1}|u_{k}|^{2},~~~~\forall~u\in\mathscr{U}_{ad},
\end{align*}
then (i) follows. The proof is finished.
\end{proof}
\noindent
\textbf{Remark 4.2.} From the proof of Theorem 4.5, we deduce that the strongly regular solution of GRE (4.7) satisfies (4.11) with the same constant $\alpha>0$ as (i) of Theorem 4.5.

\vskip3mm\noindent
By virtue of Theorems 4.4 and 4.5, we would like to mention the following corollary.

\vskip3mm\noindent
\textbf{Corollary 4.2.} Assume that ($\mathbf{A2}$) holds, then Problem (MF-LQ) is uniquely open-loop solvable for all $(l,\xi)\in \mathbb{N}_{0}\times \mathbb{R}^{n}$, and the open-loop optimal control $u^{\ast}_{k}$ is given as
\begin{align*}
u^{\ast}_{k}&=-(R_{k}+B'_{k}P_{k+1}B_{k}+D'_{k}P_{k+1}D_{k})^{\dagger}(B_{k}'P_{k+1}A_{k}+D_{k}'P_{k+1}C_{k}+S_{k}')(x_{k}^{\ast}-\mathbb{E}x_{k}^{\ast})\\
&~~~~~~~~-[R_{k}+\bar{R}_{k}+(B_{k}+\bar{B}_{k})'\Pi_{k+1}(B_{k}+\bar{B}_{k})+(D_{k}+\bar{D}_{k})'P_{k+1}(D_{k}+\bar{D}_{k})]^{\dagger}\\
&~~~~~~~~\times[(B_{k}+\bar{B}_{k})'\Pi_{k+1}(A_{k}+\bar{A}_{k})+(D_{k}+\bar{D})'P_{k+1}(C_{k}+\bar{C})+S_{k}'+\bar{S}_{k}']\mathbb{E}x_{k}^{\ast}\\
&~~~~~~~~-[R_{k}+\bar{R}_{k}+(B_{k}+\bar{B}_{k})'\Pi_{k+1}(B_{k}+\bar{B}_{k})+(D_{k}+\bar{D}_{k})'P_{k+1}(D_{k}+\bar{D}_{k})]^{\dagger}\\
&~~~~~~~~\times[(D_{k}+\bar{D}_{k})'P_{k+1}\sigma_{k}+(B_{k}+\bar{B}_{k})'(\Pi_{k+1}b_{k}+\eta_{k+1})+\rho_{k}+\bar{\rho}_{k}],
\end{align*}
where $(P_{k},\Pi_{k})$ is the unique strongly regular solution of the GRE (4.7) and $\eta_{k}$ is the solution of the LRE \eqref{16}, besides, $x_{k}^{\ast}$ satisfies the following closed-loop system
\begin{eqnarray*}
  \left\{
  \begin{array}{ll}
x_{k+1}^{\ast}=A_{k}x_{k}^{\ast}+\bar{A}_{k}\mathbb{E}x_{k}^{\ast}+B_{k}u_{k}^{\ast}+\bar{B}_{k}\mathbb{E}u_{k}^{\ast}+b_{k}+(C_{k}x_{k}^{\ast}+\bar{C}_{k}\mathbb{E}x_{k}^{\ast}+D_{k}u_{k}^{\ast}
+\bar{D}_{k}\mathbb{E}u_{k}^{\ast}+\sigma_{k})\omega_{k},\\
x_{l}^{\ast}=\xi.
\end{array}\right.
\end{eqnarray*}
Specially, if $b, \sigma, g, q, \rho=0$, the adapted solution $(\eta,\zeta)\equiv(0,0)$ and unique optimal control becomes
\begin{align*}
u^{\ast}_{k}&=-(R_{k}+B'_{k}P_{k+1}B_{k}+D'_{k}P_{k+1}D_{k})^{\dagger}(B_{k}'P_{k+1}A_{k}+D_{k}'P_{k+1}C_{k}+S_{k}')(x_{k}^{\ast}-\mathbb{E}x_{k}^{\ast})\\
&~~~~~~~~-[R_{k}+\bar{R}_{k}+(B_{k}+\bar{B}_{k})'\Pi_{k+1}(B_{k}+\bar{B}_{k})+(D_{k}+\bar{D}_{k})'P_{k+1}(D_{k}+\bar{D}_{k})]^{\dagger}\\
&~~~~~~~~\times[(B_{k}+\bar{B}_{k})'\Pi_{k+1}(A_{k}+\bar{A}_{k})+(D_{k}+\bar{D})'P_{k+1}(C_{k}+\bar{C})+S_{k}'+\bar{S}_{k}']\mathbb{E}x_{k}^{\ast}.
\end{align*}





\setcounter{equation}{0}
\section{Finiteness of Problem (MF-LQ) and convexity of cost functional}
In previous Sections, it is shown that the uniform convexity of the cost functional indicates the open-loop and closed-loop solvabilities of Problem (MF-LQ). It is expected that the convexity of the cost functional is strongly associated with the finiteness of Problem (MF-LQ), which is the main goal of this section. To begin with, let
\begin{align*}
\Gamma=\left\{(P_{k},\Pi_{k})=(P_{k}',\Pi_{k}')\Bigg|\left(
                  \begin{array}{cc}
                    W_{k} & H_{k}' \\
                    H_{k} & \Upsilon_{k} \\
                  \end{array}
                \right)\geq0,~\left(
                  \begin{array}{cc}
                    \bar{W}_{k} & \bar{H}_{k}' \\
                    \bar{H}_{k} & \bar{\Upsilon}_{k} \\
                  \end{array}
                \right)\geq0,~P_{N}\leq G,~\Pi_{N}\leq G+\bar{G}\right\},
\end{align*}
where $\Upsilon_{k}, \bar{\Upsilon}_{k}$ are given as (4.8), and
\begin{align*}
W_{k}&=Q_{k}+C_{k}'P_{k+1}C_{k}+A_{k}'\Pi_{k+1} A_{k}-P_{k},~~H_{k}=B_{k}'P_{k+1}A_{k}+D_{k}'P_{k+1}C_{k}+S_{k}',\\
\bar{W}_{k}&=Q_{k}+\bar{Q}_{k}+(C_{k}+\bar{C}_{k})'P_{k+1}(C_{k}+\bar{C}_{k})+(A_{k}+\bar{A}_{k})'\Pi_{k+1}(A_{k}+\bar{A}_{k})-\Pi_{k},\\
\bar{H}_{k}&=(B_{k}+\bar{B}_{k})'\Pi_{k+1}(A_{k}+\bar{A}_{k})+(D_{k}+\bar{D}_{k})'P_{k+1}(C_{k}+\bar{C}_{k})+S_{k}'+\bar{S}_{k}'.
\end{align*}
By Corollary 3.1, we conclude that the map $u\longmapsto J(l,\xi;u)$ is convex if and only if
\begin{align}\label{50}
\mathcal{M}\geq0,
\end{align}
which is also equivalent to ($\mathbf{A1}$). Moreover, the map $u\longmapsto J(l,\xi;u)$ is uniformly convex if and only if for some $\alpha>0$,
\begin{align}\label{51}
\mathcal{M}\geq \alpha I,
\end{align}
which is also equivalent to ($\mathbf{A2}$).

\vskip2mm\noindent
\textbf{Lemma 5.1.} Among the following statements,

\noindent
(i) Problem (MF-LQ) is finite at $l$,

\noindent
(ii) Problem (MF-LQ)$^{0}$ is finite at $l$,

\noindent
(iii) there are $P_{l}$ and $\Pi_{l}$ such that
\begin{align}\label{52}
V^{0}(l,\xi)=\langle P_{l}(\xi-\mathbb{E}\xi),\xi-\mathbb{E}\xi\rangle+\langle \Pi_{l}\mathbb{E}\xi,\mathbb{E}\xi\rangle,
\end{align}
(iv) the map $u\mapsto J(l,\xi;u)$ is convex for any $\xi\in \mathbb{R}^{n}$,

\noindent
(v) $$J^{0}(l,0;u)\geq0,~~\forall~u\in \mathscr{U}_{ad},$$

\noindent
(vi) $\Gamma\neq\emptyset$,

\noindent
the following conclusions hold:
$${\rm{(i)}}\Rightarrow {\rm{(ii)}}\Rightarrow{\rm{(iii)}}\Rightarrow{\rm{(iv)}};~~{\rm{(i)}}\Rightarrow {\rm{(v)}};~~{\rm{(vi)}}\Rightarrow {\rm{(ii)}}.$$
\begin{proof}
(i) $\Rightarrow$ (ii) It is clear.

{\rm{(ii)}} $\Rightarrow$ {\rm{(iii)}} This part of the proof could see \cite{SL21}.

{\rm{(iii)}} $\Rightarrow$ {\rm{(iv)}} Following Corollary 3.1, if $u\mapsto J(l,\xi;u)$ is not convex, then for some $u\in \mathscr{U}_{ad}$, $J^{0}(l,0;u)<0$. By Lemma 3.1, it can be seen that for any $\lambda\in \mathbb{R}^{n}$,
\begin{align*}
J(l,\xi;u_{k}+\lambda v_{k})&=J(l,\xi;u_{k})+\lambda^{2}J^{0}(l,0;v_{k})+\lambda dJ(l,\xi;u_{k};v_{k}),
\end{align*}
letting $\lambda\rightarrow\infty$, we get $V^{0}(l,\xi)\leq \lim_{\lambda\rightarrow\infty}J^{0}(l,\xi;\lambda u)=-\infty,$
which is contradiction.

(i) $\Rightarrow$ (v) Assume that $J^{0}(l,0;u)<0$ for some $u\in \mathscr{U}_{ad}$. Using Lemma 3.1 yields
\begin{align*}
J(l,\xi;\lambda u)=J(l,\xi;0)+\lambda^{2}J^{0}(l,0;u)+\lambda dJ(l,\xi;0;u),~~~~\forall~\lambda\in \mathbb{R}.
\end{align*}
Letting $\lambda\rightarrow\infty$, we have
$V(l,\xi)\leq \lim_{\lambda\rightarrow\infty}J(l,\xi;\lambda u)=-\infty,$ which is a contradiction.

{\rm{(vi)}} $\Rightarrow$ {\rm{(ii)}} If $\Gamma\neq\emptyset$, simple calculations show that
\begin{align*}
J^{0}(l,\xi;u)&=\mathbb{E}\Big\{\langle(G-P_{N})(x_{N}-\mathbb{E}x_{N}),x_{N}-\mathbb{E}x_{N}\rangle+\langle(G+\bar{G}-\Pi_{N})\mathbb{E}x_{N},\mathbb{E}x_{N}\rangle\Big\}\\
&~~~~~~+\mathbb{E}\sum\limits_{k=l}^{N-1}\left\langle \left(
                  \begin{array}{cc}
                    W_{k} & H_{k}' \\
                    H_{k} & \Upsilon_{k} \\
                  \end{array}
                \right)\left(
                         \begin{array}{c}
                           x_{k}-\mathbb{E}x_{k} \\
                           u_{k}-\mathbb{E}u_{k} \\
                         \end{array}
                       \right),\left(
                                 \begin{array}{c}
                                   x_{k}-\mathbb{E}x_{k} \\
                                   u_{k}-\mathbb{E}u_{k} \\
                                 \end{array}
                               \right)
                \right\rangle+\mathbb{E}\langle P_{l}(\xi-\mathbb{E}\xi),\xi-\mathbb{E}\xi\rangle\\
&~~~~~~~+\mathbb{E}\sum\limits_{k=l}^{N-1}\left\langle \left(
                  \begin{array}{cc}
                    \bar{W}_{k} & \bar{H}_{k}' \\
                    \bar{H}_{k} & \bar{\Upsilon}_{k} \\
                  \end{array}
                \right)\left(
                         \begin{array}{c}
                           \mathbb{E}x_{k} \\
                            \mathbb{E}u_{k} \\
                         \end{array}
                       \right),\left(
                                 \begin{array}{c}
                                    \mathbb{E}x_{k} \\
                                    \mathbb{E}u_{k} \\
                                 \end{array}
                               \right)
                \right\rangle+\langle \Pi_{l}\mathbb{E}\xi,\mathbb{E}\xi\rangle\\
&\geq \langle P_{l}(\xi-\mathbb{E}\xi),\xi-\mathbb{E}\xi\rangle+\langle \Pi_{l}\mathbb{E}\xi,\mathbb{E}\xi\rangle>-\infty,
\end{align*}
which indicates that the corresponding Problem (MF-LQ)$^{0}$ is finite at $l$.
\end{proof}

\vskip0mm\noindent
\textbf{Remark 5.1.} Lemma 5.1 shows that $J^{0}(l,0;u)\geq0,~\forall~u\in \mathscr{U}_{ad}$, is the necessary condition for open-loop solvability and finiteness of Problem (MF-LQ) at $l$.

\vskip3mm\indent
Now, we give an example as follows, which shows that the convexity of the map $u\mapsto J^{0}(l,0;u)$ is not sufficient for the finiteness of Problem (MF-LQ)$^{0}$.

\vskip3mm\noindent
\textbf{Example 5.1.} Consider the following controlled stochastic dynamic system
\begin{eqnarray*}
  \left\{
  \begin{array}{ll}
x_{k+1}=u_{k}+x_{k}\omega_{k},\\
x_{l}=\xi,~~k\in\{l,\ldots,2\},~~l\in\{0,1,2\},
\end{array}\right.
\end{eqnarray*}
and cost functional
\begin{align*}
J^{0}(l,\xi;u)=-\mathbb{E}x_{3}^{2}+\mathbb{E}\sum\limits_{k=l}^{2}u_{k}^{2}.
\end{align*}
Let $\xi=0$, one has
\begin{align*}
J^{0}(0,0;u)=-\mathbb{E}x_{3}^{2}+\mathbb{E}\sum\limits_{k=0}^{2}u_{k}^{2}=-\mathbb{E}(u_{0}^{2}+u_{1}^{2}+u_{2}^{2})+\mathbb{E}\sum\limits_{k=0}^{2}u_{k}^{2}\geq0,~~\forall~u\in \mathscr{U}_{ad},
\end{align*}
this implies that the map $u\mapsto J^{0}(0,\xi;u)$ is convex.
On the other hand, let $\xi\neq0$ and $u_{k}=\lambda\omega_{k}$, $\lambda\in \mathbb{R}$, then one obtains
\begin{align*}
J^{0}(0,0;u)=-\mathbb{E}x_{3}^{2}+\mathbb{E}\sum\limits_{k=0}^{2}u_{k}^{2}=-\mathbb{E}(3\lambda^{2}+2\lambda\xi+\xi^{2})+3\lambda^{2}=-2\lambda\xi-\xi^{2},~~\forall~u\in \mathscr{U}_{ad}.
\end{align*}
Taking $|\lambda|\rightarrow\infty$ and $\lambda\xi>0$, one finally gets $V^{0}(0,\xi)=-\infty$.

\vskip3mm\indent
By Example 5.1, to derive the finiteness of Problem (MF-LQ)$^{0}$, one needs some additional conditions except for the convexity of $u(\cdot)\mapsto J^{0}(l,\xi;u)$. To this end, we assume that $u(\cdot)\mapsto J^{0}(l,\xi;u)$ is convex, by Lemma 5.1, ($\mathbf{A1}$) holds. For any $\varepsilon>0$, we consider the state equation (1.1) ($b_{k},\sigma_{k}=0$) and the following cost functional
\begin{align*}
J^{(0,\varepsilon)}(l,\xi;u)&\doteq J^{0}(l,\xi;u)+\varepsilon\mathbb{E}\sum_{k=l}^{N-1}|u_{k}|^{2}\\
&=\mathbb{E}\Big\{\langle Gx_{N},x_{N}\rangle+\langle \bar{G}\mathbb{E}x_{N},\mathbb{E}x_{N}\rangle\Big\}+\mathbb{E}\sum_{k=l}^{N-1}\left\langle  \left(
                            \begin{array}{cc}
                              Q_{k} & S_{k}' \\
                              S_{k} & R_{k}+\varepsilon I\\
                            \end{array}
                          \right)\left(
                                   \begin{array}{c}
                                     x_{k} \\
                                     u_{k} \\
                                   \end{array}
                                 \right),\left(
                                   \begin{array}{c}
                                     x_{k} \\
                                     u_{k} \\
                                   \end{array}
                                 \right)\right\rangle\nonumber\\
&~~~~~~~~+\mathbb{E}\sum_{k=l}^{N-1}\left\langle  \left(
                            \begin{array}{cc}
                              \bar{Q}_{k} & \bar{S}_{k}' \\
                              \bar{S}_{k} & \bar{R}_{k}\\
                            \end{array}
                          \right)\left(
                                   \begin{array}{c}
                                     \mathbb{E}x_{k} \\
                                     \mathbb{E}u_{k} \\
                                   \end{array}
                                 \right),\left(
                                   \begin{array}{c}
                                     \mathbb{E}x_{k} \\
                                     \mathbb{E}u_{k} \\
                                   \end{array}
                                 \right)\right\rangle.
\end{align*}
Denote the corresponding optimal control problem and value function by Problem (MF-LQ)$^{(0,\varepsilon)}$ and $V^{(0,\varepsilon)}$, respectively, then we can immediately obtain
\begin{align*}
J^{(0,\varepsilon)}(0,0;u)=J^{0}(0,0;u)+\varepsilon\mathbb{E}\sum_{k=0}^{N-1}|u_{k}|^{2}\geq\varepsilon\mathbb{E}\sum_{k=0}^{N-1}|u_{k}|^{2},
\end{align*}
which means $u\mapsto J^{(0,\varepsilon)}(0,0;u)$ is uniformly convex. Using Theorem 4.5 and Remark 5.1, the following GRE
\begin{eqnarray}\label{53}
  \left\{
  \begin{array}{ll}
P_{k}^{\varepsilon}=Q_{k}+A_{k}'P_{k+1}^{\varepsilon}A_{k}+C_{k}'P_{k+1}^{\varepsilon}C_{k}-(B_{k}'P_{k+1}^{\varepsilon}A_{k}+D_{k}'P_{k+1}^{\varepsilon}C_{k}+S_{k}')'\\
~~~~~~~~~~~~\times(R_{k}+\varepsilon I+B_{k}'P_{k+1}^{\varepsilon}B_{k}+D_{k}'P_{k+1}^{\varepsilon}D_{k})^{-1}
(B_{k}'P_{k+1}^{\varepsilon}A_{k}+D_{k}'P_{k+1}^{\varepsilon}C_{k}+S_{k}'),\\
\Pi_{k}^{\varepsilon}=Q_{k}+\bar{Q}_{k}+(A_{k}+\bar{A}_{k})'P_{k+1}^{\varepsilon}(A_{k}+\bar{A}_{k})+(C_{k}+\bar{C}_{k})'P_{k+1}^{\varepsilon}(C_{k}+\bar{C}_{k})\\
~~~~~~~~~~~~-[(B_{k}+\bar{B}_{k})'\Pi_{k+1}^{\varepsilon}(A_{k}+\bar{A}_{k})+(D_{k}+\bar{D}_{k})'P_{k+1}^{\varepsilon}(C_{k}+\bar{C}_{k})+S_{k}'+\bar{S}_{k}']'\\
~~~~~~~~~~~~\times[R_{k}+\bar{R}_{k}+\varepsilon I+(B_{k}+\bar{B}_{k})'\Pi_{k+1}^{\varepsilon}(B_{k}+\bar{B}_{k})+(D_{k}+\bar{D}_{k})'P_{k+1}^{\varepsilon}(D_{k}+\bar{D}_{k})]^{-1}\\
~~~~~~~~~~~~\times[(B_{k}+\bar{B}_{k})'\Pi_{k+1}^{\varepsilon}(A_{k}+\bar{A}_{k})+(D_{k}+\bar{D}_{k})'P_{k+1}^{\varepsilon}(C_{k}+\bar{C}_{k})+S_{k}'+\bar{S}_{k}'],\\
P_{N}^{\varepsilon}=G,~~\Pi_{N}^{\varepsilon}=G+\bar{G},
\end{array}\right.
\end{eqnarray}
exists a unique strongly regular solution $(P_{k}^{\varepsilon},\Pi_{k}^{\varepsilon})$ such that
\begin{eqnarray}\label{54}
  \left.
  \begin{array}{ll}
R_{k}+\varepsilon I+B_{k}'P_{k+1}^{\varepsilon}B_{k}+D_{k}'P_{k+1}^{\varepsilon}D_{k}\geq \varepsilon I,\\
R_{k}+\bar{R}_{k}+\varepsilon I+(B_{k}+\bar{B}_{k})'\Pi_{k+1}^{\varepsilon}(B_{k}+\bar{B}_{k})+(D_{k}+\bar{D}_{k})'P_{k+1}^{\varepsilon}(D_{k}+\bar{D}_{k})\geq\varepsilon I.
\end{array}\right.
\end{eqnarray}
Next, we shall give the following result about the finiteness of Problem (MF-LQ)$^{0}$.

\vskip3mm\noindent
\textbf{Theorem 5.1.} Assume that ($\mathbf{A1}$) holds. For any $\varepsilon>0$, let $(P_{k}^{\varepsilon},\Pi_{k}^{\varepsilon})$ be the unique strongly regular solution of the GRE \eqref{53}. Then Problem (MF-LQ)$^{0}$ is finite if and only if $\{(P^{\varepsilon}_{0},\Pi_{0}^{\varepsilon})\}$ is bounded from below. In this case,
\begin{align}\label{55s}
\lim_{\varepsilon\rightarrow0} P^{\varepsilon}_{k}=P_{k},~~\lim_{\varepsilon\rightarrow0} \Pi^{\varepsilon}_{k}=\Pi_{k},
\end{align}
and \eqref{52} holds. Moreover,
\begin{eqnarray}\label{55}
  \left.
  \begin{array}{ll}
&R_{k}+B'_{k}P_{k+1}B_{k}+D'_{k}P_{k+1}D_{k}\geq0,\\
&R_{k}+\bar{R}_{k}+(B_{k}+\bar{B}_{k})'\Pi_{k+1}(B_{k}+\bar{B}_{k})+(D_{k}+\bar{D}_{k})'P_{k+1}(D_{k}+\bar{D}_{k})\geq0,
\end{array}\right.
\end{eqnarray}
and $P_{k}, \Pi_{k}$ are monotonically nonincreasing as well as bounded.
Specially, if Problem (MF-LQ)$^{0}$ is finite at $l=0$, then it is finite.
\begin{proof} Necessity. Assume that Problem (MF-LQ)$^{0}$ is finite, and \eqref{52} holds.
For any $\varepsilon_{2}>\varepsilon_{1}>0$, one has
$$J^{(0,\varepsilon_{2})}(l,\xi;u)\geq J^{(0,\varepsilon_{1})}(l,\xi;u)\geq J^{0}(l,\xi;u),~~\forall~(l,\xi)\in \mathbb{N}_{0}\times \mathbb{R}^{n},~\forall~u\in \mathscr{U}_{ad},$$
Therefore, for any $(l,\xi)\in \mathbb{N}_{0}\times \mathbb{R}^{n}$,
\begin{align*}
\langle P^{\varepsilon_{2}}_{l}(\xi-\mathbb{E}\xi),\xi-\mathbb{E}\xi\rangle+\langle \Pi^{\varepsilon_{2}}_{l}\mathbb{E}\xi,\mathbb{E}\xi\rangle&=V^{(0,\varepsilon_{2})}(l,\xi)=\inf_{u\in \mathscr{U}_{ad}} J^{(0,\varepsilon_{2})}(l,\xi;u)\geq
\inf_{u\in \mathscr{U}_{ad}} J^{(0,\varepsilon_{1})}(l,\xi;u)\\
&=V^{(0,\varepsilon_{1})}(l,\xi)=\langle P^{\varepsilon_{1}}_{l}(\xi-\mathbb{E}\xi),\xi-\mathbb{E}\xi\rangle+\langle \Pi^{\varepsilon_{1}}_{l}\mathbb{E}\xi,\mathbb{E}\xi\rangle\\
&\geq \inf_{u\in \mathscr{U}_{ad}} J^{0}(l,\xi;u)=V^{0}(l,\xi)
=\langle P_{l}(\xi-\mathbb{E}\xi),\xi-\mathbb{E}\xi\rangle+\langle \Pi_{l}\mathbb{E}\xi,\mathbb{E}\xi\rangle.
\end{align*}
Let $\mathbb{E}\xi=0$, we get that $\{P_{l}^{\varepsilon}\}_{\varepsilon>0}$ is a nondecreasing sequence with lower bound $P_{l}$.
Similarly, let $\xi=\mathbb{E}\xi\neq0$, we obtain that $\{\Pi_{l}^{\varepsilon}\}$ is a nondecreasing sequence with lower bound $\Pi_{l}$. We denote their limits by
$\hat{P}_{l}$ and $\hat{\Pi}_{l}$, respectively, and
\begin{align}\label{56}
\hat{P}_{l}=\lim_{\varepsilon\rightarrow0}P_{l}^{\varepsilon}\geq P_{l},~~~~~~
\hat{\Pi}_{l}=\lim_{\varepsilon\rightarrow0}\Pi_{l}^{\varepsilon}\geq \Pi_{l},~~~~~\forall~l\in \mathbb{N}_{0}.
\end{align}
On the other hand, for any $\delta>0$, there admits a $u^{\delta}\in \mathscr{U}_{ad}$ such that
\begin{align*}
V^{(0,\varepsilon)}(l,\xi)\leq J^{0}(l,\xi;u^{\delta})+\varepsilon \mathbb{E}\sum\limits_{k=l}^{N-1}|u^{\delta}_{k}|^{2}\leq V^{0}(l,\xi)+\delta+\varepsilon \mathbb{E}\sum\limits_{k=l}^{N-1}|u^{\delta}_{k}|^{2}.
\end{align*}
Taking $\varepsilon\rightarrow0$, we have
\begin{align*}
\langle \hat{P}_{l}(\xi-\mathbb{E}\xi),\xi-\mathbb{E}\xi\rangle+\langle \hat{\Pi}_{l}\mathbb{E}\xi,\mathbb{E}\xi\rangle&=\lim_{\varepsilon\rightarrow0}[\langle P^{\varepsilon}_{l}(\xi-\mathbb{E}\xi),\xi-\mathbb{E}\xi\rangle+\langle \Pi^{\varepsilon}_{l}\mathbb{E}\xi,\mathbb{E}\xi\rangle]\\
&=\lim_{\varepsilon\rightarrow0}V^{(0,\varepsilon)}(l,\xi)\leq V^{0}(l,\xi)+\delta\\
&=\langle P_{l}(\xi-\mathbb{E}\xi),\xi-\mathbb{E}\xi\rangle+\langle \Pi_{l}\mathbb{E}\xi,\mathbb{E}\xi\rangle+\delta,
\end{align*}
furthermore,
\begin{align*}
\langle \hat{P}_{l}(\xi-\mathbb{E}\xi),\xi-\mathbb{E}\xi\rangle+\langle \hat{\Pi}_{l}\mathbb{E}\xi,\mathbb{E}\xi\rangle\leq\langle P_{l}(\xi-\mathbb{E}\xi),\xi-\mathbb{E}\xi\rangle+\langle \Pi_{l}\mathbb{E}\xi,\mathbb{E}\xi\rangle+\delta.
\end{align*}
Likewise, let $\mathbb{E}\xi=0$, $P_{l}$ is bounded, and let $\xi-\mathbb{E}\xi\neq0$, $\Pi_{l}$ is bounded.
From the above arguments, we see that \eqref{55s} holds. Letting $\varepsilon\rightarrow0$ in \eqref{54}, \eqref{55} follows.

\noindent
Sufficiency. Assume that there is a $\beta\in \mathbb{R}$ such that
$$P_{0}^{\varepsilon}\geq\beta I,~~\Pi_{0}^{\varepsilon}\geq\beta I,~~\forall~\varepsilon>0,$$
then for any $u\in \mathscr{U}_{ad}$ and $\xi\in \mathbb{R}^{n}$, one obtains
\begin{align*}
J^{0}(0,\xi;u)+\varepsilon\mathbb{E}\sum\limits_{k=0}^{N-1}|u_{k}|^{2}\geq V^{(0,\varepsilon)}(0,\xi)=\langle P^{\varepsilon}_{0}(\xi-\mathbb{E}\xi),\xi-\mathbb{E}\xi\rangle+\mathbb{E}\langle \Pi^{\varepsilon}_{0}\mathbb{E}\xi,\mathbb{E}\xi\rangle\geq \beta|\xi|^{2}.
\end{align*}
Taking $\varepsilon\rightarrow0$, for any $u\in \mathscr{U}_{ad}$ and $\xi\in \mathbb{R}^{n}$, one has
$J^{0}(0,\xi;u)\geq\beta|\xi|^{2}.$
This indicates the finiteness of Problem (MF-LQ)$^{0}$ at $l=0$. Notice that for all $\xi\in \mathbb{R}^{n}$, $V^{0}(0,\xi)=\langle P_{0}(\xi-\mathbb{E}\xi),\xi-\mathbb{E}\xi\rangle+\langle \Pi_{0}\mathbb{E}\xi,\mathbb{E}\xi\rangle$,
then $P_{0}\leq P_{0}^{\varepsilon}$.
Finally, we can immediately get that Problem (MF-LQ)$^{0}$ is finite.
\end{proof}

\vskip0mm\noindent
\textbf{Corollary 5.1.} Assume that there admits a pair $(\Phi,\Psi)$ such that
\begin{align*}
&R_{k}+B_{k}'P_{k+1}B_{k}+D_{k}'P_{k+1}D_{k}\geq (B_{k}'P_{k+1}A_{k}+D_{k}'P_{k+1}C_{k}+S_{k}')\Phi^{-1}_{k}(A_{k}'P_{k+1}B_{k}+C_{k}'P_{k+1}D_{k}+S_{k}),\\
&R_{k}+\bar{R}_{k}+(B_{k}+\bar{B}_{k})'\Pi_{k+1}(B_{k}+\bar{B}_{k})+(D_{k}+\bar{D}_{k})'P_{k+1}(D_{k}+\bar{D}_{k})\geq\\ &~~~~~~~~~~~~~~~[(B_{k}+\bar{B}_{k})'\Pi_{k+1}(A_{k}+\bar{A}_{k})+(D_{k}+\bar{D}_{k})'P_{k+1}(C_{k}+\bar{C}_{k})+S_{k}'+\bar{S}_{k}']\\
&~~~~~~~~~~~~~~~\times\Psi^{-1}_{k}[(A_{k}+\bar{A}_{k})'\Pi_{k+1}(B_{k}+\bar{B}_{k})+(C_{k}+\bar{C}_{k})'P_{k+1}(D_{k}+\bar{D}_{k})+S_{k}+\bar{S}_{k}],
\end{align*}
where $(P,\Pi)$ is the solution of the following Lyapunov equation
\begin{eqnarray}\label{57}
  \left\{
  \begin{array}{ll}
\Phi_{k}+P_{k}=A_{k}'P_{k+1}A_{k}+C_{k}'P_{k+1}C_{k}+Q_{k},\\
\Psi_{k}+\Pi_{k}=(A_{k}+\bar{A}_{k})'\Pi_{k+1}(A_{k}+\bar{A}_{k})+(C_{k}+\bar{C}_{k})'P_{k+1}(C_{k}+\bar{C}_{k})+Q_{k}+\bar{Q}_{k},\\
P_{N}\leq G,~~\Pi_{N}\leq G+\bar{G}.
\end{array}\right.
\end{eqnarray}
Then Problem (MF-LQ)$^{0}$ is finite.
\begin{proof} Observe that $\Gamma\geq0$, by Theorem 2.2 in \cite{NLZ16a}, we see that Problem (MF-LQ)$^{0}$ is finite at $l=0$. Using Theorem 6.2, the desired result follows.
\end{proof}

Following Remark 4.2 and Corollary 4.2, we obtain that if the uniformly convex condition (5.2) satisfies, then Problem (MF-LQ)$^{0}$ is finite and \eqref{411} holds. Now we shall show the converse is also valid.

\vskip2mm\noindent
\textbf{Theorem 5.2.} Assume that Problem (MF-LQ)$^{0}$ is finite, and \eqref{52} holds. For some $\alpha>0$, if \eqref{411} holds, then $(P_{k},\Pi_{k})$ is the solution of GRE \eqref{47}. Therefore, the map $u\mapsto J^{0}(0,0;u)$ is uniformly convex.
\begin{proof} For any $\varepsilon>0$, suppose $(P^{\varepsilon}_{k},\Pi^{\varepsilon}_{k})$ is the unique strongly regular solution of \eqref{53}. Using Theorem 5.1, as $\varepsilon\rightarrow0$,
$P^{\varepsilon}_{k}\downarrow P_{k},~\Pi^{\varepsilon}_{k}\downarrow \Pi_{k}.$
For all $\varepsilon>0$, $P^{\varepsilon}_{k}$ and $\Pi^{\varepsilon}_{k}$ have upper bound, notice that $P_{k}$ and $\Pi_{k}$ are bounded, thus, $\{P^{\varepsilon}_{k}\}_{\varepsilon>0}$ and $\{\Pi^{\varepsilon}_{k}\}_{\varepsilon>0}$ are uniformly bounded.
Besides,
\begin{align*}
&R_{k}+B'_{k}P_{k+1}^{\varepsilon}B_{k}+D'_{k}P_{k+1}^{\varepsilon}D_{k}\geq R_{k}+B'_{k}P_{k+1}B_{k}+D'_{k}P_{k+1}D_{k}\geq \alpha I,\\
&R_{k}+\bar{R}_{k}+(B_{k}+\bar{B}_{k})'\Pi_{k+1}^{\varepsilon}(B_{k}+\bar{B}_{k})+(D_{k}+\bar{D}_{k})'P_{k+1}^{\varepsilon}(D_{k}+\bar{D}_{k})\\
&~~~~~~~~~~~~~~~~\geq R_{k}+\bar{R}_{k}+(B_{k}+\bar{B}_{k})'\Pi_{k+1}(B_{k}+\bar{B}_{k})+(D_{k}+\bar{D}_{k})'P_{k+1}(D_{k}+\bar{D}_{k})\geq \alpha I.
\end{align*}
By virtue of the dominated convergence theorem, one has
\begin{eqnarray}\label{59}
  \left\{
  \begin{array}{ll}
P_{k}^{\varepsilon}=Q_{k}+A_{k}'P_{k+1}^{\varepsilon}A_{k}+C_{k}'P_{k+1}^{\varepsilon}C_{k}-(B_{k}'P_{k+1}^{\varepsilon}A_{k}+D_{k}'P_{k+1}^{\varepsilon}C_{k}+S'_{k})'\\
~~~~~~~~~~~\times(R_{k}+\varepsilon I+B_{k}'P_{k+1}^{\varepsilon}B_{k}+D_{k}'P_{k+1}^{\varepsilon}D_{k})^{-1}
(B_{k}'P_{k+1}^{\varepsilon}A_{k}+D_{k}'P_{k+1}^{\varepsilon}C_{k}+S'_{k}),\\
\Pi_{k}^{\varepsilon}=Q_{k}+\bar{Q}_{k}+(A_{k}+\bar{A}_{k})'P_{k+1}^{\varepsilon}(A_{k}+\bar{A}_{k})+(C_{k}+\bar{C}_{k})'P_{k+1}^{\varepsilon}(C_{k}+\bar{C}_{k})\\
~~~~~~~~~~~-[(B_{k}+\bar{B}_{k})'\Pi_{k+1}^{\varepsilon}(A_{k}+\bar{A}_{k})+(D_{k}+\bar{D}_{k})'P_{k+1}^{\varepsilon}(C_{k}+\bar{C}_{k})+S'_{k}+\bar{S}'_{k}]'\\
~~~~~~~~~~~\times[R_{k}+\bar{R}_{k}+\varepsilon I+(B_{k}+\bar{B}_{k})'\Pi_{k+1}^{\varepsilon}(B_{k}+\bar{B}_{k})+(D_{k}+\bar{D}_{k})'P_{k+1}^{\varepsilon}(D_{k}+\bar{D}_{k})]^{-1}\\
~~~~~~~~~~~\times[(B_{k}+\bar{B}_{k})'\Pi_{k+1}^{\varepsilon}(A_{k}+\bar{A}_{k})+(D_{k}+\bar{D}_{k})'P_{k+1}^{\varepsilon}(C_{k}+\bar{C}_{k})+S'_{k}+\bar{S}'_{k}],
\end{array}\right.
\end{eqnarray}
converges to (4.7) in $L^{1}$ as $\varepsilon\rightarrow0$.
Thus, $P_{k}=\lim_{\varepsilon\rightarrow0}P^{\varepsilon}_{k},~\Pi_{k}=\lim_{\varepsilon\rightarrow0}\Pi^{\varepsilon}_{k}$, along with \eqref{411},
which indicates that $(P_{k},\Pi_{k})$ is a strongly regular solution of \eqref{47}. Using Theorem 4.5, we obtain that $u\mapsto J^{0}(0,0;u)$ is uniformly convex.
\end{proof}

\setcounter{equation}{0}
\section{Minimizing sequences and open-loop solvabilities}
In Section 4, we have presented Problem (MF-LQ) is uniquely open-loop solvable under the uniform convexity condition. Differently, without the uniform convexity condition, we shall discuss the open-loop solvability of Problem (MF-LQ) in this Section.

\vskip3mm\noindent
\textbf{Theorem 6.1.} Assume that Problem (MF-LQ) is finite. For any $\varepsilon>0$, if $(P^{\varepsilon}_{k},\Pi^{\varepsilon}_{k})$ is the unique strongly regular solution to the GRE \eqref{53}, and $\eta^{\varepsilon}_{k}$, $x^{\varepsilon}_{k}$ are the solutions to the following LRE
\begin{eqnarray*}
  \left\{
  \begin{array}{ll}
\eta_{k}^{\varepsilon}=(C_{k}+\bar{C}_{k})'P_{k+1}^{\varepsilon}\sigma_{k}+(A_{k}+\bar{A}_{k})'(\Pi_{k+1}^{\varepsilon}b_{k}+\eta_{k+1}^{\varepsilon})+q_{k}+\bar{q}_{k}\\
~~~~~~~~~~~~~~~+\bar{\Theta}_{k}^{\varepsilon}[(D_{k}+\bar{D}_{k})'P_{k+1}^{\varepsilon}\sigma_{k}+(B_{k}+\bar{B}_{k})'(\Pi_{k+1}^{\varepsilon}b_{k}+\eta_{k+1}^{\varepsilon})+\rho_{k}+\bar{\rho}_{k}],\\
\eta_{N}^{\varepsilon}=g+\bar{g},
\end{array}\right.
\end{eqnarray*}
and the closed-loop system, respectively,
\begin{eqnarray*}
  \left\{
  \begin{array}{ll}
x_{k+1}^{\varepsilon}=(A_{k}+B_{k}\Theta_{k}^{\varepsilon})x_{k}^{\varepsilon}+B_{k}v^{\varepsilon}_{k}+\bar{B}_{k}\mathbb{E}v^{\varepsilon}_{k}
+[\bar{A}_{k}+B_{k}\bar{\Theta}_{k}^{\varepsilon}+\bar{B}_{k}(\Theta_{k}^{\varepsilon}+\bar{\Theta}_{k}^{\varepsilon})]\mathbb{E}x_{k}^{\varepsilon}+b_{k}\\
~~~~~~~~~~~~~~~~+\Big\{(C_{k}+D_{k}\Theta_{k}^{\varepsilon})x_{k}^{\varepsilon}+C_{k}v^{\varepsilon}_{k}+\bar{C}_{k}\mathbb{E}v^{\varepsilon}_{k}
+[\bar{C}_{k}+D_{k}\bar{\Theta}_{k}^{\varepsilon}+\bar{C}_{k}(\Theta_{k}^{\varepsilon}+\bar{\Theta}_{k}^{\varepsilon})]\mathbb{E}x_{k}^{\varepsilon}\Big\}\omega_{k}+\sigma_{k},\\
x_{l}^{\varepsilon}=\xi,
\end{array}\right.
\end{eqnarray*}
where
\begin{eqnarray*}
  \left\{
  \begin{array}{ll}
\Theta_{k}^{\varepsilon}=(R_{k}+B_{k}'P_{k+1}^{\varepsilon}B_{k}+D_{k}'P_{k+1}^{\varepsilon}D_{k}+\varepsilon I)^{-1}(B_{k}'P_{k+1}^{\varepsilon}A_{k}+D_{k}'P_{k+1}^{\varepsilon}C_{k}+S'_{k}),\\
\bar{\Theta}_{k}^{\varepsilon}=
[R_{k}+\bar{R}_{k}+(B_{k}+\bar{B}_{k})'\Pi_{k+1}^{\varepsilon}(B_{k}+\bar{B}_{k})+(D_{k}+\bar{D}_{k})'P_{k+1}^{\varepsilon}(D_{k}+\bar{D}_{k})+\varepsilon I]^{-1}\\
~~~~~~~~~~~~~~~\times[(B_{k}+\bar{B}_{k})'\Pi_{k+1}^{\varepsilon}(A_{k}+\bar{A}_{k})+(D_{k}+\bar{D}_{k})'P_{k+1}^{\varepsilon}(C_{k}+\bar{C}_{k})+S'_{k}+\bar{S}'_{k}],\\
v_{k}^{\varepsilon}=-[R_{k}+\bar{R}_{k}+(B_{k}+\bar{B}_{k})'\Pi_{k+1}^{\varepsilon}(B_{k}+\bar{B}_{k})+(D_{k}+\bar{D}_{k})'P_{k+1}^{\varepsilon}(D_{k}+\bar{D}_{k})+\varepsilon I]^{-1}\\
~~~~~~~~~~~~~~~\times[(D_{k}+\bar{D}_{k})'P_{k+1}^{\varepsilon}\sigma_{k}+(B_{k}+\bar{B}_{k})'(\Pi_{k+1}^{\varepsilon}b_{k}+\eta_{k+1}^{\varepsilon})+\rho_{k}+\bar{\rho}_{k}].
\end{array}\right.
\end{eqnarray*}
Then
\begin{align}\label{61}
u^{\varepsilon}_{k}= \Theta_{k}^{\varepsilon}(x_{k}^{\varepsilon}-\mathbb{E}x_{k}^{\varepsilon})+\bar{\Theta}_{k}^{\varepsilon}\mathbb{E}x_{k}^{\varepsilon}+v_{k}^{\varepsilon},~~~~\varepsilon>0,
\end{align}
is a minimizing sequence of $u\mapsto J(l,\xi;u)$:
\begin{align}\label{62}
\lim_{\varepsilon\rightarrow0}J(l,\xi;u^{\varepsilon}_{k})=\inf_{u\in\mathscr{U}_{ad}}J(l,\xi;u_{k})=V(l,\xi).
\end{align}
\begin{proof} For all $\varepsilon>0$, consider the sate equation (1.1) and cost functional
\begin{align*}
J^{\varepsilon}(l,\xi;u_{k})= J(l,\xi;u_{k})+\varepsilon \mathbb{E}\sum\limits_{k=l}^{N-1}|u_{k}|^{2}.
\end{align*}
Denote the above problem by Problem (MF-LQ)$^{\varepsilon}$ and corresponding value functional $V^{\varepsilon}(l,\xi)$. By Corollary 4.2, $u^{\varepsilon}_{k}$ defined by \eqref{61} is the unique optimal control of Problem (MF-LQ)$^{\varepsilon}$ at $(l,\xi)$. Since
\begin{align*}
\lim_{\varepsilon\rightarrow0} \mathbb{E}\sum\limits_{k=l}^{N-1}\varepsilon|u_{k}|^{2}&=\lim_{\varepsilon\rightarrow0}J^{\varepsilon}(l,\xi;u^{\varepsilon}_{k})
-\lim_{\varepsilon\rightarrow0}J(l,\xi;u^{\varepsilon}_{k})\\
&=\lim_{\varepsilon\rightarrow0}V^{\varepsilon}(l,\xi)-\lim_{\varepsilon\rightarrow0}J(l,\xi;u^{\varepsilon}_{k})
\leq \lim_{\varepsilon\rightarrow0}V^{\varepsilon}(l,\xi)-\lim_{\varepsilon\rightarrow0}V(l,\xi)=0,
\end{align*}
hence,
\begin{align*}
\lim_{\varepsilon\rightarrow0}J(l,\xi;u^{\varepsilon}_{k})=\lim_{\varepsilon\rightarrow0}V^{\varepsilon}(l,\xi)
-\lim_{\varepsilon\rightarrow0}\mathbb{E}\sum\limits_{k=l}^{N-1}\varepsilon|u_{k}|^{2}=V(l,\xi).
\end{align*}
The proof is finished.
\end{proof}

\vskip0mm\noindent
\textbf{Theorem 6.2.} Assume that ($\mathbf{A1}$) holds, then the following statements hold:

\noindent
(i) $\lim_{\varepsilon\rightarrow0}V^{\varepsilon}(l,\xi)=V(l,\xi).$
Specially, Problem (MF-LQ) is finite at $(l,\xi)$ if and only if $\{V^{\varepsilon}(l,\xi)\}_{\varepsilon>0}$
\indent is bounded from below.

\noindent
(ii) The sequence $\{u^{\varepsilon}_{k}\}_{\varepsilon>0}$ is the minimizing sequence of $u\mapsto J(l,\xi;u)$.
\begin{proof} (i) For any $\varepsilon_{2}>\varepsilon_{1}>0$, one has
$$J^{\varepsilon_{2}}(l,\xi;u)\geq J^{\varepsilon_{1}}(l,\xi;u)\geq J(l,\xi;u),~~\forall~u\in \mathscr{U}_{ad},$$
which indicates that
$V^{\varepsilon_{2}}(l,\xi)\geq V^{\varepsilon_{1}}(l,\xi)\geq V(l,\xi).$
Hence
\begin{align}\label{63}
\bar{V}(l,\xi)\equiv \lim_{\varepsilon\rightarrow0}V^{\varepsilon}(l,\xi)\geq V(l,\xi).
\end{align}
Besides, for any $K, \beta>0$, there admits a $u^{\beta}\in \mathscr{U}_{ad}$ such that
$$V^{\varepsilon}(l,\xi)\leq J(l,\xi;u^{\beta}_{k})+\varepsilon \mathbb{E}\sum_{k=l}^{N-1}|u^{\beta}_{k}|^{2}\leq \max\{V(l,\xi),-K\}+\beta+\varepsilon \mathbb{E}\sum_{k=l}^{N-1}|u^{\beta}_{k}|^{2}.$$
Letting $\varepsilon\rightarrow0$, it yields that
$V(l,\xi)\leq \max\{V(l,\xi),-K\}+\beta.$
Furthermore, we know
\begin{align}\label{64}
\bar{V}(l,\xi)\leq V(l,\xi).
\end{align}
Combining \eqref{63} with \eqref{64}, (i) follows.

(ii) If $V(l,\xi)>-\infty$, by virtue of (i), one gets
\begin{align*}
\varepsilon \mathbb{E}\sum_{k=l}^{N-1}|u^{\varepsilon}_{k}|^{2}=J^{\varepsilon}(l,\xi;u^{\varepsilon}_{k})-J(l,\xi;u^{\varepsilon}_{k})=V^{\varepsilon}(l,\xi)-J(l,\xi;u^{\varepsilon}_{k})
\leq V^{\varepsilon}(l,\xi)-V(l,\xi)\rightarrow0,~\mbox{as}~\varepsilon\rightarrow0,
\end{align*}
therefore,
$$\lim_{\varepsilon\rightarrow0}J(l,\xi;u^{\varepsilon})
=\lim_{\varepsilon\rightarrow0}\Bigg(V^{\varepsilon}(l,\xi)-\varepsilon\mathbb{E}\sum_{k=l}^{N-1}|u^{\varepsilon}_{k}|^{2}\Bigg)=V(l,\xi).$$
If $V(l,\xi)=-\infty$, using (i) again, one has
\begin{align*}
J(l,\xi;u^{\varepsilon})\leq J^{\varepsilon}(l,\xi;u^{\varepsilon})=V^{\varepsilon}(l,\xi)\rightarrow -\infty,~\mbox{as}~\varepsilon\rightarrow0,
\end{align*}
in this case, \eqref{62} also holds. The proof is completed.
\end{proof}

By means of the minimizing sequence \eqref{61}, the open-loop solvability of Problem (MF-LQ) could be characterized as follows.

\vskip3mm\noindent
\textbf{Theorem 6.3.} Assume that ($\mathbf{A1}$) holds, then the following statements are equivalent:

\noindent
(i) Problem (MF-LQ) is open-loop solvable at $(l,\xi)$.

\noindent
(ii) The sequence $\{u^{\varepsilon}\}_{\varepsilon>0}$ exists a weakly convergent subsequence.

\noindent
(iii) The sequence $\{u^{\varepsilon}\}_{\varepsilon>0}$ exists a strongly convergent subsequence.

\noindent
In the above cases, the convergent subsequence is an open-loop optimal control of Problem (MF-LQ) at $(l,\xi)$.
\begin{proof} (iii)$\Rightarrow$ (ii) is clear.

(i)$\Rightarrow$ (ii) and (i)$\Rightarrow$ (iii). Let $v^{\ast}$ be an open-loop optimal control of Problem (MF-LQ) at $(l,\xi)$. By Corollary 4.2, for any $\varepsilon>0$, $u^{\varepsilon}_{k}$ defined by \eqref{61} is the unique optimal control of Problem (MF-LQ)$^{\varepsilon}$ at $(l,\xi)$ and
\begin{eqnarray*}
V^{\varepsilon}(l,\xi)=J^{\varepsilon}(l,\xi;u^{\varepsilon}_{k})\geq V(l,\xi)+\varepsilon \mathbb{E}\sum\limits_{k=l}^{N-1}|u^{\varepsilon}_{k}|^{2},~~
V^{\varepsilon}(l,\xi)\leq J^{\varepsilon}(l,\xi;v^{\ast}_{k})= V(l,\xi)+\varepsilon \mathbb{E}\sum\limits_{k=l}^{N-1}|v^{\ast}_{k}|^{2}.
\end{eqnarray*}
Furthermore,
\begin{align}\label{65}
\mathbb{E}\sum\limits_{k=l}^{N-1}|u^{\varepsilon}_{k}|^{2}\leq \frac{V^{\varepsilon}(l,\xi)-V(l,\xi)}{\varepsilon}\leq \mathbb{E}\sum\limits_{k=l}^{N-1}|v^{\ast}_{k}|^{2},~~~\forall~\varepsilon>0.
\end{align}
Therefore, $\{u^{\varepsilon}_{k}\}_{\varepsilon>0}$ is bounded in the Hilbert space $\mathscr{U}_{ad}$ and then exists a weakly convergent subsequence $\{u^{\varepsilon_{k}}_{k}\}_{k\geq1}$. Denote the weak limit by $u^{\ast}_{k}$. Notice that $u\mapsto J(l,\xi;u)$ is convex and continuous, then it is sequentially weakly lower semi-continuous. By \eqref{62}, letting $k\rightarrow\infty$ in the above equation, we conclude that
\begin{align*}
V(l,\xi)\leq J(l,\xi;u^{\ast}_{k})\leq \varliminf _{k\rightarrow\infty}J(l,\xi;u^{\varepsilon_{k}}_{k})\leq J(l,\xi;u_{k}),~~\forall~u_{k}\in \mathscr{U}_{ad}.
\end{align*}
Consequently, $u^{\ast}_{k}$ is also an open-loop optimal control of Problem (MF-LQ) at $(l,\xi)$. By replacing $v^\ast$ with $u^\ast$ in \eqref{65}, we get
\begin{align}\label{66}
\mathbb{E}\sum\limits_{k=l}^{N-1}|(u^{\varepsilon}_{k})^{\ast}|^{2}\leq \mathbb{E}\sum\limits_{k=l}^{N-1}|u^{\ast}_{k}|^{2},~~~\forall~\varepsilon>0.
\end{align}
Further,
\begin{align}\label{67}
\mathbb{E}\sum\limits_{k=l}^{N-1}|u^{\ast}_{k}|^{2}\leq \lim_{k\rightarrow\infty}\inf\mathbb{E}\sum\limits_{k=l}^{N-1}|u_{k}^{\varepsilon}|^{2}.
\end{align}
Together with \eqref{66}-\eqref{67}, we have
\begin{align*}
\mathbb{E}\sum\limits_{k=l}^{N-1}|u^{\ast}_{k}|^{2}=\lim_{k\rightarrow\infty}\mathbb{E}\sum\limits_{k=l}^{N-1}|u_{k}^{\varepsilon_{k}}|^{2},
\end{align*}
which implies that $\{u_{k}^{\varepsilon}\}_{k\geq1}$ is strongly convergent to $u^{\ast}$.

\noindent
(iii) $\Rightarrow$ (i) Let $\{u_{k}^{\varepsilon}\}_{k\geq1}$ be a weakly convergent subsequence of $\{u^{\varepsilon}_{k}\}_{\varepsilon>0}$ with its weak limit $u^\ast_{k}$,
then $\{u_{k}^{\varepsilon_{k}}\}_{\varepsilon>0}$ is bounded in the Hilbert space $\mathscr{U}_{ad}$. Since for any $u_{k}\in \mathscr{U}_{ad}$,
\begin{align*}
J(l,\xi;(u_{k}^{\varepsilon})^{\ast})\geq V(l,\xi)+\varepsilon_{k} \mathbb{E}\sum\limits_{k=l}^{N-1}|u^{\ast}_{k}|^{2}
=V^{\varepsilon_{k}}(l,\xi)\leq J(l,\xi;u_{k})+\varepsilon_{k} \mathbb{E}\sum\limits_{k=l}^{N-1}|u_{k}|^{2},
\end{align*}
and $u\mapsto J(l,\xi;u)$ is sequentially weakly lower semi-continuous, we can deduce
\begin{align*}
J(l,\xi;u^{\ast}_{k})\leq \varliminf _{k\rightarrow\infty}J(l,\xi;u_{k}^{\varepsilon_{k}})\leq J(l,\xi;u_{k}),~~\forall~u_{k}\in \mathscr{U}_{ad}.
\end{align*}
This implies that $u^\ast_{k}$ is an optimal control of Problem (MF-LQ) at $(l,\xi)$.
\end{proof}

\vskip0mm\noindent
\textbf{Remark 6.1.} By Theorem 6.3, we could conclude that the open-loop solvability of Problem (MF-LQ) at $(l,\xi)$ is also equivalent to the $L^{2}$-boundedness of $\{u^{\varepsilon}_{k}\}_{\varepsilon>0}$.

\vskip3mm\noindent
\textbf{Theorem 6.4.} Assume that the map $u\mapsto J^{0}(0,0;u)$ is convex, and
\begin{align}\label{68}
\sup_{\varepsilon>0}\sum\limits_{k=0}^{N-1}(|\Theta_{k}^{\varepsilon}|^{2}+|\bar{\Theta}_{k}^{\varepsilon}|^{2})<\infty,
\end{align}
then the GRE (4.7) is regularly solvable and Problem (MF-LQ)$^{0}$ is closed-loop solvable.
\begin{proof} For any $\varepsilon>0$ and $\xi\in \mathbb{R}^{n}$, let $x^{\varepsilon}$ be the solution of
\begin{eqnarray*}
  \left\{
  \begin{array}{ll}
x_{k+1}^{\varepsilon}=(A_{k}+B_{k}\Theta_{k}^{\ast})x_{k}^{\varepsilon}
+[\bar{A}_{k}+B_{k}\bar{\Theta}_{k}^{\ast}+\bar{B}_{k}(\Theta_{k}^{\ast}+\bar{\Theta}_{k}^{\ast})]\mathbb{E}x_{k}^{\varepsilon}\\
~~~~~~~~~~~~~~~~~+\Big\{(C_{k}+D_{k}\Theta_{k}^{\ast})x_{k}^{\varepsilon}
+[\bar{C}_{k}+D_{k}\bar{\Theta}_{k}^{\ast}+\bar{D}_{k}(\Theta_{k}^{\ast}+\bar{\Theta}_{k}^{\ast})]\mathbb{E}x_{k}^{\varepsilon}\Big\}\omega_{k},\\
x_{0}^{\varepsilon}=\xi,
\end{array}\right.
\end{eqnarray*}
then there exists a constant $K>0$ such that
\begin{align*}
\mathbb{E}|x_{k}^{\varepsilon}|^{2}&=\Big\{(A_{k}+B_{k}\Theta_{k}^{\varepsilon})^{2}+(C_{k}+D_{k}\Theta_{k}^{\varepsilon})^{2}\Big\}\mathbb{E}|x_{k-1}^{\varepsilon}|^{2}
+\Big\{[\bar{A}_{k}+B_{k}\bar{\Theta}_{k}^{\varepsilon}+\bar{B}_{k}(\Theta_{k}^{\varepsilon}+\bar{\Theta}_{k}^{\varepsilon})]^{2}\\
&~~~~~~~~~~~~+[\bar{C}_{k}+D_{k}\bar{\Theta}_{k}^{\varepsilon}+\bar{D}_{k}(\Theta_{k}^{\varepsilon}+\bar{\Theta}_{k}^{\varepsilon})]^{2}
+2(\bar{A}_{k}+B_{k}\Theta_{k}^{\varepsilon})[\bar{A}_{k}+B_{k}\bar{\Theta}_{k}^{\varepsilon}+\bar{B}_{k}(\Theta_{k}^{\varepsilon}+\bar{\Theta}_{k}^{\varepsilon})]\\
&~~~~~~~~~~~~+2(\bar{C}_{k}+D_{k}\Theta_{k}^{\varepsilon})[\bar{C}_{k}+D_{k}\bar{\Theta}_{k}^{\varepsilon}+\bar{D}_{k}(\Theta_{k}^{\varepsilon}+\bar{\Theta}_{k}^{\varepsilon})]\Big\}
|\mathbb{E}x_{k-1}^{\varepsilon}|^{2}\\
&\leq \xi^{2}K\sum_{j=0}^{k-1}(\Theta_{j}^{\varepsilon}+\bar{\Theta}_{j}^{\varepsilon})^{2},
\end{align*}
this combines with \eqref{68}, further indicates that $\{\Theta_{k}^{\varepsilon}x_{k}^{\varepsilon}\}_{\varepsilon>0}$ and $\{\bar{\Theta}_{k}^{\varepsilon}\mathbb{E}x_{k}^{\varepsilon}\}_{\varepsilon>0}$ are bounded. Thus, according to Theorems 5.1 and 6.3, Problem (MF-LQ)$^{0}$ is open-loop solvable at $l=0$ and finite. Let \eqref{52} hold, using Theorem 5.1 again, we have $\Upsilon_{k}\geq0$ and $\bar{\Upsilon}_{k}\geq0$.
Let $\{\Theta^{\varepsilon_{k}}_{k}\}$ and $\{\bar{\Theta}^{\varepsilon_{k}}_{k}\}$ be convergent subsequences of $\{\Theta^{\varepsilon}_{k}\}, \{\bar{\Theta}^{\varepsilon}_{k}\}$ with limits $\Theta_{k}, \bar{\Theta}_{k}$, respectively. Observe that
\begin{align*}
&R_{k}+B_{k}'P^{\varepsilon}_{k+1}B_{k}+D_{k}'P^{\varepsilon}_{k+1}D_{k}+\varepsilon I\rightarrow R_{k}+B_{k}'P^{\varepsilon}_{k+1}B_{k}+D_{k}'P^{\varepsilon}_{k+1}D_{k},~~~\mbox{as}~\varepsilon\rightarrow0,\\
&R_{k}+\bar{R}_{k}+(B_{k}+\bar{B}_{k})'\Pi^{\varepsilon}_{k+1}(B_{k}+\bar{B}_{k})+(D_{k}+\bar{D}_{k})'P^{\varepsilon}_{k+1}(D_{k}+\bar{D}_{k})+\varepsilon I\\
&~~~~~~~~~~~~~\rightarrow R_{k}+\bar{R}_{k}+(B_{k}+\bar{B}_{k})'\Pi_{k+1}(B_{k}+\bar{B}_{k})+(D_{k}+\bar{D}_{k})'P_{k+1}(D_{k}+\bar{D}_{k}),~~~\mbox{as}~\varepsilon\rightarrow0,
\end{align*}
and $\{R_{k}+B_{k}'P^{\varepsilon}_{k+1}B_{k}+D_{k}'P^{\varepsilon}_{k+1}D_{k}+\varepsilon I\}_{0<\varepsilon\leq1}$, $\{R_{k}+\bar{R}_{k}+(B_{k}+\bar{B}_{k})'\Pi^{\varepsilon}_{k+1}(B_{k}+\bar{B}_{k})+(D_{k}+\bar{D}_{k})'P^{\varepsilon}_{k+1}(D_{k}+\bar{D}_{k})+\varepsilon I\}_{0<\varepsilon\leq1}$ are uniformly bounded, we get
\begin{align*}
&B_{k}'P^{\varepsilon_{k}}_{k+1}A_{k}+D_{k}'P^{\varepsilon_{k}}_{k+1}C_{k}+S_{k}'\\
&~~~~~~=-(R_{k}+B_{k}'P^{\varepsilon_{k}}_{k+1}B_{k}+D_{k}'P^{\varepsilon_{k}}_{k+1}D_{k}+\varepsilon I)\Theta^{\varepsilon_{k}}_{k}\rightarrow-(R_{k}+B_{k}'P_{k+1}B_{k}+D_{k}'P_{k+1}D_{k})\Theta_{k},\\
&(B_{k}+\bar{B}_{k})'\Pi^{\varepsilon_{k}}_{k+1}(A_{k}+\bar{A}_{k})+(D_{k}+\bar{D}_{k})'P^{\varepsilon_{k}}_{k+1}(C_{k}+\bar{C}_{k})+S_{k}'+\bar{S}_{k}'\\
&~~~~~~=-[R_{k}+\bar{R}_{k}+(B_{k}+\bar{B}_{k})'
\Pi^{\varepsilon}_{k+1}(B_{k}+\bar{B}_{k})+(D_{k}+\bar{D}_{k})'P^{\varepsilon}_{k+1}(D_{k}+\bar{D}_{k})+\varepsilon I]
\bar{\Theta}^{\varepsilon_{k}}_{k}\\
&~~~~~~~~~~~~~~\rightarrow-[R_{k}+\bar{R}_{k}+(B_{k}+\bar{B}_{k})'\Pi_{k+1}(B_{k}+\bar{B}_{k})+(D_{k}+\bar{D}_{k})'P_{k+1}(D_{k}+\bar{D}_{k})]\bar{\Theta}_{k}.
\end{align*}
Besides, notice that
\begin{align*}
&B_{k}'P^{\varepsilon_{k}}_{k+1}A_{k}+D_{k}'P^{\varepsilon_{k}}_{k+1}C_{k}+S_{k}'\rightarrow B_{k}'P_{k+1}A_{k}+D_{k}'P_{k+1}C_{k}+S_{k}',\\
&(B_{k}+\bar{B}_{k})'\Pi^{\varepsilon_{k}}_{k+1}(A_{k}+\bar{A}_{k})+(D_{k}+\bar{D}_{k})'P^{\varepsilon_{k}}_{k+1}(C_{k}+\bar{C}_{k})+S_{k}'+\bar{S}_{k}'\\
&~~~~~~~~~~~~~~~~\rightarrow (B_{k}+\bar{B}_{k})'\Pi_{k+1}(A_{k}+\bar{A}_{k})+(D_{k}+\bar{D}_{k})'P_{k+1}(C_{k}+\bar{C}_{k})+S_{k}'+\bar{S}_{k}',
\end{align*}
thus we get
\begin{align*}
&-(R_{k}+B_{k}'P_{k+1}B_{k}+D_{k}'P_{k+1}D_{k})\Theta_{k}=B_{k}'P_{k+1}A_{k}+D_{k}'P_{k+1}C_{k}+S_{k}',\\
&-[R_{k}+\bar{R}_{k}+(B_{k}+\bar{B}_{k})'\Pi_{k+1}(B_{k}+\bar{B}_{k})+(D_{k}+\bar{D}_{k})'P_{k+1}(D_{k}+\bar{D}_{k})]\bar{\Theta}_{k}\\
&~~~~~~~~~~~~~~~~~~=(B_{k}+\bar{B}_{k})'\Pi_{k+1}(A_{k}+\bar{A}_{k})+(D_{k}+\bar{D}_{k})'P_{k+1}(C_{k}+\bar{C}_{k})+S_{k}'+\bar{S}_{k}'.
\end{align*}
These show
\begin{align*}
&\mathcal{R}(B_{k}'P_{k+1}A_{k}+D_{k}'P_{k+1}C_{k}+S_{k}')\subseteq \mathcal{R}(R_{k}+B_{k}'P_{k+1}B_{k}+D_{k}'P_{k+1}D_{k}),\\
&\mathcal{R}[(B_{k}+\bar{B}_{k})'\Pi_{k+1}(A_{k}+\bar{A}_{k})+(D_{k}+\bar{D}_{k})'P_{k+1}(C_{k}+\bar{C}_{k})+S_{k}'+\bar{S}_{k}']\\
&~~~~~~~~~~~~\subseteq\mathcal{R}[R_{k}+\bar{R}_{k}+(B_{k}+\bar{B}_{k})'\Pi_{k+1}(B_{k}+\bar{B}_{k})+(D_{k}+\bar{D}_{k})'P_{k+1}(D_{k}+\bar{D}_{k})].
\end{align*}
Because
\begin{align*}
&(R_{k}+B_{k}'P_{k+1}B_{k}+D_{k}'P_{k+1}D_{k})^{\dagger}(B_{k}'P_{k+1}A_{k}+D_{k}'P_{k+1}C_{k}+S_{k}')\\
&~~~~=-(R_{k}+B_{k}'P_{k+1}B_{k}+D_{k}'P_{k+1}D_{k})^{\dagger}
(R_{k}+B_{k}'P_{k+1}B_{k}+D_{k}'P_{k+1}D_{k})\Theta_{k},\\
&[R_{k}+\bar{R}_{k}+(B_{k}+\bar{B}_{k})'\Pi_{k+1}(B_{k}+\bar{B}_{k})+(D_{k}+\bar{D}_{k})'P_{k+1}(D_{k}+\bar{D}_{k})]^{\dagger}\\
&~~~~~~~~~~\times[(B_{k}+\bar{B}_{k})'\Pi_{k+1}(A_{k}+\bar{A}_{k})+(D_{k}+\bar{D}_{k})'P_{k+1}(C_{k}+\bar{C}_{k})+S_{k}'+\bar{S}_{k}']\\
&~~~~=-[R_{k}+\bar{R}_{k}+(B_{k}+\bar{B}_{k})'\Pi_{k+1}(B_{k}+\bar{B}_{k})+(D_{k}+\bar{D}_{k})'P_{k+1}(D_{k}+\bar{D}_{k})]^{\dagger}\\
&~~~~~~~~~~\times[R_{k}+\bar{R}_{k}+(B_{k}+\bar{B}_{k})'\Pi_{k+1}(B_{k}+\bar{B}_{k})+(D_{k}+\bar{D}_{k})'P_{k+1}(D_{k}+\bar{D}_{k})]\bar{\Theta}_{k},
\end{align*}
and $(R_{k}+B_{k}'P_{k+1}B_{k}+D_{k}'P_{k+1}D_{k})^{\dagger}(R_{k}+B_{k}'P_{k+1}B_{k}+D_{k}'P_{k+1}D_{k})$, $[R_{k}+\bar{R}_{k}+(B_{k}+\bar{B}_{k})'\Pi_{k+1}(B_{k}+\bar{B}_{k})+(D_{k}+\bar{D}_{k})'P_{k+1}(D_{k}+\bar{D}_{k})]^{\dagger}
[R_{k}+\bar{R}_{k}+(B_{k}+\bar{B}_{k})'\Pi_{k+1}(B_{k}+\bar{B}_{k})+(D_{k}+\bar{D}_{k})'P_{k+1}(D_{k}+\bar{D}_{k})]$
are orthogonal projections, we obtain
\begin{align*}
&(R_{k}+B_{k}'P_{k+1}B_{k}+D_{k}'P_{k+1}D_{k})^{\dagger}(B_{k}'P_{k+1}A_{k}+D_{k}'P_{k+1}C_{k}+S_{k}')\in L^{2}(\mathbb{N};\mathbb{R}^{m\times n}),\\
&[R_{k}+\bar{R}_{k}+(B_{k}+\bar{B}_{k})'\Pi_{k+1}(B_{k}+\bar{B}_{k})+(D_{k}+\bar{D}_{k})'P_{k+1}(D_{k}+\bar{D}_{k})]^{\dagger}\\
&~~~~~~~\times[(B_{k}+\bar{B}_{k})'\Pi_{k+1}(A_{k}+\bar{A}_{k})+(D_{k}+\bar{D}_{k})'P_{k+1}(C_{k}+\bar{C}_{k})+S_{k}'+\bar{S}_{k}']\in L^{2}(\mathbb{N};\mathbb{R}^{m\times n}).
\end{align*}
Finally, taking $k\rightarrow\infty$ and combining with the above arguments, we see that the GRE (4.7) is regularly solvable. The proof is finished.
\end{proof}

\setcounter{equation}{0}
\section{Examples}
In this section, we give two examples to elaborate some results. In the first example, the standard condition (1.5) does not satisfy, but the corresponding GRE is still solvable and cost functional is uniformly convex. This deduces that the uniform convexity condition is weaker than (1.5).

\vskip3mm\noindent
\textbf{Example 7.1.} Consider the following discrete-time stochastic dynamical system
\begin{eqnarray*}
  \left\{
  \begin{array}{ll}
x_{k+1}=\sqrt{2}u_{k}+(x_{k}+u_{k}+\mathbb{E}u_{k})\omega_{k},\\
x_{0}=\xi,
\end{array}\right.
\end{eqnarray*}
and the cost functional
\begin{align*}
J(0,\xi;u)=\mathbb{E}(4x_{2}^{2}-3\mathbb{E}x_{2}^{2})-\sum\limits_{k=0}^{1}(u_{k}^{2}+\mathbb{E}u_{k}^{2}).
\end{align*}
The corresponding GRE is obtained as
\begin{eqnarray*}
  \left\{
  \begin{array}{ll}
P_{k}=P_{k+1}-\frac{P_{k+1}^{2}}{3P_{k+1}-1},~~~P_{2}=4,\\
\Pi_{k}=P_{k+1}-\frac{2P_{k+1}^{2}}{\Pi_{k+1}+2P_{k+1}-1},~~~\Pi_{2}=1.
\end{array}\right.
\end{eqnarray*}
Since $R_{k}\ll0,~R_{k}+\bar{R}_{k}\ll0$, the standard condition (1.5) does not hold. In this case, one could prove that the above GRE is solvable at $\{0,1,2\}$ with
\begin{eqnarray*}
  \left\{
  \begin{array}{ll}
P_{2}=4,~~~P_{1}=\frac{28}{11},~~~P_{0}=\frac{1260}{803},\\
\Pi_{2}=1,~~\Pi_{1}=0,~~\Pi_{0}=\frac{23268}{15917}.
\end{array}\right.
\end{eqnarray*}
It is not difficult to see
\begin{eqnarray*}
  \left\{
  \begin{array}{ll}
R_{k}+B_{k}'P_{k+1}B_{k}+D_{k}'P_{k+1}D_{k}\geq1,\\
R_{k}+\bar{R}_{k}+(B_{k}+\bar{B}_{k})'\Pi_{k+1}(B_{k}+\bar{B}_{k})+(D_{k}+\bar{D}_{k})'P_{k+1}(D_{k}+\bar{D}_{k})\geq1,
\end{array}\right.
\end{eqnarray*}
then using Theorem 4.5, we deduce that the cost functional $J(0,\xi;u)$ is uniformly convex and this problem exists an unique optimal control
\begin{align*}
u_{k}^{\ast}=-\frac{P_{k+1}}{3P_{k+1}-1}(x_{k}^{\ast}-\mathbb{E}x_{k}^{\ast})-\frac{P_{k+1}}{\Pi_{k+1}+2P_{k+1}-1}\mathbb{E}x_{k}^{\ast},
\end{align*}
where $x_{k}^{\ast}$ satisfies the following closed-loop system
\begin{eqnarray*}
  \left\{
  \begin{array}{ll}
x_{k+1}^{\ast}=-\sqrt{2}\frac{P_{k+1}}{3P_{k+1}-1}x_{k}^{\ast}+\sqrt{2}\left(\frac{P_{k+1}}{3P_{k+1}-1}-\frac{P_{k+1}}{\Pi_{k+1}+2P_{k+1}-1}\right)\mathbb{E}x_{k}^{\ast}\\
~~~~~~~~~~~~~+\left[\left(1-\frac{P_{k+1}}{3P_{k+1}-1}\right)x_{k}^{\ast}+\left(\frac{P_{k+1}}{3P_{k+1}-1}
-\frac{2P_{k+1}}{\Pi_{k+1}+2P_{k+1}-1}\right)\mathbb{E}x_{k}^{\ast}\right]\omega_{k},\\
x_{0}^{\ast}=\xi.
\end{array}\right.
\end{eqnarray*}

The following example tells us that the existence of a open-loop solvability does not imply the existence of a closed-loop solvability, since the conditions (4.9)-(4.10) does not satisfy and the solution of the corresponding GRE is not regular.

\vskip3mm\noindent
\textbf{Example 7.2.} Consider the following discrete-time stochastic controlled system
\begin{eqnarray*}
  \left\{
  \begin{array}{ll}
x_{k+1}=x_{k}+\mathbb{E}x_{k}+u_{k}^{(1)}-u_{k}^{(2)}+\mathbb{E}u_{k}^{(1)}+\mathbb{E}u_{k}^{(2)}
+\Big(u_{k}^{(1)}+u_{k}^{(2)}-\mathbb{E}u_{k}^{(1)}-\mathbb{E}u_{k}^{(2)}\Big)\omega_{k},\\
x_{0}=\xi,
\end{array}\right.
\end{eqnarray*}
and the cost functional
\begin{align*}
J(0,\xi;u)=\mathbb{E}\sum\limits_{k=0}^{4}(x_{k}^{2}+2\mathbb{E}x_{k}^{2})+\mathbb{E}(x_{5}^{2}+2\mathbb{E}x_{5}^{2}),
\end{align*}
where $u_{k}=(u_{k}^{(1)},u_{k}^{(2)})'$. The corresponding GRE is obtained as
\begin{eqnarray}
  \left\{
  \begin{array}{ll}
P_{k}=1+P_{k+1}-P_{k+1}^{2}\left(
                             \begin{array}{cc}
                               1 & -1 \\
                             \end{array}
                           \right)\left(
                                    \begin{array}{cc}
                                      2P_{k+1} & 0 \\
                                      0 & 2P_{k+1} \\
                                  \end{array}
                                 \right)^{\dagger}\left(
                                                    \begin{array}{c}
                                                       1 \\
                                                       -1 \\
                                                    \end{array}
                                                   \right)=1,~~~P_{5}=1,\\
\Pi_{k}=3+\Pi_{k+1}-\Pi_{k+1}^{2}\left(
                             \begin{array}{cc}
                               1 & 0 \\
                             \end{array}
                           \right)\left(
                                    \begin{array}{cc}
                                     \Pi_{k+1} & 0 \\
                                      0 & 0 \\
                                    \end{array}
                                  \right)^{\dagger}\left(
                                                     \begin{array}{c}
                                                       1 \\
                                                       0 \\
                                                     \end{array}
                                                   \right)=3,~~~\Pi_{5}=3.
\end{array}\right.
\end{eqnarray}
Clearly, $(P_{k},\Pi_{k})=(1,3)$ is the unique solution of the GRE (7.1). Notice that
\begin{align*}
&\mathcal{R}(B_{k}'P_{k+1}A_{k}+D_{k}'P_{k+1}C_{k}+S_{k})=\mathcal{R}((1~~-1)')=\{(a~~-a)':a\in \mathbb{R}\},\\
&\mathcal{R}[R_{k}+B_{k}'P_{k+1}B_{k}+D_{k}'P_{k+1}D_{k}]=\mathcal{R}\left(\left(
                                                            \begin{array}{cc}
                                                              2 & 0 \\
                                                             0 & 2 \\
                                                           \end{array}
                                                          \right)\right)=\{(a~~b)':a,~b\in \mathbb{R}\},\\
&\mathcal{R}[(B_{k}+\bar{B}_{k})'\Pi_{k+1}(A_{k}+\bar{A}_{k})+(D_{k}+\bar{D}_{k})'P_{k+1}(C_{k}+\bar{C}_{k})+S_{k}+\bar{S}_{k}]\\
&~~~~~~~~~~~~~=\mathcal{R}((12~~0)')=\{(c~~0)':c\in \mathbb{R}\},\\
&\mathcal{R}[R_{k}+\bar{R}_{k}+(B_{k}+\bar{B}_{k})'\Pi_{k+1}(B_{k}+\bar{B}_{k})+(D_{k}+\bar{D}_{k})'P_{k+1}(D_{k}+\bar{D}_{k})]\\
&~~~~~~~~~~~~~=\mathcal{R}\left(\left(
                                                            \begin{array}{cc}
                                                              12 & 0 \\
                                                              0 & 0 \\
                                                           \end{array}
                                                          \right)\right)
=\{(c~~0)':c\in \mathbb{R}\},
\end{align*}
the range conditions (4.9)-(4.10) are not satisfied, and $(P_{k},\Pi_{k})$ is not regular. According to Theorem 4.4, we see that this problem is not closed-loop solvable.
In order to discuss the open-loop solvability of the above problem, for any $\varepsilon>0$, we consider the following cost functional
\begin{align*}
J^{\varepsilon}(0,\xi;u)=J(0,\xi;u)+\varepsilon\mathbb{E}\sum\limits_{k=0}^{4}|u_{k}|^{2}
=\mathbb{E}\sum\limits_{k=0}^{4}(x_{k}^{2}+2\mathbb{E}x_{k}^{2}+\varepsilon u_{k}^{2})+\mathbb{E}(x_{5}^{2}+3\mathbb{E}x_{5}^{2}).
\end{align*}
The corresponding GRE is given as
\begin{eqnarray*}
  \left\{
  \begin{array}{ll}
P_{k}^{\varepsilon}=1+P_{k+1}^{\varepsilon}-\frac{2(P_{k+1}^{\varepsilon})^{2}}{2P_{k+1}^{\varepsilon}+\varepsilon},~~~P_{5}^{\varepsilon}=1,\\
\Pi_{k}^{\varepsilon}=3+\Pi_{k+1}^{\varepsilon}-\frac{(\Pi_{k+1}^{\varepsilon})^{2}}{\Pi_{k+1}^{\varepsilon}+\varepsilon},~~~\Pi_{5}^{\varepsilon}=3.
\end{array}\right.
\end{eqnarray*}
By taking $\varepsilon\rightarrow0$, we get $\hat{P}_{0}=\lim_{\varepsilon\rightarrow0}P_{k}^{\varepsilon}=1,~\hat{\Pi}_{0}=\lim_{\varepsilon\rightarrow0}\Pi_{k}^{\varepsilon}=3.$
Using Theorem 6.2, the Problem is finite and the value function
\begin{align*}
V(0,\xi)=\mathbb{E}\langle (\xi-\mathbb{E}\xi),\xi-\mathbb{E}\xi\rangle+
\langle 3\mathbb{E}\xi,\mathbb{E}\xi\rangle=\mathbb{E}\xi^{2}+2(\mathbb{E}\xi)^{2}.
\end{align*}
Besides,
let
\begin{align*}
&\Theta_{k}^{\varepsilon}=-(R_{k}+B_{k}'P_{k+1}^{\varepsilon}B_{k}+D_{k}'P_{k+1}^{\varepsilon}D_{k}+\varepsilon I)^{-1}(B_{k}'P_{k+1}^{\varepsilon}A_{k}+D_{k}'P_{k+1}^{\varepsilon}C_{k}+S_{k}')=\left(
                                                                         \begin{array}{cc}
                                                                           \frac{1}{2P_{k+1}+\varepsilon} & -\frac{1}{2P+\varepsilon} \\
                                                                        \end{array}
                                                                       \right)',\\
&\bar{\Theta}_{k}^{\varepsilon}=-[R_{k}+\bar{R}_{k}+(B_{k}+\bar{B}_{k})'\Pi_{k+1}^{\varepsilon}(B_{k}+\bar{B}_{k})+(D_{k}+\bar{D}_{k})'P_{k+1}^{\varepsilon}(D_{k}+\bar{D}_{k})+\varepsilon I]^{-1}\\
&~~~~~~~~~~~~~~~~\times[(B_{k}+\bar{B}_{k})'\Pi_{k+1}^{\varepsilon}(A_{k}+\bar{A}_{k})+(D_{k}+\bar{D}_{k})'P_{k+1}^{\varepsilon}(C_{k}+\bar{C}_{k})+S_{k}'+\bar{S}_{k}']=\left(
                                                                         \begin{array}{cc}
                                                                           \frac{1}{\Pi_{k+1}+\varepsilon} & 0 \\
                                                                         \end{array}
                                                                      \right)',
\end{align*}
then the solution of
\begin{eqnarray*}
  \left\{
  \begin{array}{ll}
x_{k+1}^{\varepsilon}=(A_{k}+B_{k}\Theta_{k})x_{k}^{\varepsilon}+[\bar{A}_{k}+B_{k}\bar{\Theta}_{k}+\bar{B}_{k}(\Theta_{k}+\bar{\Theta}_{k})]\mathbb{E}x_{k}^{\varepsilon}+b_{k}\\
~~~~~~~~~~~~+\Big\{(C_{k}+D_{k}\Theta_{k})x_{k}^{\varepsilon}+[\bar{C}_{k}+D_{k}\bar{\Theta}_{k}+\bar{D}_{k}(\Theta_{k}+\bar{\Theta}_{k})]\mathbb{E}x_{k}^{\varepsilon}+\sigma_{k}\Big\}\omega_{k},\\
x_{0}^{\varepsilon}=\xi.
\end{array}\right.
\end{eqnarray*}
is obtained as
\begin{align*}
x_{k+1}^{\varepsilon}=\bigg(1+\frac{2}{2P_{k+1}+\varepsilon}\bigg)x_{k}^{\varepsilon}+\bigg(1+\frac{2}{\Pi_{k+1}+\varepsilon}\bigg)\mathbb{E}x_{k}^{\varepsilon},
\end{align*}
and as $\varepsilon\rightarrow0$,
\begin{align*}
u_{k}^{\varepsilon}=\left(
                        \begin{array}{c}
                          \frac{\xi-\mathbb{E}\xi}{2P_{k+1}+\varepsilon}\Big(1+\frac{2}{2P_{k+1}+\varepsilon}\Big)^{k}+\frac{\mathbb{E}\xi}{\Pi_{k+1}+\varepsilon}
                          \Big(2+\frac{2}{2P_{k+1}+\varepsilon}+\frac{2}{\Pi_{k+1}+\varepsilon}\Big)^{k} \\
                          -\frac{\xi-\mathbb{E}\xi}{2P_{k+1}+\varepsilon}\Big(1+\frac{2}{2P_{k+1}+\varepsilon}\Big)^{k} \\
                        \end{array}
                      \right)
\rightarrow \left(
              \begin{array}{c}
                2^{k-1}(\xi-\mathbb{E}\xi)+\frac{1}{3}(\frac{11}{3})^{k}\mathbb{E}\xi \\
                -2^{k-1}(\xi-\mathbb{E}\xi) \\
              \end{array}
            \right),
\end{align*}
therefore, according to Theorem 6.3, the original problem is open-loop solvable and an open-loop optimal control
\begin{align*}
u_{k}^{\ast}=\left(
                               \begin{array}{cc}
                                 2^{k-1}(\xi-\mathbb{E}\xi)+\frac{1}{3}(\frac{11}{3})^{k}\mathbb{E}\xi & -2^{k-1}(\xi-\mathbb{E}\xi) \\
                               \end{array}
                             \right)'.
\end{align*}

\setcounter{equation}{0}
\section{Perspectives and open problems}
In this paper, we have discussed the open-loop and closed-loop solvabilities for the mean-field stochastic LQ problems. As we expected, these two solvabilities are intrinsically different. A key point that make our method valid is the equivalence of the strongly regular solvability of the GRE and the uniform convexity of the cost functional.
In this section, we give a brief exposition on the prospects that are open to the researchers. The following topics shall be explored in our future publications.

(1) In this paper, we have discussed the stochastic LQ problems with deterministic coefficients, in which the fact $\mathbb{E}(A_{k}x_{k})=A_{k}\mathbb{E}(x_{k})$ has played a crucial role. Natural, we hope that a theory could also be established for problems with random coefficients. In that setting, such a fact would be no longer valid, and the GRE would turn to be a nonlinear backward stochastic difference equation whose solvability is very tricky and challenging.

(2) In general, the process $\{v_{k}=\sum_{l=0}^{k}\omega_{l}\}$ does not possess the so-called martingale characterization, as the one for the Brownian motion. In the future, it is expected to explore the infinite horizon discrete-time MF-LQ optimal control problems by infinite horizon MF-FBSDEs.

(3) We would like to mention the so-called time inconsistency in a near future.



\vskip8mm
\renewcommand\baselinestretch{0.9}


{\small\begin{thebibliography}{99}


\bibitem{ACZ02} {\sc M. Ait Rami, X. Chen and X.Y. Zhou}, {\em Discrete-time indefinite LQ control with state and control dependent noises}, J. Global Optim., 23 (2002), pp.~ 245--265.

\bibitem{BC20} {\sc F. Barbieri and O.L.V. Costa}, {\em Mean-field formulation for the infinite-horizon mean-variance control of discrete-time linear systems with multiplicative noises}, IET Control Theory Appl., 14(17) (2020), pp.~2600--2612.

\bibitem{BDT19} {\sc J. Barreiro-Gomez, T.E. Duncan and H. Tembine}, {\em Co-opetitive linear-quadratic mean-field-type games}, IEEE Trans. Cybern., 50(12) (2019), pp.~5089--5098.

\bibitem{BDT20} {\sc J. Barreiro-Gomez, T.E. Duncan and H. Tembine}, {\em Discrete-time linear-quadratic mean-field-type repeated games: Perfect, incomplete, and imperfect information}, Automatica, 112 (2020): 108647.

\bibitem{B23} {\sc N. B\"{a}uerle}, {\em Mean field Markov decision processes}, Appl. Math. \& Optim., 88(1) (2023): 12.

\bibitem{BFY13} {\sc A. Bensoussan, J. Frehse and P. Yam}, {\em Mean Field Games and Mean Field Type Control Theory}, Springer, New York, 2013.

\bibitem{BSYY16} {\sc A. Bensoussan, K.C.J. Sung, S.C.P. Yam and S.P. Yung}, {\em Linear-quadratic mean field games}, J. Optim. Theory Appl., 169 (2016), pp.~496--529.

\bibitem{CLZ98} {\sc S. Chen, X. Li and X.Y. Zhou}, {\em Stocahstic linear quadratic regulators with indefinite control weight costs}, SIAM J. Control Optim., 36 (1998), pp.~ 1685--1702, doi:10.1137/S0363012996310478.

\bibitem{CLZ98} {\sc J.Y. Di, C. Tan, Z.Q. Zhang and W.S. Wong}, {\em Stabilisation for discrete-time mean-field stochastic Markov jump systems with multiple delays}, IET Control Theory \& Appl., (2023), https://doi.org/10.1049/cth2.12477.

\bibitem{DNW22} {\sc B. Dong, T. Nie and Z. Wu}, {\em Maximum principle for discrete-time stochastic control problem of mean-field type}, Automatica, 144 (2022): 110497.

\bibitem{D15} {\sc K. Du}, {\em Solvability conditions for indefinite linear quadratic optimal stochastic control problems and associated stochastic Riccati equations}, SIAM J. Control Optim., 53 (2015), pp.~3673--3689.

\bibitem{ELN13} {\sc R. Elliott, X. Li and Y.H. Ni}, {\em Discrete time mean-field stochastic linear-quadratic optimal control problems}, Automatica, 49 (2013), pp.~3222--3233.

\bibitem{FZZ19} {\sc H. Frankowska, H.S. Zhang and X. Zhang}, {\em Necessary optimality conditions for local minimizers of stochastic optimal control problems with state constraints}, Trans. Amer. Math. Soc., 372(2) (2019), pp.~1289--1331.

\bibitem{GLZ23} {\sc W.H. Gao, Y.N. Lin, and W.H. Zhang}, {\em Incentive feedback Stackelberg strategy for the discrete-time stochastic systems}, J. Frankl. Inst., 360(3) (2023), pp.~2404--2420.

\bibitem{G16} {\sc P.J. Graber}, {\em Linear quadratic mean field type control and mean field games with common noise with application to production of an exhaustible resource}, Appl. Math. Opim., 74 (2016), pp.~459--486.

\bibitem{HLY15} {\sc J.H. Huang, X. Li and J.M. Yong}, {\em A linear-quadratic optimal control problem for mean-field stochastic differential equations in infinite horizon}, Math. Control Relat. Fields, 5 (2015), pp.~97--139.

\bibitem{K60} {\sc R.E. Kalman}, {\em Contribution to the theory of optimal control}, Bol. Soc. Mat. Mexicana, 5(2) (1960), pp.~102--119.

\bibitem{LLY20} {\sc N. Li, X. Li and Z.Y. Yu}, {\em Indefinite mean-field type linear-quadratic stochastic optimal control problems}, Automatica, 122 (2020): 109267.

\bibitem{L19a} {\sc Y.N. Lin}, {\em Feedback Stackelberg strategies for the discrete-time mean-field stochastic systems in infinite horizon}, J. Frankl. Inst., 356(10) (2019), pp.~5222--5239.

\bibitem{LJZ18} {\sc Y.N. Lin, X.S. Jiang and W.H. Zhang}, {\em An open-loop Stackelberg strategy for the linear quadratic mean-field stochastic differential game}, IEEE Trans. Autom. Control, 64(1) (2018), pp.~97--110.

\bibitem{LZ01} {\sc A.E.B. Lim and X.Y. Zhou}, {\em Linear-quadratic control of backward stochastic differential equations}, SIAM J. Control Optim., 40 (2001), pp.~450--474.

\bibitem{LLL18} {\sc R.R. Liu, Y. Li and X.K. Liu}, {\em Linear-quadratic optimal control for unknown mean-field stochastic discrete-time system via adaptive dynamic programming approach}, Neurocomputing, 282 (2018), pp.~16--24.

\bibitem{LLL20} {\sc X.K. Liu, Q.M. Liu and Y. Li}, {\em Finite-time guaranteed cost control for uncertain mean-field stochastic systems}, J. Frankl. Inst., 357(5) (2020), pp.~2813--2829.

\bibitem{LSX19} {\sc X. Li, J. Sun and J. Xiong}, {\em Linear quadratic optimal control problems for mean-field backward stochastic differential equations}, Appl. Math. Optim., 80 (2019), pp.~223--250.

\bibitem{LSY16} {\sc X. Li, J.R. Sun and J.M. Yong}, {\em Mean-field stochastic linear quadratic optimal control problems: closed-loop solvability}, Probab. Uncertain. Quant. Risk, 1(1) (2016), pp.~1--24.

\bibitem{L19b} {\sc Q. L\"{u}}, {\em Well-posedness of stochastic Riccati equations and closed-loop solvability for stochastic linear quadratic optimal control problems}, J. Differential Equations, 267 (2019), pp.~180--227.

\bibitem{MZZ19} {\sc L.M. Ma, W.H. Zhang and Y. Zhao}, {\em Study on stability and stabilizability of discrete-time mean-field stochastic systems}, J. Frankl. Inst., 356(4) (2019), pp.~2153--2171.

\bibitem{M20} {\sc J. Moon}, {\em Linear-quadratic mean field stochastic zero-sum differential games}, Automatica, 120 (2020): 109067.

\bibitem{NCZZ17} {\sc  Y.H. Ni, K.F. Cedric Yiu, H.S. Zhang and J.F. Zhang}, {\em Delayed optimal control of stochastic LQ problem}, SIAM J. Control Optim., 55(5) (2017), pp.~3370--3407.

\bibitem{NEL15} {\sc Y.H. Ni, R. Elliott and X. Li}, {\em Discrete time mean-field stochastic linear-quadratic optimal control problems}, II: Infinite horizon case, Automatica, 57 (2015), pp.~65--77.

\bibitem{NLZ16a} {\sc Y.H. Ni, X. Li and J.F. Zhang}, {\em Indefinite mean-field stochastic linear-quadratic optimal control: from finite horizon to infinite horizon}, IEEE Trans. Autom. Control, 61 (11) (2016), pp.~3269--3284.

\bibitem{NLZ16b} {\sc Y.H. Ni, X. Li and J.F. Zhang}, {\em Mean-field stochastic linear-quadratic optimal control with Markov jump parameters}, Systems Control Lett., 93 (2016), pp.~69--76.

\bibitem{NZL15} {\sc Y.H. Ni, J.F. Zhang and X. Li}, {\em Indefinite mean-field stochastic linear-quadratic optimal control}, IEEE Trans. Autom. Control, 60 (7) (2015), pp.~1786--1800.

\bibitem{S23} {\sc N. Saldi}, {\em Linear Mean-Field Games with Discounted Cost}, arXiv preprint arXiv:2301.06074, (2023).

\bibitem{SL21} {\sc T. Song and B. Liu}, {\em Discrete-time mean-field stochastic linear quadratic optimal control problem with finite horizon}, Asian J. Control, 23(2) (2021), pp.~979--989.

\bibitem{SL23} {\sc T. Song and B. Liu}, {\em First- and second-order maximum principles for discrete-time stochastic optimal control with recursive utilities}, IEEE Trans. Autom. Control, (2023), https://doi.org/10.1109/TAC.2023.3261345.

\bibitem{S17} {\sc J.R. Sun}, {\em Mean-field stochastic linear quadratic optimal control problems: Open-loop solvabilities}, ESAIM Control Optim. Calc. Var., 23(3) (2017), pp.~1099--1127.

\bibitem{S21} {\sc J.R. Sun}, {\em Two-person zero-sum stochastic linear-quadratic differential games}, SIAM J. Control Optim., 59 (3) (2021), pp.~1804--1829.

\bibitem{SLY16} {\sc J.R. Sun, X. Li and J.M. Yong}, {\em Open-loop and closed-loop solvabilities for stochastic linear quadratic optimal control problems}, SIAM J. Control Optim., 54 (5) (2016), pp.~2274--2308.

\bibitem{SXY21} {\sc J.R. Sun, J. Xiong and J.M. Yong}, {\em Indefinite stochastic linear-quadratic optimal control problems with random coefficients: closed-loop representation of open-loop optimal controls}, Ann. Appl. Probab., 31(1) (2021), pp.~460--499.

\bibitem{SY14} {\sc J.R. Sun and J.M. Yong}, {\em Linear quadratic stochastic differential games: open-loop and closed-loop saddle points}, SIAM J. Control Optim., 52(6) (2014), pp.~4082--4121.

\bibitem{SY20} {\sc J.R. Sun and J.M. Yong}, {\em Stochastic Linear-Quadratic Optimal Control Theory: Open-Loop and Closed-Loop Solutions}, Springer Nature, 2020.

\bibitem{SWW21} {\sc J.R. Sun, H.X. Wang and Z. Wu}, {\em Mean-field linear-quadratic stochastic differential games}, J. Differential Equations, 296 (2021), pp.~299--334.

\bibitem{T03} {\sc S.J. Tang}, {\em General linear quadratic optimal stochastic control problems with random coefficients: Linear stochastic Hamilton systems and backward stochastic Riccati equations}, SIAM J. Control Optim., 42 (2003), pp.~53--75.

\bibitem{Y13} {\sc J.M. Yong}, {\em Linear-quadratic optimal control problems for mean-field stochastic differential equations}, SIAM J. Control Optim., 51 (2013), pp.~2809--2838.

\bibitem{Y17} {\sc J. M. Yong}, {\em Linear-quadratic optimal control problems for mean-field stochastic differential equations-time-consistent solutions}, Trans. Amer. Math. Soc., 369 (2017), pp.~5467--5523.

\bibitem{WSY19} {\sc H.X. Wang, J.R. Sun and J.M. Yong}, {\em Weak closed-loop solvability of stochastic linear-quadratic optimal control problems}, Discrete Contin. Dyn. Syst. Ser. A, 39 (2019), pp.~2785--2805.

\bibitem{WZL18} {\sc T. Wang, H.G. Zhang and Y.H. Luo}, {\em Stochastic linear quadratic optimal control for model-free discrete-time systems based on Q-learning algorithm}, Neurocomputing, 312 (2018), pp.~1--8.

\bibitem{ZDS23} {\sc T.L. Zhang, F.Q. Deng and P. Shi}, {\em Non-fragile finite-time stabizazation for discrete mean-field stochastic systems}, IEEE Trans. Autom. Control, (2023), https://10.1109/TAC.2023.3238849.

\bibitem{ZQF18} {\sc H.S. Zhang, Q.Y. Qi and M.Y. Fu}, {\em Optimal stabilization control for discrete-time mean-field stochastic systems}, IEEE Trans. Autom. Control, 64(3) (2018), pp.~1125--1136.



\end{thebibliography}}
\end{document}